\documentclass[leqno]{article}
\usepackage{amssymb,latexsym,amscd}
\newcommand{\se}[1]{{\section{#1}} {\setcounter{equation}{0}}}
\newtheorem{theorem}{Theorem}[section]
\newtheorem{lm}{Lemma}[section]
\newtheorem{prop}{Proposition}[section]
\newtheorem{de}{Definition}[section]
\newtheorem{co}{Corollary}[section]

\def\k{{K\"{a}hler }}

\def\cy{{Calabi-Yau }}
\def\l{{Lagrangian }}

\input epsf
\begin{document}
\hbadness=10000
\title{{\bf Lagrangian torus fibration of quintic Calabi-Yau hypersurfaces II:}\\
{\Large {\bf Technical results on gradient flow construction}}}
\author{Wei-Dong Ruan\\
Department of Mathematics\\
University of Illinois at Chicago\\
Chicago, IL 60607\\}
\date{}
\footnotetext{Partially supported by NSF Grant DMS-9703870 and DMS-0104150.}
\maketitle
\tableofcontents
\se{Introduction and background}
This paper is a sequel to my recent paper \cite{lag1}. It will provide technical details of our gradient flow construction and related problems, which are essential for our construction of Lagrangian torus fibrations in \cite{lag1} and subsequent papers \cite{lag3,tor,ci}.\\

\subsection{Background}
The motivation of our work on Lagrangian torus fibrations of \cy manifolds comes from the Strominger-Yau-Zaslow conjecture of mirror symmetry. According to their conjecture, on each Calabi-Yau manifold there should exist a special Lagrangian torus fibration. This conjectural special Lagrangian torus fibration structure is used to give a construction of the mirror \cy manifold and even a possible explanation of mirror symmetry. Despite its great potential in solving the mirror symmetry conjecture, there are very few known examples of special Lagrangian submanifolds or special Lagrangian fibrations for dimension $n\geq 3$. Given our lack of knowledge for special Lagrangian, one may consider relaxing the requirement to Lagrangian fibrations, which is largely unexplored and interesting in its own right. Special Lagrangians are very rigid, on the other hand, Lagrangian submanifolds are more flexible and can be modified locally by Hamiltonian deformation. So it is a reasonable first step to take. For many applications to mirror symmetry, especially those concerning (symplectic) topological structure of fibrations, Lagrangian torus fibrations will provide quite sufficient information. In this paper, as in the previous paper \cite{lag1}, we will mainly concern Lagrangian torus fibrations of Calabi-Yau hypersurfaces in toric varieties. Aside from its application to mirror symmetry, from purely mathematical point of view, the construction of Lagrangian torus fibrations for Calabi-Yau manifolds is clearly important for understanding the topology and geometry of Calabi-Yau manifolds. It is also of independent interest in symplectic geometry.\\

In our paper \cite{lag1}, we described a very simple and natural construction of Lagrangian torus fibrations via gradient flow which in principle will be able to produce Lagrangian torus fibrations for general Calabi-Yau hypersurfaces in toric varieties. For simplicity, we described in great detail the case of Fermat type quintic Calabi-Yau threefold family $\{X_\psi\}$ in ${\mathbb{CP}^4}$ defined by\\
\[
p_{\psi}=\sum_1^5 z_k^5 - 5\psi \prod_{k=1}^5 z_k=0
\]\\
near the large complex limit $X_{\infty}$\\
\[
p_{\infty}=\prod_{k=1}^5 z_k=0.
\]\\
Most of the essential features for more general cases already showed up there. We also discussed the so-called {\it expected} special Lagrangian torus fibration structure, especially the monodromy transformations of the expected fibration and the expected singular fibre structures implied by monodromy information in this case. Then we compared the Lagrangian fibrations we constructed with the {\it expected} special Lagrangian torus fibrations. Finally, we discussed its relavence to mirror construction for Fermat type quintic Calabi-Yau hypersurfaces.\\

Due to its position in mathematics and its origin from physics, our work on Lagrangian torus fibrations is pursuing two (sometimes rather different) goals. From physics point of view, the construction of Lagrangian torus fibration is intended to be used to uncover the symplectic topological structure of the special Lagrangian torus fibration in SYZ conjecture. For such purpose,  the Lagrangian torus fibrations that reveal the structures of the conjectured special Lagrangian fibrations are prefered. From mathematical point of view, the construction of Lagrangian torus fibrations can be used to understand the symplectic topology of Calabi-Yau manifolds. For such purpose, it is prefered that the Lagrangian fibrations have simple and well hehaved singular locus and singular fibres. These two points of view coincide miraculously for two-dimensional Calabi-Yau manifolds (K3 surfaces). In this case special Lagrangian fibrations for K3 surfaces under hyper\k twist reduce to the classical elliptic fibrations, which generically have 24 nodal $\mathbb{CP}^1$ singular fibres. Historically, after the SYZ conjecture was proposed, it was a common belief that 3-dimensional Calabi-Yau manifolds would exhibit similar elegant special Lagrangian torus fibration structure. More precisely, the conjectured special Lagrangian torus fibration maps for 3-dimensional Calabi-Yau manifolds should be $C^\infty$ with 1-dimensional singular locus. (It was not even clear if the singular locus was supposed to be knot, link or graph.) Our gradient flow construction in \cite{lag1} naturally produced Lagragian torus fibrations with codimension 1 singular locus that is a fattenning of a {\bf graph}. At the time they were considered {\it wrong} fibrations as far as special Lagrangian torus fibrations are concerned. Our discussion of the so-called ``expected special Lagrangian fibration structure" in \cite{lag1} based on monodromy computation was an attempt to conform to the conventional wisdom.\\

When we talked about the ``expected special Lagrangian fibration structure'', we were actually refering to the expected behavior of special Lagrangian torus fibrations that would enable SYZ mirror construction as originally proposed to work. The key ingredient of such is to expect the singular locus of the fibration to be of codimension 2, i.e., a 1-dimensional graph. In fact such expected structure may not coincide with the actual behavior of the special Lagrangian torus fibration. (Then it is necessary for the SYZ mirror construction as originally proposed to be modified.) Indeed, our gradient flow construction illustrates that Lagrangian torus fibrations of Calabi-Yau manifolds more naturally are non-$C^{\infty}$ with codimension 1 singular locus. Recent examples of D. Joyce \cite{Joyce} further indicate that such codimension 1 singular locus might be the generic behavior of special Lagrangian torus fibrations. The construction of such fibrations with graph singular locus (``expected special Lagrangian fibration structure'') is likely to be possible only in the symplectic category as shown in our construction. It is now commonly believed that such fibration structure with codimension 2 singular locus will only appear as the limiting structure of special Lagrangian torus fibrations when the Calabi-Yau manifolds approach the large complex and large radius limit.\\

With these new developments that are changing the conventional point of view, it is quite clear that our two goals (the physical goal prefering the Lagrangian fibrations that reveal the structure of conjectured special Lagrangian fibrations and the mathematical goal prefering simple and well hehaved singular locus and singular fibres) diverge somewhat for 3-dimensional Calabi-Yau manifolds. We will try to reassess these two goals more clearly. In our work, we will construct two kinds of Lagrangian torus fibrations for Calabi-Yau manifolds, which we refer to as Lagrangian torus fibration with codimension 2 singular locus (graph) and Lagrangian torus fibration with codimension 1 singular locus (which is a fattening of the graph in the previous kind).\\

From physical point of view, we believe that the Lagrangian torus fibrations with codimension 1 singular locus, which naturally come out of our gradient flow construction, will reveal the symplectic topological structures of special Lagrangian torus fibrations in SYZ conjecture, and may serve as a starting point to deform toward the actual special Lagrangian torus fibrations. As pointed out by M. Gross based on the well-known fact (Corollary \ref{ab}), if the fibration is $C^\infty$, minimality of fibres will guarantee that the singular locus is of dimension one. Discussion in section 2 further indicates that even as Lagrangian fibrations, the Lagrangian torus fibrations with codimension 1 singular locus constructed in \cite{lag1} can not be $C^\infty$ along the topological singular set (which is a union of 10 genus 6 curves in \cite{lag1}) of the fibration map. Therefore concerning smoothness, the best one can hope is to make the Lagrangian fibration map $C^\infty$ away from the topological singular set. Let $F: X \rightarrow B$ be the Lagrangian fibration map. To talk about smoothness of $F$, it is very important to determine the smooth structure of the base $B$ of the fibration (which is apriori not determined by the fibration). Some version of the SYZ conjecture requires $X$ to possess a horizontal section $S$ that intersects each fibre transversally at a regular point of the fibre with intersection number $1$. In such case, we can get around this difficulty without mentioning the smooth structure on $B$. We will call the special Lagrangian fibration $F: X\rightarrow B$ $C^\infty$ if $i\circ F: X \rightarrow X$ is a $C^\infty$ map, where $i: B\cong S\rightarrow X$ is the natural embedding of the section $S\subset X$. $i\circ F$ can be understood as mapping each fiber to its intersection point with $S$. This definition will be equivalent to specifying the smooth sructure of $B$ to ensure that the fibration map $F$ induces a diffeomorphism from $S$ to $B$. Slightly more general cases that do not require a global horizontal section are discussed in section 2. Notice that the natural Lagrangian torus fibration $F$ with codimension 1 singular locus constructed via our gradient flow method is non-$C^\infty$ on a bigger set than the topological singular set of $F$. In \cite{smooth}, we will discuss methods to make $F$ $C^\infty$ away from the topological singular set of $F$.\\

From mathematical point of view in relation to symplectic topology of Calabi-Yau manifolds, we would very much like to have an analogous picture as the beautiful case of K3 surfaces in the case of Lagrangian torus fibrations for Calabi-Yau 3-folds. Namely Lagrangian fibrations with graph singular locus and well behaved singular fibres. Such fibrations are more convenient to work with for the purposes of performing topological computations or constructing mirror manifolds symplectic topologically. Our construction of Lagrangian torus fibrations with codimension 2 singular locus (which ironically was called ``expected special Lagrangian fibration structure" in \cite{lag1}) serves such purpose. Much more technical difficulties are involved here than in the construction with codimension 1 singular locus. This construction is done in section 9 with help from section 6. From the discussion in section 2, one can see that our Lagrangian torus fibrations with codimension 2 singular locus potentially can be much smoother. In principle there is no topological obstruction away from the topological singular set of the so-called type $II$ singular fibres which are isolated. On the other hand, how smooth such fibrations should be is not entirely clear. We will discuss some partial smoothing results of such fibration in \cite{smooth}.\\

\subsection{Introduction}
In this paper we address several technical aspects that are involved in our gradient flow construction, especially those involved in \cite{lag1}. The first aspect is the dynamics of singular vector fields. The gradient vector fields we use have highly degenerate singularities and have poles. Sections 3, 4 and 5 are dedicated to a close study of the gradient flow we use. Gradient vector field is a very classical object and has been very useful in many ways, such as in Morse theory. But the kind of gradient vector field we encounter is not the ordinary nice non-degenerate gradient vector field. In our case, critical points are usually not non-degenerate, critical set is usually not even isolated. Worst of all, the function that produces the gradient vector field is not even well defined everywhere, and has infinity at some lower dimensional subsets. Clearly, special care is needed to accommodate these complexities to make sure the gradient flow behaves in the way we wanted. It turns out that two simple observations make it possible for us to handle this kind of vector fields. First, although the function and the gradient vector field could have infinities (poles), the gradient vector field can be reduced to a multiple of a $C^\infty$ vector field by a positive function (that could have infinities). Namely, as direction field, this vector field is $C^\infty$. Therefore the dynamics could be understood by analyzing the dynamics of the corresponding $C^\infty$ vector field and the positive multiple function. Second, although the singular set is usually degenerate and not isolated, the dominating term of the vector field at a singular set is homogeneous in a non-degenerate way, and have a certain structural stability.\\

The second aspect is the construction of toroidal \k metric. As we know the behavior of the flow of a vector field is very sensitive to perturbation of the vector field around singularities of the vector field, especially for those highly degenerate vector fields. For our gradient flow to behave the way we wanted it is very essential for our \k form to be toroidal as defined in section 7. (Counter-examples exist for our results on gradient flow if the \k form is not toroidal.) The construction of toroidal \k metric compatible with the underlying algebraic toroidal structure is clearly a problem of independent interest in \k geometry. In section 7, we deal with the special case of normal crossing. We prove that any \k metric can be perturbed into a toroidal metric in normal crossing case.\\

The third aspect is the results on deformation of symplectic manifold and sympletic submanifold structures as discussed in section 6. Our results here are very explicit in nature. Such results are very useful for constructing symplectomorphisms preserving certain symplectic submanifold structures, which will be used in different forms in many steps of our construction. Results in section 6 are clearly of independent fundamental interest in symplectic geometry in addition to the application to our work and worth further exploration.\\

The fourth aspect is the symplectic deformation of the Lagrangian torus fibration with codimension 1 singular locus into a Lagrangian torus fibration with codimension 2 singular locus, which serves our mathematical goal. The deformation process reduces to some miraculous computations on an explicit deformation for $\mathbb{CP}^2$ that happens to work (section 9). Such miraculous computations can be generalized to more general curves in more general toric surfaces as discussed in \cite{N} and to higher dimensions, which we hope to discuss in the future.\\

The general scheme of our construction is to start with a Lagrangian torus fibration of the large complex limit $(X_\infty, \omega_{\rm FS})$. Using results in section 7, we can construct a toroidal \k form $\omega_T$ with respect to $X_\infty \cup X_\psi$ as a small perturbation of $\omega_{\rm FS}$. With results in section 6, we may construct symplectomorphisms from $(X_\infty, X_\psi, \omega_{\rm FS})$ to $(X_\infty, X_\psi, \omega_T)$. Using such symplectomorphisms, we may first push the Lagrangian torus fibration of $(X_\infty, \omega_{\rm FS})$ to the Lagrangian torus fibration of $(X_\infty, \omega_T)$. Then we use the gradient flow (developed in section 3,4,5) under the toroidal \k metric to construct the Lagrangian torus fibration on $(X_\psi, \omega_T)$. Using those symplectomorphisms again, we can push back to get the Lagrangian torus fibration on $(X_\psi, \omega_{\rm FS})$. If we start with the natural Largrangian torus fibration of $(X_\infty, \omega_{\rm FS})$ determined by the moment map, we will get the Lagrangian torus fibration of $(X_\psi, \omega_{\rm FS})$ with codimension 1 singular locus. If we want to construct the Lagrangian torus fibration of $(X_\psi, \omega_{\rm FS})$ with codimension 2 singular locus, we have to start with the Largrangian torus fibration of $(X_\infty, \omega_{\rm FS})$ constructed explicitly in section 9.\\

In Section 2 we clarify some philosophical points that mainly serve our physical goal. For example the advantage of considering Lagrangian fibration in comparison to non-Lagrangian topological fibration, which is even easier to construct, and the difference between $C^\infty$ Lagrangian fibration and piecewise smooth Lipschitz continuous Lagrangian fibration. More precisely, the constraints on the topological type of singular fibres of $C^\infty$ Lagrangian fibrations are discussed. As an application, the Lagrangian torus fibration for Fermat type quintics with codimension 1 singular locus constructed in \cite{lag1} can not be $C^\infty$. Section 2 is logically independent of the rest of this paper.\\

Sections 3, 4 and 5 are devoted to the discussion of the behavior of our gradient flow under toroidal \k metric. In section 3, we discuss the direction field, which is very helpful in understanding the dynamics of vector fields with infinity. We also give explicit solutions of several local examples that serve as local models of our gradient flow. Through the explicit solutions of these local models, we can already see the non-smoothness of the Lagrangian fibration constructed via gradient flow. In section 4, we discuss perturbation of nondegenerate homogeneous hyperbolic vector fields in general. In section 5, we discuss perturbation of local models discussed in section 3 and \cite{lag1}. Although our primary interest here is the gradient flow associated with a family of Calabi-Yau hypersurfaces, our work on gradient flow can be formulated into a general theorem (theorem \ref{de}) that can be applied to more general situations.\\

The content of section 6 is already mentioned in ``the third aspect". There are two places where we will apply the results in this section. The first application is to construct a $C^{0,1}$ symplectomorphism from $\mathbb{CP}^4$ with the Fubini-Study metric to $\mathbb{CP}^4$ with toroidal metric with respect to $X_\infty\cup X_\psi$ that maps $X_\infty\cup X_\psi$ to itself (theorem \ref{eg}). The second application is to deform the symplectic curves in $\mathbb{CP}^2$'s to achieve graph image (theorem \ref{ee} and corollary \ref{er}), and extend to $X_\infty \subset \mathbb{CP}^4$ (theorem \ref{ef}).\\

In Section 7 we discuss the construction of the so-called toroidal \k metric as mentioned in ``the second aspect". Section 7 does not depend on the rest of the paper. The reader should feel free to refer to it when necessary.\\

In section 8, we formulate and prove one of the main theorems (theorem \ref{ha}) of this paper. This theorem establishes a symplecticomorphism from a smoothing of a normal crossing variety to the normal crossing variety itself through gradient flow deformation under fairly general conditions. This theorem has wide potential of applications, for example, to construct Lagrangian fibration on a smoothing of a normal crossing variety based on a Lagrangian fibration structure of the normal crossing variety. Such applications particularly include the construction of Lagrangian torus fibration for Calabi-Yau hypersurfaces in toric varieties. In section 8, as application of theorem \ref{ha}, we construct the Lagrangian torus fibrations for Fermat type quintic Calabi-Yau hypersurfaces in $\mathbb{CP}^4$ with codimension 1 singular locus that was discussed in \cite{lag1}.\\

In section 9, we embark on the construction of Lagrangian torus fibration for Fermat type quintic Calabi-Yau in $\mathbb{CP}^4$ with codimension 2 singular locus. The key ingredient is the computations that lead to lemma \ref{gd}. Everything else is more or less routine given the results proved in previous sections.\\

{\bf Remark:} Our gradient flow method also provides a direct approach to construct {\bf non-Lagrangian} torus fibrations on quintic Calabi-Yau manifolds with codimension 2 singular locus {\bf constructively}. Such topological construction should be much easier than our symplectic construction. (In topological situation, instead of more difficult arguments of perturbed dynamical systems, topological cut and paste based on the local models should be sufficient. In deformation to codimension 2 singular locus, modification of ${\cal F}_t$ is unnecessary and corresponding extension to non-Lagrangian diffeomorphism is almost trivial to do.)\\

{\bf Remark:} Although the original purpose of this paper is to deal with problems of analysis, symplectic geometry and dynamical systems nature which are necessary for our gradient flow method of constructing Lagrangian torus fibrations on Calabi-Yau hypersurfaces, many of the results and methods developed on gradient flow, Hamiltonian deformation of submanifold of symplectic manifold, smoothing of Lagrangian fibration are also interesting in their own right. We hope to have more discussion of these methods and their applications in the future.\\

As we claimed, our method can be used to construct Lagrangian torus fibrations for general quintic Calabi-Yau hypersurfaces in ${\mathbb{CP}^4}$ and more generally for Calabi-Yau hypersurfaces, even Calabi-Yau complete intersections in toric variety. These constructions and their applications to SYZ mirror construction will be discussed in subsequent papers (\cite{lag1,tor,ci}).\\

\se{Smoothness of Lagrangian fibration}
Let us start with a well-known fact for Lagrangian fibrations. Let $(X,\omega)$ be a smooth symplectic manifold. A fibration $F: X\rightarrow B$ is called a $C^l$-Lagrangian fibration, if $F$ is a $C^l$ map and the smooth part of each fibre is Lagrangian. The following result is well known.\\
\begin{theorem}
Let $F: X\rightarrow B$ be a $C^{1,1}$-Lagrangian fibration, then for any $b\in B$, there is an action of $T^*_bB$ on $F^{-1}(b)$.\\
\end{theorem}
{\bf Proof:} For any closed 1-form $\alpha$ on $B$, $F^*\alpha$ is a closed $C^{0,1}$ 1-form, there is a corresponding $C^{0,1}$ Hamiltonian vector field $ H_{\alpha}$ that satisfies $F^*\alpha = i(H_{\alpha})\omega$. $H_{\alpha}$ will be along the fibres, because for any $W$ along the fibres, we have\\
\[
\omega(H_{\alpha}, W) = \langle F^*\alpha, W\rangle = 0.
\]
Then for closed 1-forms $\alpha_1$ $\alpha_2$ on $B$, we have\\
\[
i([H_{\alpha_1},H_{\alpha_2}])\omega = {\cal L}_{H_{\alpha_1}}(i(H_{\alpha_2})\omega) - i(H_{\alpha_2}){\cal L}_{H_{\alpha_1}}\omega = d(\omega(H_{\alpha_1},H_{\alpha_2}))=0.
\]
Namely\\
\[
[H_{\alpha_1},H_{\alpha_2}] =0.
\]
Therefore $T^*_bB$ acts on $F^{-1}(b)$ as Lie algebra. Exponentiating the Lie algebra action, which amounts to considering the flow corresponding to Hamiltonian vector field $H_{\alpha}$, will give us the action in the theorem.
\begin{flushright} $\Box$ \end{flushright}
{\bf Remark:} Here it is very crucial for $H_{\alpha}$ to be $C^{0,1}$, which will ensure that the solution curves of the Hamiltonian vector field $H_{\alpha}$ are uniquely determined by their initial values.\\

From this fact, it is rather tempting to expect that for any $b \in B$, the corresponding fibre $F^{-1}(b)$ is stratified by a disjoint union of finitely many orbits $\{O_l\}$ of $T^*_bB$, where each $O_l$ is diffeomorphic to $(S^1)^k\times\mathbb{R}^m$ for some $k+m\leq 3$. Unfortunately, this is not the case in general. For example the $I^*_0$ singular fibre in the Kodaira's classification of singular fibres of elliptic surfaces (as pointed out to me by M. Gross) is a counter example. The main reason is that one component of $I^*_0$ singular fibre has multiplicity 2, therefore is a union of infinitely many 0-dimensional orbits. To avoid this annoying irregularity, it is necessary to require the fibration to be of ``multiplicity one'' in a certain sense. Recall that the SYZ torus fibration of a \cy manifold requires the existence of a section for the fibration. This requires the component of the fibre intersecting the section to be of multiplicity one. In particular, for the generic case where each fibre only contains one irreducible component, every fibre is required to be of multiplicity one. From this perspective, multiplicity one condition is very natural for our application.\\
\begin{de}
\label{aa}
For a $C^l$-Lagrangian fibration $F: X\rightarrow B$, $x\in X$ is called {\bf regular} if there exists a small neighborhood $U_x$ of $x$ such that $F^{-1}(b)\cap U_x$ is a manifold for every $b$. Each connected component of the open set of regular points in $F^{-1}(b)$ is called a {\bf regular component} of $F^{-1}(b)$.\\
A regular point $x\in X$ is called {\bf simple} if there is a local section $(x\in)S_x\subset U_x$ transversal to the fibration such that $F|_{S_x}: S_x \rightarrow F(S_x) \subset B$ is a $C^l$-diffeomorphism.\\
A regular component of $F^{-1}(b)$ is called simple if it consists of only simple regular points.\\
$F$ is called a simple $C^l$-Lagrangian fibration if any regular point $x\in X$ is simple.\\
\end{de}
{\bf Remark:} Differentiability of $F$ very much depends on the smooth structure of $B$. Even for a trivial family over a 1-dimensional base $B$, change of smooth structure of $B$ via the homeomorphism $y=x^{\frac{1}{3}}$ will clearly destroy the smoothness of $F$. To avoid this kind of artificial non-smoothness, it is crucial to require the base $B$ to possess smooth structures inherited from the smooth structure of the smooth horizontal section that intersects each fibre transversally with intersection number $1$. {\bf Simple} condition in the previous definition is just a combination of such canonical choice of smooth structure on $B$ and the multiplicity $1$ condition of the fibres.\\

\begin{prop}
A regular component of $F^{-1}(b)$ is simple if and only if there exists one simple regular point in it.
\end{prop}
\begin{flushright} $\Box$ \end{flushright}
From now on, when we talk about a Lagrangian fibration $F: X\rightarrow B$, we will always assume that at least one component of the fibre is of multiplicity $1$ and the smooth structure on $B$ (which apriori is not determined by the fibration) is determined by the following assumption$\star$. When each fibre is irreducible, such smooth structure of $B$ is uniquely determined.\\

{\bf Assumption$\star$:} Locally under $F$ the base $B$ is diffeomorphic to certain smooth horizontal section that intersects each fibre transversally with intersection number $1$.\\

\begin{co}
\label{ab}
Assume that $F: X\rightarrow B$ is a simple $C^l$-Lagrangian fibration. Then for any $b \in B$,\\
\[
{\rm Reg}(F^{-1}(b)) = \bigcup_l O_l
\]
is a disjoint union of orbits of $T^*_bB$, where each $O_l$ is diffeomorphic to $(S^1)^k\times\mathbb{R}^m$ for some $k+m = 3$.
\end{co}
\begin{flushright} $\Box$ \end{flushright}
{\bf Remark:} When $F$ is generic in certain sense, we expect $\displaystyle F^{-1}(b) = \bigcup_l O_l$ to be a disjoint union of orbits of $T^*_bB$, where each $O_l$ is diffeomorphic to $(S^1)^k\times\mathbb{R}^m$ for some $k+m \leq 3$.\\

Recall from \cite{lag1}, the flow of $V=\frac{\nabla f}{|\nabla f|^2}$ induces Lagrangian fibration $F: X_{\psi} \rightarrow \partial \Delta$. According to theorem 3.1 in \cite{lag1}, singular locus of the fibration in $\partial \Delta$ is\\
\[
\tilde{\Gamma} = \tilde{\Gamma}^0 \cup \tilde{\Gamma}^1 \cup \tilde{\Gamma}^2.
\]
For $p\in \tilde{\Gamma}^0$, $F^{-1}(p)$ is a Lagrangian 3-torus with $5$ two-tori collapsed to $5$ singular points. It can be thought of as a union of 5 points and 5 $(S^1)^2\times\mathbb{R}^1$'s. This is a possible singular Lagrangian fibre type of $C^{1,1}$ Lagrangian fibration as described in the previous corollary. For $p\in \tilde{\Gamma}^2$, $F^{-1}(p)$ is a Lagrangian 3-torus with $50$ circles collapsed to $50$ singular points. For $p\in \tilde{\Gamma}^1$, $F^{-1}(p)$ is a Lagrangian 3-torus with $25$ circles collapsed to $25$ singular points. It is not hard to see that these two types of singular fibres are irreducible and their regular components do not admit an action by $\mathbb{R}^3$, therefore cannot be expressed as union of orbits described in the previous corollary. With the understanding that the manifold structure of $B$ is determined by assumption$\star$, this analysis implies that\\
\begin{co}
The Lagrangian fibration $F$ for Fermat type quintics with codimension 1 singular locus constructed in \cite{lag1} is not $C^{1,1}$, in particular, $F$ is not a smooth ($C^{\infty}$) map.
\end{co}
\begin{flushright} $\Box$ \end{flushright}
In the case of SYZ conjecture, the special Lagrangian torus fibration $F: X\rightarrow B$ is required to possess a canonical special Lagrangian horizontal section $i: B \rightarrow S \subset X$ that intersects each fibre transversely at a regular point of the fibre with intersection number $1$. This will enable us to get around the mentioning of smooth structure of $B$ entirely. We will call the special Lagrangian fibration $F: X\rightarrow B$ $C^{l,\alpha}$ if $i\circ F: X \rightarrow X$ is a $C^{l,\alpha}$ map. $i\circ F$ can be understood as mapping each fiber to its intersection point with $S$. This definition will be equivalent to specifying the smooth sructure of $B$ using assumption$\star$ with respect to $S$ and requiring $F: X\rightarrow B$ to be a simple $C^{l,\alpha}$-Lagrangian fibration in the sense of definition \ref{aa}.\\

It was commonly believed (before the examples of Joyce!) that the Lagrangian fibrations in SYZ construction are $C^{\infty}$ maps. In particular the singular fibres should have structures as described in previous theorem. It is not hard to see that the singular fibres in the discussion of expected special Lagrangian fibration in \cite{lag1} have the topological types of singular Lagrangian fibres of $C^{1,1}$ Lagrangian fibrations as described in corollary \ref{ab}. The type $I_5$ fibre is a union of 5 circles and 5 $(S^1)^2\times\mathbb{R}^1$'s, with Euler number equal to zero. The type $II_{5\times 5}$ fibre is a union of 50 points, 75 $\mathbb{R}^1$'s and 25 $S^1\times\mathbb{R}^2$'s, with Euler number equal to $-25$. The type $III_5$ fibre is a union of 5 points and 5 $(S^1)^2\times\mathbb{R}^1$'s, with Euler number equal to $5$.\\

Indeed, in theorem \ref{ga}, we are able to construct a piecewise smooth Lagrangian torus fibration with graph singular locus for Fermat type quintic with such topological structure. One may try to use this construction as a starting point to deform to a smooth ($C^{\infty}$) Lagrangian fibration. It is not clear how smooth a Lagrangian fibration with such topological structure can be. Our guess is that the Lagrangian fibration at least can be made $C^{\infty}$ away from type $II_{5\times 5}$ fibres. We will explore this further in \cite{smooth}.\\

In common sense, it is not hard to modify a non-smooth map to a smooth map by small perturbation. In fact, if one does not worry about Lagrangian condition, it is very easy to modify our fibration map $F$ to a smooth fibration map with the same singular fibre structure. On the other hand, if one wants to keep the Lagrangian condition, getting a smooth Lagrangian fibration can not be achieved by just a small perturbation! The reason is that the singular fibres in our construction are topologically different from the possible singular fibres in a smooth ($C^{\infty}$) Lagrangian fibration as described in corollary \ref{ab}. Clearly a small perturbation will not be enough to achieve this kind of topological change. This observation shows a major difference between general smooth maps and smooth Lagrangian maps. Denote ${\cal LF}^{l,\alpha}$ as the set of simple Lagrangian fibration maps, whose $l$-th derivatives are ${\alpha}$-H\"{o}lder continuous. Under this notation, the above observation can be rephrased more generally as\\
\begin{co}
The closure of ${\cal LF}^{1,1}$ in ${\cal LF}^{0,1}$ with respect to the $C^0$-topology is not equal to ${\cal LF}^{0,1}$.
\end{co}
\begin{flushright} $\Box$ \end{flushright}
{\bf Remark:} Even if we remove the simple condition in the definition of ${\cal LF}^{l,\alpha}$, we believe the statement in the above corollary should still be correct. Although this more general situation will not concern our application.\\

With this observation in mind, we believe a general strategy to deform a non-$C^\infty$ Lagrangian fibration to a $C^\infty$ agrangian fibration should contain two steps. In the first step we modify the Lagrangian fibration into a Lagrangian fibration with the expected topological structure usually by a non-small change (which usually is still non-smooth). Namely, make the singular locus to be of codimension 2 and the singular fibres to be of the types described in Corollary \ref{ab}. In the second step we modify the (possibly non-smooth) Lagrangian fibration with expected topological structure into a smooth Lagrangian fibration with the same topological structure by a small perturbation. We will call the first step {\bf topological modification} and the second step {\bf analytical modification}. In  section 9, we will achieve the topological modification via the approach of Hamiltonian deformation of submanifolds of symplectic manifold.\\

\se{Direction field and local models}
In this section we will discuss local examples that serve as local models of our gradient vector fields. In section 5, we will prove perturbation stability of these local models, which will enable us to understand precisely the behavior of our global gradient flow and the Lagrangian fibration structure produced via gradient flow. In principle, the readers may be able to skip the discussion of local models in this section, which logically is not absolutely needed for the more general discussions in later sections. In reality, the structure of our gradient flow around singularities is very complicated, while our local models can be solved explicitly and provide clear pictures of the local structure of our gradient flow around singularities. Furthermore, the solutions of these local models will provide guidance and motivation for many rather involved and technical arguments in section 5. We also want to point out that these local examples have independent interest of their own as dynamical systems of homogeneous vector fields.\\

We will solve these local models explicitly in this section. As we claimed previously the fibrations we get are only piecewise smooth. This effect can already be seen through the explicit solutions of the local models computed here.\\

In this section, we will also introduce the concept of direction field. In our situation, although the ambient space ${\mathbb{CP}^4}$ is smooth, the gradient vector field $ \nabla f$ is very singular. It is singular when $f=\rm{Re}(s)$ is singular (namely, along $X_\infty$). Even $V=\frac{\nabla f}{|\nabla f|^2}$ (the vector field we actually use) is singular at the singular set of $X_{\infty}$. The corresponding dynamical system is better understood through direction field. In general the concept of direction field is a more suitable setting for discussing the dynamical systems of this kind of singular vector fields.\\

\subsection{Direction field}
Let $X$ be a smooth manifold. By a {\bf direction field} $[V]$ on $X$, we mean the equivalence class of its representatives. $\{V_\alpha\}$ is called a {\bf representative} of $[V]$ if there is a coordinate chart $\{U_{\alpha}\}$ and vector field $V_{\alpha}$  defined in an open and dense subset of $U_{\alpha}$ such that $V_{\alpha}=\rho_{\alpha \beta}V_{\beta}$ with $\rho_{\alpha \beta}>0$ in an open and dense subset of $(U_{\alpha}\cap U_{\beta})$. Two representatives $\{V^i_\alpha\}$ with respect to the covers $\{U^i_\alpha\}$ for $i=1,2$ are said to be equivalent if $V^1_{\alpha}=\rho_{\alpha \beta}V^2_{\beta}$ with $\rho_{\alpha \beta}>0$ in an open and dense subset of $(U^1_{\alpha}\cap U^2_{\beta})$. Namely a direction field $[V]$ is roughly a vector field up to scale multiplication by a positive function. Let $E_V$ be the {\bf singular set} of $[V]$. $p\in X\backslash E_V$ if and only if there exists a neighborhood $U_p$ and coordinate $x$ and representative $V_p$ of $[V]$ on $U_p$ such that $V_p=\frac{\partial}{\partial x_1}$. Clearly, this is an open condition that will guarantee $E_V$ to be closed. We will always assume that the {\bf regular set} $X\backslash E_V$ is open and dense. $[V]$ is called {\bf regular} if $E_V = \emptyset$. $\{V_\alpha\}$ is called a {\bf smooth representative} of $[V]$ if each $V_{\alpha}$ is a smooth vector field on $U_{\alpha}$ that is nonvanishing on $U_{\alpha}\backslash E_V$ and $V_{\alpha}=\rho_{\alpha \beta}V_{\beta}$ with $\rho_{\alpha \beta}>0$ on $(U_{\alpha}\cap U_{\beta})\backslash E_V$. A direction field with smooth representatives will be called {\bf smooth direction field}. We will mainly concern smooth direction fields in our work. Later when we talk about a direction field, we always refer to a particular representative of it.\\

{\bf Remark:} As we mentioned, direction field is a more proper setting to understand the essence of a dynamical system. For instance, the dynamical systems of two vector fields corresponding to the same direction field will have the same orbit structure. Their difference only comes from reparametrization. Direction field is extremely suitable for the discussion of dynamics of singular vector fields with poles.\\

{\bf Example:} By our definition, for $V=|\frac{x_2}{x_3}|\frac{\partial}{\partial x_1}$, $E_V$ is empty. And for $V=x_1\frac{\partial}{\partial x_2}+x_2\frac{\partial}{\partial x_1}$, $E_V=(0,0)$.
\begin{flushright} $\Box$ \end{flushright}
{\bf Remark:} There is a related concept called {\bf line field}, which can be similarly defined as direction field, by replacing every ``$>0$" in the above definition by ``$\not=0$". There are similarly the concepts of smooth line field, singular set of line field, etc. For example $V=\frac{x_2}{x_3}\frac{\partial}{\partial x_1}$ is a smooth line field, but not a smooth direction field. The term direction field is employed in \cite{A} to refer to a regular line field $[V]$ in our notation (with the singular set $E_V = \emptyset$).\\

Given a direction field $[V]$, we can consider its orbits. $\phi:\mathbb{R}\rightarrow X$ is called an orbit of $[V]$ if there is a representative $V$ of $[V]$ such that $\frac{d\phi}{dt} = V(t)$. $\phi$ is called a complete orbit if $\phi(+\infty),\phi(-\infty)\in E_V$ (when the limits exist). A non-complete orbit can always be extended to a complete one. (Sometimes reparametrization is needed.)\\

Direction fields and their orbits are quite useful for us to understand the dynamic systems of vector fields with singularities. Around a regular point of a direction field, the orbit structure is very simple. The main difficulty is to understand the behavior of a direction field and its orbits near the singular set.\\

A particular type of direction field is the gradient vector field $\nabla f$ of a real singular function $f$ with respect to some metric (e.g. flat metric). It is usually convenient to consider $V=\frac{\nabla f}{|\nabla f|^2}$, since $V$ moves level sets of $f$ to level sets. A particularly interesting case is when $f$ is the real part of a meromorphic function $s=f+ig$ and the metric is K\"{a}hler. Then the gradient of $f$ with respect to the metric is the same as the Hamiltonian vector field of $g$ with respect to the K\"{a}hler form. In particular, $g$ is constant along the flow. (See \cite{lag1}.)\\

More precisely, we can write $f(z)={\rm Re}(\frac{p(z)}{q(z)})$, where $p(z),q(z)$ are holomorphic functions on ${\mathbb{C}^n}$. We have the direction field $V=\frac{\nabla f}{|\nabla f|^2}$. Consider fibration $F: {\mathbb{C}^n}\backslash \{p=q=0\} \rightarrow {\mathbb{CP}^1}$, where $u=F(z)= \frac{p(z)}{q(z)}$. Then we have\\
\begin{lm}
$V$ is perpendicular to the fibres of $F$ and $F_*V={\rm Re}(\frac{\partial}{\partial u})$.
\end{lm}
{\bf Proof:} $V$ as direction field of $\nabla f$ clearly is perpenticular to fibres of $F$, which belong to level sets of $f$. By lemma 3.1 in \cite{lag1}, $V$ will leave $\{{\rm Im}(u)={\rm constant}\}$ invariant. Clearly, $V(u) = 1$. Therefore $F_*V={\rm Re}(\frac{\partial}{\partial u})$.
\begin{flushright} $\Box$ \end{flushright}
\begin{lm}
$V$ is a smooth direction field.
\end{lm}
{\bf Proof:}
\[
V = \frac{1}{|\nabla f|^2|q|^4}(|q|^2\nabla ({\rm Re}(p\bar{q})) - {\rm Re}(p\bar{q})\nabla (|q|^2)).
\]
It is clearly a smooth direction field.
\begin{flushright} $\Box$ \end{flushright}
In the following, we will discuss several examples that will serve as local models of our gradient flow. Later in sections 4 and 5, we will discuss the perturbation stability of these local models that will enable us to derive the structure theorems of the Lagrangian fibrations constructed via gradient flow. To illustrate the idea, we  start with the simplest non-trival example.\\

\subsection{The baby case}
{\bf Example 3.1:} On ${\mathbb{C}^3}$, consider $\displaystyle f(z)={\rm Re}\left(\frac{z_2 z_3}{z_1}\right)$, then\\
\[
\nabla f ={\rm Re}\left(\frac{\bar{z}_1\bar{z}_2 \frac{\partial}{\partial z_3} + \bar{z}_1\bar{z}_3 \frac{\partial}{\partial z_2}-\bar{z}_2\bar{z}_3 \frac{\partial}{\partial z_1}}{\bar{z}_1^2}\right),
\]
\[
V={\rm Re}\left(z_1^2\frac{\bar{z}_1\bar{z}_2 \frac{\partial}{\partial z_3} + \bar{z}_1\bar{z}_3 \frac{\partial}{\partial z_2}- \bar{z}_2\bar{z}_3 \frac{\partial}{\partial z_1}}{|z_2|^2|z_1|^2 + |z_2|^2|z_3|^2 + |z_3|^2|z_1|^2}\right).
\]
Vector field $V$ deforms the hypersurface $z_2z_3=0$ to the hypersurface $z_1=0$. As we know, $V$ leaves $g(z)={\rm Im}(\frac{z_2 z_3}{z_1})$ invariant. Restricted to the real hypersurface $g(z)={\rm Im}(\frac{z_2 z_3}{z_1})=0$, $V$ will be moving among complex hypersurfaces $X_c$ defined by $z_2z_3-cz_1=0$ for $c$ real.\\

We are interested in looking for invariant functions of $V$ that restrict to certain initial values on $z_2z_3=0$. An obvious invariant function is
\[
\phi_1(z)= |z_2|^2-|z_3|^2.
\]
$(|z_2|^2-|z_3|^2,{\rm Re}(z_1))$
gives a Lagrangian fibration of $X_0$. Vector field $V$ induces symplectic morphism (see definition \ref{df}) from $X_c$ for $c$ real to $X_0$ and will induce Lagrangian fibration on $X_c$ for $c$ real. To make the fibration map explicit, we need to find in addition to $\phi_1$ another $V$-invariant function $\phi_2$ such that $\phi_2|_{X_0}={\rm Re}(z_1)$. It is rather difficult to compute $\phi_2$ from the expression of $V$ from above. To be able to compute $\phi_2$ explicitly, it is important to notice that only the restriction of $\phi_2$ to $\{g=0\}=\cup_{c\in\mathbb{R}}X_c$ is determined by the initial value. When restricted to $\{g=0\}=\cup_{c\in\mathbb{R}}X_c$, we may assume that $z_1\bar{z}_2\bar{z}_3\in \mathbb{R}$. This gives us the following crucial simplification of $V$ which enables us to compute $\phi_2$.

\[
V=\frac{z_1 \bar{z}_2\bar{z}_3}{|z_2|^2|z_1|^2 + |z_2|^2|z_3|^2 + |z_3|^2|z_1|^2}{\rm Re}\left( \frac{|z_1|^2}{\bar{z}_3} \frac{\partial}{\partial z_3} + \frac{|z_1|^2}{\bar{z}_2}  \frac{\partial}{\partial z_2}- z_1\frac{\partial}{\partial z_1}\right).
\]

With moments of thought, we get

\[
\phi_2(z) = \left\{
\begin{array}{l}
{\rm Re}(z_1)\left(1+\left|\frac{z_2}{z_1}\right|^2\right)^{\frac{1}{2}}\ {\rm when}\ |z_2|\leq|z_3|,\\
{\rm Re}(z_1)\left(1+\left|\frac{z_3}{z_1}\right|^2\right)^{\frac{1}{2}}\ {\rm when}\ |z_3|\leq|z_2|.
\end{array}
\right.
\]

In general, there is a more systematical way of finding invariant functions. Notice that from the first expression of $V$ we can see that on ${\mathbb{C}^3}$ we have the invariant functions

\[
|z_2|^2-|z_3|^2,\ \ |z_1|^2+|z_2|^2,\ \ |z_1|^2 +|z_3|^2.
\]

When restricted to $\{g=0\}$, from the second expression of $V$, on $\{g=0\}$ we have the invariant functions

\[
{\rm Im}(\log z_i)=-i(\log z_i - \log |z_i|),\ \frac{z_i}{|z_i|},\ \ {\rm for}\ \ i=1, 2, 3.
\]

On the hyperplane $z_3=0$, we have coordinate $(z_1,z_2)$. Coordinate functions can be invariantly extended to

\begin{eqnarray*}
Z_1 &=& z_1\left(1+\left|\frac{z_3}{z_1}\right|^2\right)^{\frac{1}{2}} = \frac{z_1}{|z_1|}(|z_1|^2 + |z_3|^2)^{\frac{1}{2}} \ \ {\rm when}\ |z_3|\leq|z_2|,\\
Z_2 &=& z_2\left(1-\left|\frac{z_3}{z_2}\right|^2\right)^{\frac{1}{2}} = \frac{z_2}{|z_2|}(|z_2|^2 - |z_3|^2)^{\frac{1}{2}} \ \ {\rm when}\ |z_3|\leq|z_2|.
\end{eqnarray*}
On the hyperplane $z_2=0$, we have coordinate $(z_1,z_3)$. Coordinate functions can be invariantly extended to\\
\begin{eqnarray*}
Z_1 &=& z_1\left(1+\left|\frac{z_2}{z_1}\right|^2\right)^{\frac{1}{2}} = \frac{z_1}{|z_1|}(|z_1|^2 + |z_2|^2)^{\frac{1}{2}} \ \ {\rm when}\ |z_2|\leq|z_3|,\\
Z_3 &=& z_3\left(1-\left|\frac{z_2}{z_3}\right|^2\right)^{\frac{1}{2}} = \frac{z_3}{|z_3|}(|z_3|^2 - |z_2|^2)^{\frac{1}{2}} \ \ {\rm when}\ |z_2|\leq|z_3|.
\end{eqnarray*}\\
$\phi =(\phi_1, \phi_2)$ determines a Lagrangian fibration $\phi: X_c =\{z_2z_3 =cz_1\} \rightarrow {\mathbb{R}^2}$ for any $c\in \mathbb{R}$. When $c\not=0$, $X_c$ is smooth. For $c(\not=0)\in \mathbb{R}$, $\phi|_{X_c}$ is the \l fibration of $X_c$ induced from the \l fibration $(|z_2|^2-|z_3|^2,{\rm Re}(z_1))$ of $X_0$ by the flow of $V$. $\phi$ is not a smooth map, it is only piecewise smooth (Lipschitz). More precisely, $\phi_1$ is a smooth function on $X_c$. $\phi_2$ is smooth on $X_c$ away from the real hypersurface $\{\phi_1^{-1}(0)\}\cap X_c$. It is easy to see that the restriction of $\phi_2$ to the real hypersurface $\{\phi_1^{-1}(0)\}\cap X_c$ is also smooth. Therefore non-smoothness of $\phi$ comes from non-smoothness of $\phi_2$ on the normal direction of the hypersurface $\{\phi_1^{-1}(0)\}\cap X_c$. Since each fiber is a smooth submanifold, one may expect that reparametrization of coordinates in the image ${\mathbb{R}^2}$ will make the fibration map smooth. It turns out this is not the case.\\
\begin{prop}
\label{bb}
The tangent spaces of the fibres of the fibration\\
\[
\phi: X_c =\{z_2z_3 =cz_1\} \rightarrow {\mathbb{R}^2}
\]
vary Lipschitz continuously. But for any coordinate change on ${\mathbb{R}^2}$, the map $\phi$ can at most be made Lipschitz continuous.\\
\end{prop}
The first statement is rather natural. Since the fibration map is piecewise smooth, it is reasonable to expect that the tangent space of the fibres of the fibration vary piecewise smoothly, in particular vary Lipschitz continuously.\\

The second statement at first glance looks a little surprising. Since tangent map is the first derivative of the fibration map, naively, one may expect that the Lipschitz $(C^{0,1})$ continuity when varying the tangent spaces of fibres will make fibration map at least $C^{1,1}$ continuous after maybe suitable change of coordinates on the image ${\mathbb{R}^2}$. We will show by explicit computation that this is not the case.\\

{\bf Proof of proposition \ref{bb}:} Choose coordinate $(r_1, \theta_1, r_2, \theta_2)$ on $X_c$, then we have\\
\[
\phi_1 = r_2^2 - \frac{c^2r_1^2}{r_2^2}.
\]\\
\[
\phi_2 = \left\{
\begin{array}{l}
\cos \theta_1 (r_1^2 + r_2^2)^{\frac{1}{2}}\ {\rm when}\ r_2^2 \leq cr_1,\\\\
\cos \theta_1 (r_1^2 + \frac{c^2r_1^2}{r_2^2})^{\frac{1}{2}}\ {\rm when}\ r_2^2 \geq cr_1.
\end{array}
\right.
\]
When $r_2^2 \leq cr_1$,\\
\[
d\phi_2 = (r_1^2 + r_2^2)^{\frac{1}{2}}d(\cos \theta_1) + \frac{1}{2}(r_1^2 + r_2^2)^{-\frac{1}{2}}\cos \theta_1 d(r_1^2 + r_2^2).
\]
When $r_2^2 \geq cr_1$,\\
\[
\phi_2 = \cos \theta_1 (r_1^2 + r_2^2 - \phi_1)^{\frac{1}{2}}.
\]\\
\[
d\phi_2 = (r_1^2 + r_2^2 - \phi_1)^{\frac{1}{2}}d(\cos \theta_1) + \frac{1}{2}(r_1^2 + r_2^2 - \phi_1)^{-\frac{1}{2}}\cos \theta_1 d(r_1^2 + r_2^2 - \phi_1).
\]

Let

\[
\alpha = \left\{
\begin{array}{ll}
(r_1^2 + r_2^2)^{\frac{1}{2}}d(\cos \theta_1) + \frac{1}{2}(r_1^2 + r_2^2)^{-\frac{1}{2}}\cos \theta_1 d(r_1^2 + r_2^2)\ &{\rm when}\ r_2^2 \leq cr_1\\\\
(r_1^2 + r_2^2 - \phi_1)^{\frac{1}{2}}d(\cos \theta_1) + \frac{1}{2}(r_1^2 + r_2^2 - \phi_1)^{-\frac{1}{2}}\cos \theta_1 d(r_1^2 + r_2^2)\ &{\rm when}\ r_2^2 \geq cr_1
\end{array}
\right.
\]\\\\
Since $\phi_1$ vanishes at the boundary of the two regions, it is easy to see that $\alpha$ is Lipschitz continuous. On the other hand, $\alpha$ can also be written as\\\\
\[
\alpha = \left\{
\begin{array}{ll}
d\phi_2\ &{\rm when}\ r_2^2 \leq cr_1,\\\\
d\phi_2 - \frac{1}{2}(r_1^2 + r_2^2 - \phi_1)^{-\frac{1}{2}}\cos \theta_1 d\phi_1\ &{\rm when}\ r_2^2 \geq cr_1.
\end{array}
\right.
\]\\
Therefore,\\
\[
\alpha = d\phi_2\ \ \ ({\rm mod}\ d\phi_1).
\]
The distribution of tangent spaces of fibres of the fibration is determined by ${\rm span}\langle d\phi_1,d\phi_2\rangle = {\rm span}\langle d\phi_1, \alpha\rangle $. Since $d\phi_1$ is smooth and $\alpha$ is Lipschitz continuous, the tangent spaces of fibres of the fibration vary Lipschitz continuously.\\\\
The second statement of the proposition can be seen as follows. Take a smooth horizontal section (for example the real locus, which under coordinate $(z_1,z_2)$ is characterized as $\theta_1=\theta_2=0$ and can be parametrized by $(\rho_1,\rho_2)$, where $\rho_i=|z_i|$), this section will determine a natural coordinate for the base via fibration map. (In the case of real locus, denote the corresponding coordinate $(\rho_1,\rho_2)$.) If there exists a coordinate on the image ${\mathbb{R}^2}$ such that the fibration map is $C^{l,a}$, then by definition, the transform between $(\rho_1,\rho_2)$ and this coordinate will be of class $C^{l,a}$. Namely smooth section coordinate $(\rho_1,\rho_2)$ will always exhibit optimal smoothness of the fibration map. Under $(\rho_1,\rho_2)$ coordinate, the fibration map can be written as $\rho (z) = (\rho_1(z),\rho_2(z))$ that satisfies $\phi(\rho(z)) = \phi(z)$, which can be reduced to

\[
\rho_2^2 - \frac{c^2\rho_1^2}{\rho_2^2} = r_2^2 - \frac{c^2r_1^2}{r_2^2}.
\]\\
\[
\rho_1^2 + \rho_2^2  = \left\{
\begin{array}{ll}
(r_1^2 + r_2^2) (\cos \theta_1)^2\ &{\rm when}\ r_2^2 \leq cr_1,\\\\
(r_1^2 + r_2^2)(\cos \theta_1)^2 + (r_1^2 - \frac{c^2r_1^2}{r_2^2}) (\sin \theta_1)^2\ &{\rm when}\ r_2^2 \geq cr_1.
\end{array}
\right.
\]\\
From these expressions, it is not hard to see that $\rho(z_1,z_2)$ is at most Lipschitz continuous when $\sin \theta_1\not=0$ (or $z_1$ is not real) and $r_2^2 = cr_1$.
\begin{flushright} $\Box$ \end{flushright}
It is interesting to observe that $(z_2,z_3)$ is a nice global coordinate for $X_c$, since we can always solve for $z_1=z_2z_3/c$. Under this coordinate $\phi_2$ can simply be expressed as\\
\[
\phi_2 = {\rm Re}\left(\frac{z_2z_3}{r_{23}}\right)\left(1+\frac{r_{23}^2}{c^2}\right)^{\frac{1}{2}},
\]
where $r_{23} = \min(r_2,r_3)$. From this expression, it is easy to observe that\\
\begin{prop}
$\phi_2$ is smooth when $|z_2|\not=|z_3|$. $\phi_2$ is merely Lipschitz continuous $(C^{0,1})$ when $|z_2|=|z_3|$. But at the point $z_2=z_3=0$, $\phi_2$ is actually $C^{1,1}$ continuous.
\end{prop}
\begin{flushright} $\Box$ \end{flushright}
In relation to the last conclusion of the proposition, it is interesting to understand the singular fibre $\phi^{-1}(0)$. This singular fiber actually has a very simple description. Under the coordinate $(z_2,z_3)$,\\
\[
\phi^{-1}(0) = \{z_2=i\bar{z}_3\}\cup \{z_2=-i\bar{z}_3\}.
\]
This is just the union of two Lagrangian planes intersecting at the origin.\\\\
Under the coordinate $(z_2,z_3)$ we can also express invariant functions\\
\[
Z_2 = z_2\left(1-\frac{r_{23}^2}{|z_2|^2}\right)^{\frac{1}{2}} = \frac{z_2}{|z_2|}(|z_2|^2 - r_{23}^2)^{\frac{1}{2}},
\]
\[
Z_3 = z_3\left(1-\frac{r_{23}^2}{|z_3|^2}\right)^{\frac{1}{2}} = \frac{z_3}{|z_3|}(|z_3|^2 - r_{23}^2)^{\frac{1}{2}}.
\]
From these expressions we can see that $Z_2$, $Z_3$ are merely $C^{0,\frac{1}{2}}$ continuous when $|z_2|=|z_3|$. But at the point $z_2=z_3=0$, $Z_2$, $Z_3$ are actually $C^{0,1}$ continuous.\\

This construction can be easily generalized to higher dimensions.\\

\subsection{The general case and local models}
{\bf Example 3.2:} In general, on ${\mathbb{C}^{m+n}}$, consider $f(z) = {\rm Re}(s)$, where\\
\[
s= \left.\left(\prod_{i=1}^mz_i\right)\right/\left(\prod_{i=m+1}^{m+n}z_i\right).
\]
Then\\
\[
\nabla f = {\rm Re}\left(\bar{s}\left(\sum_{i=1}^m \frac{z_i}{|z_i|^2} \frac{\partial}{\partial z_i} - \sum_{i=m+1}^{m+n} \frac{z_i}{|z_i|^2} \frac{\partial}{\partial z_i} \right) \right),
\]
\[
V= {\rm Re}\left(\left.\left( \sum_{i=1}^m \frac{z_i}{|z_i|^2} \frac{\partial}{\partial z_i} - \sum_{i=m+1}^{m+n} \frac{z_i}{|z_i|^2} \frac{\partial}{\partial z_i} \right) \right/s\sum_{i=1}^{m+n} \frac{1}{|z_i|^2} \right).
\]\\
Vector field $V$ deforms the hypersurface $\prod_{i=1}^mz_i=0$ to $\prod_{i=m+1}^{m+n}z_i=0$. Notice that from the above expression of $V$ we can see that on ${\mathbb{C}^{m+n}}$ the invariant functions are\\
\[
\rho_{ij} = |z_i|^2-|z_j|^2,\ \ {\rm for}\ 1\leq i,j \leq m\ {\rm or}\ m+1 \leq i,j \leq m+n.
\]
\[
\rho_{ij} = |z_i|^2+|z_j|^2,\ \ {\rm for}\ 1\leq i \leq m\ {\rm and}\ m+1 \leq j \leq m+n.
\]\\
As in the previous example, $V$ leaves $g(z)={\rm Im}(s)$ invariant. Restricted to the real hypersurface $g(z)={\rm Im}(s)=0$, $V$ can be simplified as\\
\[
V= {\rm Re}\left.\left( \sum_{i=1}^m \frac{z_i}{|z_i|^2} \frac{\partial}{\partial z_i} - \sum_{i=m+1}^{m+n} \frac{z_i}{|z_i|^2} \frac{\partial}{\partial z_i} \right) \right/s\sum_{i=1}^{m+n} \frac{1}{|z_i|^2}.
\]\\
Using this expression of $V$, on $\{g=0\}$ we have the invariant functions\\
\[
{\rm Im}(\log z_i)=-i(\log z_i - \log |z_i|),\ \frac{z_i}{|z_i|},\ \ {\rm for}\ \ 1\leq i \leq m+n.
\]
For $1\leq l \leq m$, on the hyperplane $z_l=0$, we have coordinate $(z_1, \cdots, \hat{z}_l, \cdots, z_{m+n})$. Coordinate functions can be invariantly extended to\\
\begin{eqnarray*}
Z_i &=& z_i\left(1+\left|\frac{z_l}{z_i}\right|^2\right)^{\frac{1}{2}} = \frac{z_i}{|z_i|}(|z_i|^2 + |z_l|^2)^{\frac{1}{2}} \ \ {\rm when}\ |z_l|=\min_{1\leq k \leq m}|z_k|,\ m+1\leq i \leq m+n,\\
Z_i &=& z_i\left(1-\left|\frac{z_l}{z_i}\right|^2\right)^{\frac{1}{2}} = \frac{z_i}{|z_i|}(|z_i|^2 - |z_l|^2)^{\frac{1}{2}} \ \ {\rm when}\ |z_l|= \min_{1\leq k \leq m}|z_k|,\ 1\leq i \leq m.
\end{eqnarray*}
\begin{flushright} $\Box$ \end{flushright}
In example 3.2, the special cases of $n=0,1$ exactly correspond to the two general local models for deforming the large complex limit to smooth Calabi-Yau hypersurfaces in toric varieties. (The $n=0$ local model was already discussed in \cite{lag1}.) We now briefly describe these local models.\\

{\bf Local model I:} On ${\mathbb{C}^{n+1}}$, consider $f(z) = {\rm Re}(s)$, where $s=\frac{\prod_{i=1}^nz_i}{z_0}$. Then\\
\[
\nabla f = {\rm Re}\left(\bar{s}\left(\sum_{i=1}^n \frac{z_i}{|z_i|^2} \frac{\partial}{\partial z_i} - \frac{z_0}{|z_0|^2} \frac{\partial}{\partial z_0} \right) \right),
\]
\[
V= {\rm Re}\left(\left.\left( \sum_{i=1}^n \frac{z_i}{|z_i|^2} \frac{\partial}{\partial z_i} - \frac{z_0}{|z_0|^2} \frac{\partial}{\partial z_0} \right) \right/s\sum_{i=0}^n \frac{1}{|z_i|^2} \right).
\]\\
Vector field $V$ deforms the hypersurface $\prod_{i=1}^nz_i=0$ to the hypersurface $z_0=0$. The invariant functions of $V$ on ${\mathbb{C}^{n+1}}$ are\\
\[
\rho_{ij} = |z_i|^2-|z_j|^2,\ \rho_{0i} = |z_0|^2+|z_i|^2\ {\rm for}\ i,j\geq 1.
\]\\
As we know, $V$ leaves $g(z)={\rm Im}(s)$ invariant. The invariant functions of $V$ on $\{g=0\}$ are\\
\[
{\rm Im}(\log z_i)=-i(\log z_i - \log |z_i|),\ \frac{z_i}{|z_i|},\ \ {\rm for}\ \ 0\leq i \leq n.
\]
For $1\leq m \leq n$, on the hyperplane $z_m=0$, we have coordinate $(z_0, \cdots, \hat{z}_m, \cdots, z_n)$. Coordinate functions can be invariantly extended to\\
\begin{eqnarray*}
Z_0 &=& z_0\left(1+\left|\frac{z_m}{z_0}\right|^2\right)^{\frac{1}{2}} = \frac{z_0}{|z_0|}(|z_0|^2 + |z_m|^2)^{\frac{1}{2}} \ \ {\rm when}\ |z_m|=\min_{1\leq k \leq n}|z_k|,\\
Z_i &=& z_i\left(1-\left|\frac{z_m}{z_i}\right|^2\right)^{\frac{1}{2}} = \frac{z_i}{|z_i|}(|z_i|^2 - |z_m|^2)^{\frac{1}{2}} \ \ {\rm when}\ |z_m|= \min_{1\leq k \leq n}|z_k|,\ 1\leq i \leq n.
\end{eqnarray*}
\begin{flushright} $\Box$ \end{flushright}
{\bf Local model II:} (Example in the subsection 3.3 of \cite{lag1}) On $\mathbb{C}^{{n}}$, consider $f(z) = {\rm Re}(s)$, where $s=\displaystyle \prod_{i=1}^nz_i$. Then\\
\[
\nabla f = {\rm Re}\left(\bar{s}\sum_{i=1}^n \frac{z_i}{|z_i|^2} \frac{\partial}{\partial z_i} \right),
\]
\[
V= {\rm Re}\left(\left.\left( \sum_{i=1}^n \frac{z_i}{|z_i|^2} \frac{\partial}{\partial z_i}\right) \right/s\sum_{i=1}^n \frac{1}{|z_i|^2} \right).
\]\\
Vector field $V$ deforms the hypersurface $\prod_{i=1}^nz_i=0$ to the hypersurface $\prod_{i=1}^nz_i= c >0$. The invariant functions of $V$ on $\mathbb{C}^{{n}}$ are\\
\[
\rho_{ij} = |z_i|^2-|z_j|^2\ {\rm for}\ 1\leq i,j\leq n.
\]\\
As we know, $V$ leaves $g(z)={\rm Im}(s)$ invariant. The invariant functions of $V$ on $\{g=0\}$ are\\
\[
{\rm Im}(\log z_i)=-i(\log z_i - \log |z_i|),\ \frac{z_i}{|z_i|},\ \ {\rm for}\ \ 1\leq i \leq n.
\]
For $1\leq m \leq n$, on the hyperplane $z_m=0$, we have coordinate $(z_1, \cdots, \hat{z}_m, \cdots, z_n)$. Coordinate functions can be invariantly extended to\\
\[
Z_i = z_i\left(1-\left|\frac{z_m}{z_i}\right|^2\right)^{\frac{1}{2}} = \frac{z_i}{|z_i|}(|z_i|^2 - |z_m|^2)^{\frac{1}{2}} \ \ {\rm when}\ |z_m|= \min_{1\leq k \leq n}|z_k|,\ 1\leq i \leq n.
\]
\begin{flushright} $\Box$ \end{flushright}
In example 3.2, the remaining cases of $m,n>1$ correspond to local models of some non-generic situations of singular hypersufaces in toric varieties and are not necessary for our application to the case of generic smooth hypersurfaces in toric varieties.\\

A crucial step in our construction is to show that the gradient flow will fix $X_\infty\cap X_\psi$ and flow $X_\infty\backslash (X_\infty\cap X_\psi)$ to $X_\psi\backslash (X_\infty\cap X_\psi)$. This fact for the local models can be shown using the following proposition and its proof. The general case will be proved in section 5 using perturbation argument.\\
\begin{prop}
On ${\mathbb{C}^{m+n}}$, consider $f(z) = {\rm Re}(s)$, where\\
\[
s= \left.\left(\prod_{i=1}^mz_i\right)\right/\left(\prod_{i=m+1}^{m+n}z_i\right).
\]
Then when $n,m>0$, no non-constant solution curve of the direction field $\nabla f$ will approach the origin.\\
\end{prop}
{\bf Proof:} Assume a solution curve $\phi(t)$ approaches the origin, then the invariant functions\\
\[
\rho_{ij} = |z_i|^2+|z_j|^2,\ \ {\rm for}\ 1\leq i \leq m\ {\rm and}\ m+1 \leq j \leq m+n
\]\\
will vanish along the solution curve $\phi(t)$. This implies that\\
\[
z_i(\phi(t))=0\ \ {\rm for}\ 1\leq i \leq m+n.
\]\\
Namely $\phi(t)= 0$ is constant solution.
\begin{flushright} $\Box$ \end{flushright}

\se{Hyperbolic homogeneous vector fields}
As is well known, the behavior of dynamical systems near singular points (even non-degenerate ones) of the vector fields is very sensitive to local perturbation. For singular points with higher degeneracy, such dependence is even more delicate. In our situation, the gradient vector field is highly degenerate and have poles. Perturbation stability of such vector fields is crucial for us to ensure our dynamical systems behave according to the local models discussed in the previous section near the singular points of the gradient vector fields.\\

We will discuss perturbation stability in this and the next sections. In this section, we will start with the general discussion of homogeneous vector fields. The key property that ensures the stability of our system is certain hyperbolic property of our vector field. We will explain such hyperbolicity in detail. The general discussion of stability of hyperbolic homogeneous vector fields in this section will ensure the stability of the unstable (stable) cones of the our vector field near singularities. In the next section, using more detailed information of our vector field, we will be able to prove stronger stability for our particular vector field, further ruling out the possibility that solution curves may spiral down to the singular points.\\

For the discussion of homogeneous vector fields, let us use polar coordinate $(r,\theta)$. Let $V=V_0+V_1$, where $V_0$ is the degree $d$ homogeneous term of $V$ and $V_1$ is the higher order term. Then

\[
V_0 = r^{d-1}v_0(\theta)\frac{\partial}{\partial \theta} + r^dp_0(\theta)\frac{\partial}{\partial r},
\]

where $v_0 = v_0(\theta)\displaystyle\frac{\partial}{\partial \theta} = \displaystyle\sum_{i=1}^{n-1}v_0^i(\theta)\frac{\partial}{\partial \theta_i}$ is the expression of $v_0 = V_0|_{S^{n-1}}$ under the coordinate $\theta =(\theta_1,\cdots,\theta_{n-1})$ of $S^{n-1}$. Similarly,

\[
V = r^{d-1}v(r,\theta)\frac{\partial}{\partial \theta} + r^dp(r,\theta)\frac{\partial}{\partial r},
\]

where $v(r,\theta) = v_0(\theta) + rv_1(r,\theta)$ and $p(r,\theta) = p_0(\theta) + rp_1(r,\theta)$. Use $r$ as parameter, then $\theta(r)$ satisfies

\[
r\frac{d\theta}{dr} = \frac{v(r,\theta)}{p(r,\theta)}.
\]

We say a solution curve $C:\ (r(t),\theta(t))$ comes out of (or into) the origin, if $C$ is a continuous curve in $(r,\theta)$ space with one end point being $(0,\theta_0)$. (For simplicity of notation, we will express such fact by $\displaystyle \lim_{r\rightarrow 0} \theta(r) = \theta_0$, although $\theta(r)$ is usually not a single valued function on $r$.) It is easy to see such curve will be tangent to the ray $\{\theta = \theta_0\}$ at the origin. There is another way a solution curve $C$ approaches the origin, which we refer to as $C$ spirals down to the origin, when $\displaystyle \lim_{r\rightarrow 0} \theta(r)$ does not exist. For the general discussion in this section, we will only deal with solution curves coming out of (or into) the origin. The more subtle solution curves spiraling down to the origin will be ruled out in the next section for the specific vector fields we need.\\

\begin{de}
For a vector field $V$ defined near the origin, the unstable (stable) variety $S^+$ ($S^-$) of $V$ is defined to be the union of solution curves coming out of (into) the origin.
\end{de}

\begin{prop}
When $V_0$ is homogeneous, its unstable (stable) variety $S^+$ ($S^-$) is a cone, and is called the unstable (stable) cone.
\end{prop}
{\bf Proof:} Assume

\[
V_0 = r^{d-1}v_0(\theta)\frac{\partial}{\partial \theta} + r^dp_0(\theta)\frac{\partial}{\partial r}.
\]

Use $r$ as parameter, then $\theta(r)$ satisfies

\[
r\frac{d\theta}{dr} = \frac{v_0(\theta)}{p_0(\theta)}.
\]

It is easy to verify that if $(r,\theta(r))$ is a solution curve coming out of (into) the origin, then $(cr,\theta(r))$ is also a solution curve coming out of (into) the origin for any constant $c$. Therefore $S^+$ ($S^-$) is a cone.
\begin{flushright} $\Box$ \end{flushright}
The behavior of a general homogeneous vector field $V_0$ could be quite messy. To make the discussion meaningful, in this paper, we will concentrate on the set of ``good" homogeneous vector fields satisfying certain reasonable constraints that will include most interesting homogeneous vector fields. More precisely, when $p_0(\theta) \not\equiv 0$, we assume that the domain of the function $v_0(\theta)/p_0(\theta)$ is open and dense in $S^{n-1}$, and $v_0(\theta)/p_0(\theta)$ can be extended to a smooth function on $\{\theta |\displaystyle\lim_{\theta'\rightarrow \theta}|v_0(\theta')/p_0(\theta')|\not= + \infty\}$. It is straightforward to check that the homogeneous vector field in Example 3.2 of section 3 and consequently the homogeneous vector fields in local model I and local model II in section 3 satisfy our restriction here. Define

\[
Z = \{\theta\in S^{n-1}|\lim_{\theta'\rightarrow \theta}v_0(\theta')/p_0(\theta')=0\}.
\]

\[
Z_+ = \{\theta\in Z|p_0\geq 0\ {\rm (not\ identically\ 0)\ in\ a\ small\ neighborhood\ of\ }\theta\}.
\]

Similarly, we can define $Z_-$, and $Z_0= Z\backslash(Z_+\cup Z_-)$. Readers may find explicit computations of $Z_\pm$ for the gradient vector field of local model II in the example at the end of this section.\\

{\bf Remark:} The properties of the homogeneous vector field $V_0$ we will need are that $v_0(\theta)/p_0(\theta)$ is smooth in a small neighborhood $U_Z \subset S^{n-1}$ of $Z$ and $|v_0(\theta)/p_0(\theta)|$ is bounded from below by a positive constant in $S^{n-1}\backslash U_Z$. It is easy to see that our ``good" homogeneous vector fields will exhibit such properties. For our specific application, $Z$ is a smooth manifold and $p_0$ is non-vanishing when restricted to $Z$.\\

Heuristically speaking, for $\theta_0\not\in Z$, near $(0,\theta_0)$ the spherical component of $V$ will dominate the radial part of $V$. So it will not be possible for $\log r$ to approach $-\infty$ (or $r$ to approach $0$) while $\theta$ approaches $\theta_0$. More precisely, {\it no solution curve will approach $(0,\theta_0)$ for $\theta_0\not\in Z$}. There are general arguments to show this fact under certain constraints that will be satisfied by the particular vector fields for our application. As a corrollary, {\it if $Z=\emptyset$, then no solution curve of original homogeneous system will approach the origin.}\\

Instead of presenting the general argument here, we will give an alternative argument in the next section which is special and more elegant for our particular vector fields. The special argument will prove stronger results. In particular, it will rule out solution curves that will spiral down and eventually aproach the origin.\\

Let us now turn our attention to $Z$. We are mostly interested in the case of non-degenerate homogeneous vector fields defined as follows.\\

\begin{de}
\label{ca}
A homogeneous vector field $V_0$ is called {\bf non-degenerate} if $Z_0=\emptyset$, $Z_+$ and $Z_-$ are disjoint smooth submanifolds in $S^{n-1}$, and $\frac{v_0}{p_0}$ induces non-degenerate bilinear form\\
\[
\langle n_1,n_2\rangle = n_1\left(\frac{v_0}{p_0}, n_2\right)
\]
on the normal bundle of $Z_+$ and $Z_-$. A non-degenerate homogeneous vector field $V_0$ is called {\bf hyperbolic} if the real part of the eigenvalues of the bilinear form $\langle n_1,n_2\rangle$ with respect to the metric pairing $(n_1,n_2)$ are all strictly negative.\\
\end{de}

{\bf Remark:} In the above definition, for $\theta \in Z_+$, $n_1,n_2\in (N_{Z_+})_\theta$, we may extend $n_1,n_2$ to vector fields on $S^{n-1}$. Then $(\frac{v_0}{p_0}, n_2)$ is a function on $S^{n-1}$. Its derivative in the $n_1$ direction $n_1(\frac{v_0}{p_0}, n_2)$ is also a function on $S^{n-1}$. In the above definition, $\langle n_1,n_2\rangle$ is defined as the value of $n_1(\frac{v_0}{p_0}, n_2)$ at $\theta$. It is easy to see that this definition is independent of how $n_1,n_2$ are extended.\\

We will also introduce the concept of geometric hyperbolicity that is more natural geometrically. Let $U_{Z_-}$ be a small tubular neighborhood of $Z_-$. $U_{Z_-}$ can be naturally identified with a small tubular neighborhood of the zero section of the normal bundle $N_{Z_-}$. We have

\[
\begin{array}{ccc}
i:\ U_{Z_-}&\hookrightarrow& N_{Z_-} \\
&&\\
\ \ \downarrow\pi&&\\
&&\\
Z_- &&\\
\end{array}
\]

According to this diagram, we can introduce the position vector field $e_-$ on $U_{Z_-}$. For $\theta\in \pi^{-1}(\pi(\theta))\hookrightarrow (N_{Z_-})_{\pi(\theta)}$, the tangent space $T_{\theta}\pi^{-1}(\pi(\theta))$ is naturally identified with $(N_{Z_-})_{\pi(\theta)}$. Under this identification, $e_-(\theta)$ is defined to correspond to $i(\theta)$. Let $U_{Z_-}^* = U_{Z_-}\backslash Z_-$ ($U_{Z_+}^* = U_{Z_+}\backslash Z_+$).

\begin{de}
\label{cb}
A homogeneous vector field $V_0$ is called {\bf geometrically hyperbolic} if $Z_0=\emptyset$, $Z_+$ and $Z_-$ are disjoint smooth submanifolds in $S^{n-1}$ and for $\theta\in U_{Z_-}^*$ ($\theta\in U_{Z_+}^*$), we have

\[
(e_-(\theta), v_0(\theta))_{S^{n-1}}>0\ \ ((e_+(\theta), v_0(\theta))_{S^{n-1}}<0).
\]

\end{de}
Let $S_0^+$ ($S_0^-$) denote the cone over $Z_+$ ($Z_-$), and $U_{S_0^+}$ ($U_{S_0^-})$ denote the cone over $U_{Z_+}$ ($U_{Z_-}$). Then $U_{S_0^+}$ ($U_{S_0^-})$ is a small cone neighborhood of the cone $S_0^+$ ($S_0^-$).\\

\begin{theorem}
Assume that $V_0$ is a geometrically hyperbolic homogeneous vector field of degree $d$ and $S^+$ ($S^-$) is the unstable (stable) cone of $V_0$. Then $S^+ \cap U_{S_0^+} = S_0^+$ ($S^- \cap U_{S_0^-} = S_0^-$) is a cone over $Z_+$ ($Z_-$).
\end{theorem}
{\bf Proof:} Clearly $S_0^+ \subset S^+$ ($S_0^- \subset S^-$). Without loss of generality, we will concentrate on $S^-$. Recall $U_{S_0^-}$ is a cone neighborhood of $S_0^-$. The geometric hyperbolicity condition on $V_0$ is equivalent to the condition that $V_0$ is always pointing toward outside of $U_{S_0^-}$ when restricted to $\partial U_{S_0^-}$.\\

Assume that $x(t)$ ($t\leq 0$, $x(0)=0$) inside $U_{S_0^-}$ is a solution curve of $V_0$ coming into the origin from the angle $\theta_0$. Then $cx(t)$ are also solutions for any constant $c$. Recall that $v_0(\theta)/p_0(\theta)$ is assumed to be smooth near $Z_-$. Consequently, the ray $\{\theta = \theta_0\} = \displaystyle \lim_{c\rightarrow +\infty} cx(t)$ is a solution and belongs to $S_0^-$. In particular, $x(t)$ is tangent to $S_0^-$ at the origin. If $x(t)$ is not entirely in $S_0^-$, by suitably shrinking $U_{S_0^-}$, we may assume that $x(t)$ intersects $\partial U_{S_0^-}$ at $x(t_0)$. Since $V_0$ is pointing toward outside of $U_{S_0^-}$ at $x(t_0)$, $x(t)$ will be outside of $U_{S_0^-}$ for all $t>t_0$ because no solution curve can enter $\partial U_{S_0^-}$ from outside. This contradicts with the conclusion that $x(t)$ is tangent to $S_0^-$ at the origin. Therefore, $x(t)$ is entirely in $S_0^-$.
\begin{flushright} $\Box$ \end{flushright}
It will be interesting to relate the concepts of hyperbolicity and geometric hyperbolicity. We will start with a linear algebra lemma.\\

\begin{lm}
\label{ci}
Assume that $A$ is a symmetric positive definite real matrix and $B$ is an anti-symmetric real matrix. Then the real parts of eigenvalues of $A+B$ are all positive.
\end{lm}
{\bf Proof:} Assume $(A+B)v = \lambda v$. Then $v^*(A+B)v = v^*\lambda v = \lambda |v|^2$. (Recall $v^* = \bar{v}^T$.) It is easy to verify that $v^*Av$ is real and $v^*Bv$ is purely imaginary. Therefore ${\rm Re}(\lambda) |v|^2 = v^*Av> 0$. Therefore ${\rm Re}(\lambda) > 0$.
\begin{flushright} $\Box$ \end{flushright}
\begin{lm}
\label{ck}
The symmetrization of the bilinear form $\langle n_1,n_2\rangle$ introduced in definition \ref{ca} is exactly the Hessian form of $\frac{1}{2}(e_-(\theta), \frac{v_0}{p_0}(\theta))_{S^{n-1}}$ along the fibre $\pi^{-1}(\pi(\theta))\hookrightarrow (N_{Z_-})_{\pi(\theta)}$ at zero.
\end{lm}
{\bf Proof:} Recall that the Hessian for a function $f$ with respect to the metric is defined as

\[
{\rm Hess}(f)(n_1,n_2) = n_1(n_2(f)) - \nabla_{n_1}n_2 (f).
\]

At a point where $f$ vanishes to order 2, ${\rm Hess}(f)(n_1,n_2) = n_1(n_2(f))$. Notice that $f = (e_-(\theta), \frac{v_0}{p_0}(\theta))_{S^{n-1}}$ vanishes to order 2 on $\pi(\theta)$. For $n_1,n_2 \in (N_{Z_-})_{\pi(\theta)} \cong T_{\pi(\theta)}\pi^{-1}(\pi(\theta))$ we have

\begin{eqnarray*}
{\rm Hess}(f)|_{\pi(\theta)}(n_1,n_2) &=& (\nabla_{n_1}e_-(\theta), \nabla_{n_2}\frac{v_0}{p_0}(\theta))_{S^{n-1}} + (\nabla_{n_2}e_-(\theta), \nabla_{n_1}\frac{v_0}{p_0}(\theta))_{S^{n-1}}\\
&=& (n_1, \nabla_{n_2}\frac{v_0}{p_0}(\theta))_{S^{n-1}} + (n_2, \nabla_{n_1}\frac{v_0}{p_0}(\theta))_{S^{n-1}}\\
&=& \langle n_2,n_1\rangle + \langle n_1,n_2\rangle.
\end{eqnarray*}
\begin{flushright} $\Box$ \end{flushright}
\begin{prop}
\label{cj}
Assume that $V_0$ is a geometrically hyperbolic homogeneous vector field of degree $d$. If the symmetrization of the bilinear form $\langle n_1,n_2\rangle$ in definition \ref{ca} (equivalently the Hessian form of $(e_{\pm}(\theta), \frac{v_0}{p_0}(\theta))_{S^{n-1}}$ along the fibre $\pi^{-1}(\pi(\theta))\hookrightarrow (N_{Z_{\pm}})_{\pi(\theta)}$ at zero) is non-degenerate, then $V_0$ is hyperbolic.
\end{prop}
{\bf Proof:} Without loss of generality, we will concentrate on $Z_-$. According to lemma \ref{ck}, the symmetrization of the bilinear form $\langle n_1,n_2\rangle$ equals to the Hessian form of $\frac{1}{2}(e_-(\theta), \frac{v_0}{p_0}(\theta))_{S^{n-1}}$ along the fibre $\pi^{-1}(\pi(\theta))\hookrightarrow (N_{Z_-})_{\pi(\theta)}$ at zero. $V_0$ being geometrically hyperbolic implies that $(e_-(\theta), \frac{v_0}{p_0}(\theta))_{S^{n-1}}$ is negative. This together with the assumption that the Hessian of $(e_-(\theta), \frac{v_0}{p_0}(\theta))_{S^{n-1}}$ is non-degenerate imply that the Hessian of $(e_-(\theta), \frac{v_0}{p_0}(\theta))_{S^{n-1}}$ is negative definite. Hence the symmetrization of the bilinear form $\langle n_1,n_2\rangle$ is negative definite. By lemma \ref{ci}, we conclude that the real parts of the eigenvalues of the bilinear form $\langle n_1,n_2\rangle$ with respect to the metric pairing $(n_1,n_2)$ are all strictly negative. Therefore $V_0$ is hyperbolic.
\begin{flushright} $\Box$ \end{flushright}
\begin{co}
\label{cl}
Assume that $V_0$ is a homogeneous vector field that is equivalent to a gradient vector field $\nabla f$ as a direction field. Then $V_0$ is hyperbolic if and only if $V_0$ is non-degenerate and geometrically hyperbolic.
\end{co}
{\bf Proof:} Under the assumption of the corollary, the bilinear form $\langle n_1,n_2\rangle$ is symmetric. By proposition \ref{cj}, the conclusion is immediate.
\begin{flushright} $\Box$ \end{flushright}
When $V_0$ is equivalent to a gradient vector field $\nabla f$ as a direction field, hyperbolicity is easier to check. Assume that $f$ is homogeneous of degree $d+1$. Then $f(r,\theta) = r^{d+1} \hat{f}(\theta)$, where $\hat{f}(\theta) = f(1,\theta)$. Under the standard flat metric, we have

\[
\nabla f = (d+1) r^d \hat{f}(\theta)\frac{\partial}{\partial r} + r^{d-1} \nabla_\theta \hat{f}(\theta)\frac{\partial}{\partial \theta}.
\]

Namely,

\[
v_0(\theta) = \nabla_\theta \hat{f}(\theta),\ \ p_0(\theta) = (d+1) \hat{f}(\theta),\ \ \frac{v_0}{p_0} = \frac{\nabla_\theta \hat{f}}{(d+1) \hat{f}}.
\]

\begin{prop}
\label{cc}
A non-degenerate homogeneous vector field $V_0=\rho \nabla f$ is hyperbolic if

\[
e_-(\theta)(\hat{f})>0\ \ (e_+(\theta)(\hat{f})<0)
\]

for $\theta\in U_{Z_-}^*$ ($\theta\in U_{Z_+}^*$), where $\hat{f}(\theta) = f(1,\theta)$.
\end{prop}

{\bf Proof:} Since $v_0(\theta) = \nabla_\theta \hat{f}(\theta)$, we have

\[
(e_\pm(\theta), v_0(\theta))_{S^{n-1}} = e_\pm(\theta)(\hat{f}).
\]

By the definition of geometric hyperbolicity and corollary \ref{cl}, the conclusion is immediate.
\begin{flushright} $\Box$ \end{flushright}
\begin{co}
\label{cd}
A non-degenerate homogeneous vector field $V_0=\rho \nabla f$ is hyperbolic if $\hat{f} = f|_{S^{n-1}}$ achieves its maximum in $Z_+$ and achieves its minimum in $Z_-$.
\end{co}

{\bf Proof:} Since $\displaystyle \frac{v_0}{p_0} = \frac{\nabla_\theta \hat{f}}{(d+1) \hat{f}}$ and both $\nabla_\theta \hat{f}$ and $\hat{f}$ vanish when restrict to $Z_+$ or $Z_-$. It is easy to verify that along $Z_+$ or $Z_-$ the bilinear form

\[
\langle n_1,n_2\rangle = \frac{{\rm Hess}(\hat{f})(n_1,n_2)}{(d+1) \hat{f}}.
\]

By the definition of hyperbolicity, the conclusion is immediate.
\begin{flushright} $\Box$ \end{flushright}
{\bf Remark:} $V_0=\rho \nabla f$ is non-degenerate at $Z_+$ and $Z_-$ if and only if $Z_+$ and $Z_-$ are non-degenerate critical sets of $f|_{S^{n-1}}$. Proposition \ref{cc} and corollary \ref{cd} for gradient vector field are still meaningful even if the non-degeneracy condition for the bilinear form in definition \ref{ca} is removed.\\

The perturbation stability of unstable (stable) cones can be reduced to the following singular ODE problems.

\begin{prop}
For the initial value problem\\
\[
\left\{
\begin{array}{l}
\displaystyle\frac{d\theta_i}{dr} = -\displaystyle\frac{\lambda_i}{r}\theta_i + h_i(r,\theta),\ \ {\rm for}\ i=1,\cdots,m.\\\\
\theta(0)=0
\end{array}
\right.
\]\\
with $\lambda_i \geq 0$ for $i=1,\cdots,m$, there exists a unique solution. When $-1 <\lambda_i < 0$ for $i=1,\cdots,m$, there exists a unique solution satisfying $\theta(r) = O(r)$.\\
\end{prop}
{\bf Proof:} Under the condition $\lambda_i \geq 0$
(or $-1 <\lambda_i < 0$ and $\theta(r) = O(r)$) for $i=1,\cdots,m$, the equation can be reduced to\\
\[
\frac{dr^{\lambda_i}\theta_i}{dr} = r^{\lambda_i} h_i(r,\theta).
\]\\
Hence\\
\[
\theta_i(r) = r^{-\lambda_i}\int_0^r t^{\lambda_i} h_i(t,\theta)dt.
\]\\
In a neighborhood of the origin, assume that\\
\[
\left|\frac{\partial h}{\partial \theta}\right| \leq M.
\]\\
Then we have\\
\[
\delta\theta_i(r) \leq \frac{M}{1+\lambda_i}r|\delta\theta|.
\]\\
The equations implies that\\
\[
\frac{d\theta_i}{dr}(0) = \frac{1}{1+\lambda_i} h_i(0,0).
\]\\
Let\\
\[
B = \{\theta(r) \in C[0,r_0]|\theta(0)=0,\ \frac{d\theta_i}{dr}(0) = \frac{1}{1+\lambda_i} h_i(0,0)\}.
\]\\
Consider the operator $T: B\rightarrow B$ defined as\\
\[
(T\theta)_i(r) = r^{-\lambda_i}\int_0^r t^{\lambda_i} h_i(t,\theta)dt.
\]\\
Solutions of our ordinary differential equation exactly correspond to fixed points of $T$. Clearly\\
\[
\|\delta(T\theta)\| \leq Mr_0 \max_{1\leq i\leq n}\left(\frac{1}{1+\lambda_i}\right)\|\delta\theta\|.
\]\\
Choose $r_0$ such that\\
\[
Mr_0 \max_{1\leq i\leq n}\left(\frac{1}{1+\lambda_i}\right) = \frac{1}{2}.
\]\\
Then $T$ is a contraction map and $\|\delta(T\theta)\| \leq \frac{1}{2}\|\delta\theta\|$. By contraction mapping theorem, $T$ has one unique fixed point in $B$, which is the unique solution of our equation.
\begin{flushright} $\Box$ \end{flushright}
In general, consider $\theta = (\theta_1,\cdots, \theta_{m+n}) = (\alpha, \beta)$ where $\alpha = (\theta_1,\cdots, \theta_m)$, $\beta = (\theta_{m+1},\cdots, \theta_{m+n})$. \\
\begin{prop}
\label{ce}
Assume that the real parts of all the eigenvalues of $H(0,\beta)$ are strictly positive and $H(\theta),h_1(r,\theta),h_2(r,\theta)$ are smooth. Then the initial value problem
\[
\left\{
\begin{array}{l}
\displaystyle\frac{d\alpha}{dr} = -\frac{1}{r}H(\theta)\alpha + h_1(r,\theta)\\\\
\displaystyle\frac{d\beta}{dr} = h_2(r,\theta)\\\\
\theta(0)=(0,\beta_0)
\end{array}
\right.
\]

has a unique solution that  depends on the initial value continuously.\\
\end{prop}
{\bf Proof:} The first equation can be rewritten as

\[
\frac{d\alpha}{dr} = -\frac{1}{r}H(0,\beta_0)\alpha + \frac{1}{r}H_1(\theta) + h_1(r,\theta),
\]

where $H_1(\theta)$ is smooth and vanishes to the second order at $(0,\beta_0)$. Denote $H_{\beta_0} = H(0,\beta_0)$, the equation can be reduced to\\
\[
\frac{dr^{H_{\beta_0}}\alpha}{dr} = r^{H_{\beta_0}} \left(\frac{1}{r}H_1(\theta) + h_1(r,\theta)\right).
\]

Hence

\[
\alpha(r) = r^{-H_{\beta_0}}\int_0^r t^{H_{\beta_0}} \left(\frac{1}{t}H_1(\theta) + h_1(t,\theta)\right)dt.
\]

Let

\[
B = \left\{\theta(r) \in C[0,r_0]\left| \theta(0)=(0,\beta_0),\ |\theta(r)-\theta(0)|\leq Mr\right.\right\},
\]

where $M>0$ is a constant. Consider the operator $T$ defined as

\[
(T\theta)(r) = \left( r^{-H_{\beta_0}}\int_0^r t^{H_{\beta_0}} \left(\frac{1}{t}H_1(\theta) + h_1(t,\theta)\right)dt, \beta_0 + \int_0^r h_2(t,\theta)dt\right).
\]

For $r_0$ small and $M$ suitably chosen depending on the bounds of $h_1,h_2$ and eigenvalues of $H_{\beta_0}$, it is easy to verify that $T: B\rightarrow B$. Solutions of our ordinary differential equation exactly correspond to fixed points of $T$. Since $H,h_1,h_2$ are smooth, it is easy to derive that\\
\[
\|\delta(T\theta)\| \leq \left(CMr_0\|(I + H_{\beta_0})^{-1}\|\right)\|\delta\theta\|.
\]\\
Choose $r_0$ such that

\[
CMr_0\|(I + H_{\beta_0})^{-1}\| \leq \frac{1}{2}.
\]

Then $T$ is a contraction map and $\|\delta(T\theta)\| \leq \frac{1}{2}\|\delta\theta\|$. By contraction mapping theorem, $T$ has a unique fixed point in $B$, which is the unique solution of our equation.\\

Similarly, for solutions $\theta(r)$ depending on initial value $\beta_0$, with the same choice of $r_0$, we have

\[
\|\delta\theta\| \leq \frac{1}{2}\|\delta\theta\| + (C_1 + C_2|\log r_0|)r_0|\delta \beta_0| + |\delta \beta_0|.
\]

Choose smaller $r_0$ if necessary, we have

\[
\|\delta\theta\| \leq 3|\delta \beta_0|.
\]

This shows the continuous dependence of solutions on the initial values.\\

According to the integral version of the equation, we can easily derive that

\[
\frac{d\alpha}{dr}(0) = (I + H_{\beta_0})^{-1} h_1(0,0,\beta_0).
\]

\[
\frac{d\beta}{dr}(0) = h_2(0,0,\beta_0).
\]
\begin{flushright} $\Box$ \end{flushright}
\begin{prop}
\label{ch}
With the same conditions as in proposition \ref{ce}, the system of equations in proposition \ref{ce} will have no additional solution $\theta(r) = (\alpha(r),\beta(r))$, $r\in(0,r_0]$ satisfying $|\alpha(r)| < \epsilon$ for $r\in(0,r_0]$ and $\epsilon >0$ sufficiently small. In particular, there is no solution with initial value $\theta(0) = (\alpha_0,\beta_0)$ when $\alpha_0\not=0$ is small.
\end{prop}

{\bf Proof:} For $r_1<r$ in $(0,r_0]$, by the integral version of the equation for $\beta$, we have

\[
\beta(r_1) = \beta(r) - \int_{r_1}^{r} h_2(t,\theta)dt.
\]

Since $h_2(t,\theta)$ is bounded, the limit

\[
\beta_0 = \lim_{r_1\rightarrow 0}\beta(r_1) = \beta(r) - \int_{0}^{r} h_2(t,\theta)dt
\]

exists and

\[
\beta(r) = \beta_0 + \int_{0}^{r} h_2(t,\theta)dt
\]

This equation implies the estimate $|\beta(r) - \beta_0| \leq Cr$.\\

By the integral version of the equation for $\alpha$, we have

\[
\alpha(r) - \left(\frac{r_1}{r}\right)^{H_{\beta_0}}\alpha(r_1) = r^{-H_{\beta_0}}\int_{r_1}^r t^{H_{\beta_0}} \left(\frac{1}{t}H_1(\theta) + h_1(t,\theta)\right)dt.
\]

Since $|\alpha(r_1)| < \epsilon$ and the eigenvalues of $H_{\beta_0}$ are strictly positive, take limit when $r_1\rightarrow 0$, we have

\[
\alpha(r) = r^{-H_{\beta_0}}\int_0^r t^{H_{\beta_0}} \left(\frac{1}{t}H_1(\theta) + h_1(t,\theta)\right)dt.
\]

Since $H_1(\alpha,\beta)$ is of second order on $(\alpha,\beta-\beta_0)$ and $h_1(t,\theta)$ is bounded, we have the estimate

\[
|\alpha(r)| \leq C_1 \sup_{t\in (0,r]}|\alpha(t)|^2 + C_2 r.
\]

Hence

\[
\sup_{t\in (0,r]}|\alpha(t)| \leq C_1 \sup_{t\in (0,r]}|\alpha(t)|^2 + C_2 r.
\]

At this point if we take $\epsilon = \frac{1}{2C_2}$, we have

\[
|\alpha(r)| \leq C r.
\]

Consequently, $\displaystyle\lim_{r\rightarrow 0}\alpha(r) = 0$. Namely $\theta(r) = (\alpha(r),\beta(r))$ is a solution in proposition \ref{ce}.
\begin{flushright} $\Box$ \end{flushright}
\begin{theorem}
\label{cf}
Assume that $V_0$ is a hyperbolic homogeneous vector field of degree $d$. $V=V_0+V_1$, where $|V_1| \leq Cr|V_0|$. And $S^+$ ($S^-$) is the unstable (stable) variety of $V$. Then $S^+ \cap U_{S_0^+}$ ($S^- \cap U_{S_0^+}$) is homeomorphic and tangent to the component $S_0^+$ ($S_0^-$) of the unstable (stable) cone of $V_0$, which is a cone over $Z_+$ ($Z_-$). More precisely, the solutions of the homogeneous system and the perturbed system starting from (ending at) the origin near $S_0^+$ ($S_0^-$) are naturally 1-1 correspondent and the corresponding pair of curves are tangent at the origin and naturally identified according to the distance to the origin.
\end{theorem}
{\bf Proof:}
For a solution ray of $V_0$ in $S_0^+$, choose the polar coordinate $(r,\theta)$ and $\theta = (\alpha,\beta)$ such that along the solution ray $\theta=0$ and locally $S_0^+ = \{\alpha=0\}$. Under such coordinate,

\[
V = r^{d-1}v(r,\theta)\frac{\partial}{\partial \theta} + r^dp(r,\theta)\frac{\partial}{\partial r},
\]

\[
V_0 = r^{d-1}v_0(\theta)\frac{\partial}{\partial \theta} + r^dp_0(\theta)\frac{\partial}{\partial r},
\]

\[
V_1 = r^dv_1(r,\theta)\frac{\partial}{\partial \theta} + r^{d+1}p_1(r,\theta)\frac{\partial}{\partial r},
\]

\[
v(r,\theta) = v_0(\theta) + rv_1(r,\theta),\ \ p(r,\theta) = p_0(\theta) + rp_1(r,\theta).
\]

The condition $|V_1| \leq Cr|V_0|$ implies that

\[
|v_1|^2 + p_1^2 \leq C(|v_0|^2 + p_0^2) = Cp_0^2(\left|\frac{v_0}{p_0}\right|^2 + 1) \leq Cp_0^2.
\]

\[
\frac{v}{p} = \frac{v_0}{p_0} + O(r).
\]

Locally we can write $\displaystyle\frac{v_0}{p_0} = -H(\theta)\alpha$. The hyperbolicity of $V_0$ implies that the real parts of all eigenvalues of $H(0,\beta)$ are strictly positive. Solution rays of the homogeneous system near the particular solution ray are the solutions of the following initial value problem.

\[
\left\{
\begin{array}{l}
\displaystyle\frac{d\alpha}{dr} = -\frac{1}{r}H(\theta)\alpha\\\\
\displaystyle\frac{d\beta}{dr} = 0\\\\
\theta(0)=(0,\beta_0)
\end{array}
\right.
\]

Solution rays of the inhomogeneous system as higher order perturbation of the homogeneous system are the solutions of the following perturbed initial value problem.

\[
\left\{
\begin{array}{l}
\displaystyle\frac{d\alpha}{dr} = -\frac{1}{r}H(\theta)\alpha + h_1(r,\theta)\\\\
\displaystyle\frac{d\beta}{dr} = h_2(r,\theta)\\\\
\theta(0)=(0,\beta_0)
\end{array}
\right.
\]

Here $h_1,h_2$ are bounded smooth functions on $(r,\theta)$. By proposition \ref{ce}, for each $\beta_0$ there exists a unique solution to this initial value problem that continuously depends on the initial value $\beta_0$. By proposition \ref{ch}, there are no other nearby solutions. The solutions of the two systems can be naturally identified, which gives us the homeomorphism we need.
\begin{flushright} $\Box$ \end{flushright}
For gradient vector field of a homogeneous function, we have the perturbation stability of stable and unstable manifolds under higher order perturbation of the standard flat metric.

\begin{co}
\label{cg}
Assume that $f$ is a homogneous function and $Z_+\subset S^{n-1}$ ($Z_-\subset S^{n-1}$) is a non-degenerate critical submanifold where $f|_{S^{n-1}}$ achieves maximum (minimum). For metric $g_{ij} = \delta_{ij} + O(r)$, let $S^+$ ($S^-$) denote the unstable (stable) variety of $V=\nabla f$. Then $S^+ \cap U_{S_0^+}$ ($S^- \cap U_{S_0^+}$) is homeomorphic and tangent to the cone $S_0^+$ ($S_0^-$) over $Z_+$ ($Z_-$) near the origin.
\end{co}

{\bf Proof:} Let $V_0=\nabla_0 f$ be the gradient vector field of $f$ with respect to the flat metric $\delta_{ij}$. Then $V_0$ is the leading homogeneous part of $V$. By applying theorem \ref{cf}, the result in this corollary is immediate.
\begin{flushright} $\Box$ \end{flushright}
{\bf Example:} On $\mathbb{C}^{{n+m}}$, consider $f(z) = {\rm Re}(s)$, where\\
\[
s= \prod_{i=1}^nz_i.
\]
Let $V_0=\nabla_0 f$, $V=\nabla f$ be the gradient vector field of $f$ with respect to the standard metric $\delta_{ij}$ and perturbed metric $g_{ij} = \delta_{ij} + a_{ij}$ respectively. Here $a_{ij} = O(|z|)$.\\

$f = {\rm Re}(s)$ when restricted to $S^{2n+2m-1}$ achieves its maximal on $Z_+$ and achieves its minimal on $Z_-$, where\\
\[
Z_\pm= \{z\in S^{2n-1}|s(z)=\pm 1/n^{n/2},\ |z_i|^2 = 1/n\ {\rm for\ } i\leq n,\ z_i=0\ {\rm for\ } i> n\}
\]\\
are two $(n-1)$-tori. Therefore, $V_0$ is hyperbolic along\\
\[
S_0^\pm= \{z|s(z)\in \mathbb{R}_\pm,\ |z_i|^2 = |z_j|^2\ {\rm for\ all\ } 1\leq i,j \leq n,\ z_i=0\ {\rm for\ } i> n\}.
\]\\
According to corollary \ref{cg}, (coupled with the results from the next section) the unstable (stable) variety $S^+$ ($S^-$) of $V=\nabla f$ are homeomorphic and tangent to the cone $S_0^+$ ($S_0^-$) over $Z_+$ ($Z_-$) near the origin. Topologically, it is a cone over a $(n-1)$-torus. We can move the origin within the subspace $\{z_i=0\ {\rm for\ } i> n\}$ to get the family of cones. This example will take care of the behavior of our gradient vector field near all stratas of $X_\infty \backslash X_\psi$.
\begin{flushright} $\Box$ \end{flushright}

\se{Singular gradient vector fields}
The above discussion on the behavior of the dynamical systems of higher order perturbation of a hyperbolic non-degenerate homogeneous vector field near the stable and unstable cones is quite complete. But the discussion of the behavior away from the stable and unstable cones is a bit short. For instance, it is not clear whether there are solution curves that spiral down and eventually approach the origin. In general this kind of behavior is hard to rule out. But for the kind of gradient homogeneous vector fields we are interested in, we can do better. The arguments are based on the detailed computation of the invariant functions of the homogeneous examples in section 3. There are essentially two cases and their parametrized versions corresponding to the two local models solved in section 3.\\

Recall that $f={\rm Re}(s)$ and $V=\frac{\nabla f}{|\nabla f|^2}$. Consider the flow\\
\[
\frac{dz}{dt} = V(z).
\]\\
This flow has two nice properties: $f(z(t))= f(z(0))+t$ and ${\rm Im}(s)$ is constant along the flow.\\

{\bf Case 1:} We will deal with the parametrized case directly. On $\mathbb{C}^{{m+n+l}}$, with coordinate $z\in {\mathbb{C}^{m+n}}\times \mathbb{C}^{l} \cong \mathbb{C}^{{m+n+l}}$, consider $f(z) = {\rm Re}(s)$, where

\[
s= \left.\left(\prod_{i=1}^mz_i\right)\right/\left(\prod_{i=m+1}^{m+n}z_i\right).
\]

Then with the standard flat metric,

\[
\nabla_0 f = {\rm Re}\left(\bar{s}\left(\sum_{i=1}^m \frac{z_i}{|z_i|^2} \frac{\partial}{\partial z_i} - \sum_{i=m+1}^{m+n} \frac{z_i}{|z_i|^2} \frac{\partial}{\partial z_i} \right) \right).
\]
\[
V_0= {\rm Re}\left(\left.\left( \sum_{i=1}^m \frac{z_i}{|z_i|^2} \frac{\partial}{\partial z_i} - \sum_{i=m+1}^{m+n} \frac{z_i}{|z_i|^2} \frac{\partial}{\partial z_i} \right) \right/s\sum_{i=1}^{m+n} \frac{1}{|z_i|^2} \right).
\]

Under the perturbed metric $g^{\bar{i}j} = \delta_{ij} + a_{ij}$, where $a_{ij}=O(|z|)$,

\[
V = {\rm Re}\left(\left.\left( \sum_{i=1}^m \frac{z_i}{|z_i|^2} g^{\bar{i}j}\frac{\partial}{\partial z_j} - \sum_{i=m+1}^{m+n} \frac{z_i}{|z_i|^2} g^{\bar{i}j}\frac{\partial}{\partial z_j} \right) \right/s\sum_{i=1}^{m+n} \frac{1}{|z_i|^2} \right)(1+O(|z|)).
\]\\

\begin{lm}
\label{dc}
When $m,n>0$, under the flow of V, a solution curve that starts away from the origin will always stay away from the origin.
\end{lm}

{\bf Proof:}
Let\\
\[
\rho = \sum_{i=1}^{m+n} \lambda_i|z_i|^2 + \sum_{i=m+n+1}^{m+n+l} |z_i|^2,
\]
with $\lambda_i>0$ and\\
\[
\sum_{i=1}^m \lambda_i =\sum_{i=m+1}^{m+n} \lambda_i.
\]

$\rho(z)$ is equivalent to $|z|^2$. We only need to show that along a solution curve that starts with $\rho\not= 0$, $\rho$ will not approach 0.

\begin{eqnarray*}
\frac{d\log \rho}{dt} &=& {\rm Re}\left(\left.\left( \sum_{i=1}^m \frac{z_i}{|z_i|^2} a_{ij}\frac{\partial \rho}{\partial z_j} - \sum_{i=m+1}^{m+n} \frac{z_i}{|z_i|^2} a_{ij}\frac{\partial \rho}{\partial z_j} \right) \right/s\rho\sum_{i=1}^{m+n} \frac{1}{|z_i|^2} \right)(1+O(|z|))\\
&\leq& C\left.\left(\sum_{i=1}^{m+n} \frac{1}{|z_i|} \right)|z|^2 \right/s\rho\sum_{i=1}^{m+n} \frac{1}{|z_i|^2} \leq C\left/s\left(\sum_{i=1}^{m+n} \frac{1}{|z_i|^2}\right)^{\frac{1}{2}}\right.\\
&\leq& C\left.\left(\prod_{i=1}^{m} z_i \right)^{\frac{1}{m}} \right/s \leq C\frac{1}{s^{1-\frac{1}{m}}} .
\end{eqnarray*}
Recall that we may take $t={\rm Re}(s)$. Then\\
\[
\lim_{t\rightarrow 0} \log \rho(t) \geq \log \rho(t_0) - C\int_0^{t_0}\frac{1}{t^{1-\frac{1}{m}}}dt \geq C.
\]\\
Namely\\
\[
\lim_{t\rightarrow 0} \rho(t) \not= 0.
\]\\
For the direction of $t\rightarrow +\infty$, we have\\
\begin{eqnarray*}
\frac{d\log \rho}{dt} &\leq& C\left/s\left(\sum_{i=1}^{m+n} \frac{1}{|z_i|^2}\right)^{\frac{1}{2}}\right.\\
&\leq& C\left.\left(\prod_{i=1}^{n} z_{m+i} \right)^{\frac{1}{n}} \right/s \leq C\frac{1}{s^{1+\frac{1}{n}}} \leq C\frac{1}{t^{1+\frac{1}{n}}}.
\end{eqnarray*}
Then\\
\[
\lim_{t\rightarrow +\infty} \log \rho(t) \geq \log \rho(t_0) - C\int_{t_0}^{+\infty}\frac{1}{t^{1+\frac{1}{n}}}dt \geq C.
\]\\
Namely\\
\[
\lim_{t\rightarrow +\infty} \rho(t) \not= 0.
\]
\begin{flushright} $\Box$ \end{flushright}

Let $X_c = \{s=c\}$ denote the level sets of $s$. In particular,\\
\[
X_0 = \left\{z\left|\prod_{i=1}^mz_i =0 \right.\right\},\ \ X_\infty = \left\{z\left| \prod_{i=m+1}^{m+n}z_i \right.\right\}.
\]

The flow of $V = \frac{\nabla f}{|\nabla f|^2}$ will move level sets to level sets. For details on toroidal \k metric mentioned in the following theorem, please refer to section 7.

\begin{theorem}
\label{da}
For $m,n>0$, on $\mathbb{C}^{{m+n+l}}$ with \k metric $g$ that is toroidal with respect to $X_0 \cup X_\infty$, where

\[
X_0 = \left\{z\left|\prod_{i=1}^mz_i =0 \right.\right\},\ \ X_\infty = \left\{z\left| \prod_{i=m+1}^{m+n}z_i \right.\right\},
\]

\[
X_{{\rm inv}} = X_0 \cap X_\infty = \left\{z\left|\prod_{i=1}^mz_i = \prod_{i=m+1}^{m+n}z_i =0 \right.\right\}
\]
is invariant under the flow of $V = \frac{\nabla f}{|\nabla f|^2}$, where $f(z) = {\rm Re}(s)$,

\[
s= \left.\left(\prod_{i=1}^mz_i\right)\right/\left(\prod_{i=m+1}^{m+n}z_i\right).
\]

A solution curve that starts from $X_0\backslash X_{{\rm inv}}$ will always stay away from $X_{{\rm inv}}$. More precisely, the flow of $V$ move $X_0\backslash X_{{\rm inv}}$ to $X_c\backslash X_{{\rm inv}}$.
\end{theorem}
{\bf Proof:}
The invariance of $X_{{\rm inv}}$ is quite obvious. Since the \k metric $g$ is toroidal with respect to $X_0 \cup X_\infty$, near any point $p\in X_{{\rm inv}}$, locally we get a situation as in lemma \ref{dc}. According to lemma \ref{dc} a solution curve starting from $X_0\backslash X_{{\rm inv}}$ will always stay away from $p$, therefore stay away from $X_{{\rm inv}}$.
\begin{flushright} $\Box$ \end{flushright}

{\bf Case 2:} To illustrate the idea more clearly, we will start with the non-parametrized version as an example.\\

{\bf Example:} On ${\mathbb{C}^n}$, consider $f(z) = {\rm Re}(s)$, where\\
\[
s= \prod_{i=1}^nz_i.
\]
Under the standard flat metric,\\
\[
\nabla_0 f = {\rm Re}\left(\bar{s}\sum_{i=1}^n \frac{z_i}{|z_i|^2} \frac{\partial}{\partial z_i} \right).
\]
\[
V_0 = {\rm Re}\left(\left.\left( \sum_{i=1}^n \frac{z_i}{|z_i|^2} \frac{\partial}{\partial z_i}\right) \right/s\sum_{i=1}^{n} \frac{1}{|z_i|^2} \right).
\]\\
Let\\
\[
\rho = \max_{i,j}\left||z_i|^2-|z_j|^2\right|.
\]
Then $\rho$ is a Lipschitz function. $U_\epsilon = \{\rho < \epsilon |z|^2\}$ is a cone neighborhood of the homogeneous unstable cone $S_0^+ = \{\rho=0\}$. Theorem \ref{cf} proved that a solution curve staying within $U_\epsilon$ when $r \rightarrow 0$ is in $S^+$, which is tangent and homeomorphic to $S_0^+$. Away from $U_\epsilon$,\\
\[
\epsilon |z|^2 \leq \rho \leq |z|^2.
\]
Due to the hyperbolicity of the vector field, $V$ is always pointing toward outside of $U_\epsilon$ along $\partial U_\epsilon$. Therefore a solution curve that starts from outside of $U_\epsilon$ will stay away from $U_\epsilon$.\\\\
Consider a solution curve away from $U_\epsilon$. Assume that $g^{i\bar{j}} = \delta_{ij} + a_{ij}$, where $a_{ij}=O(|z|)$. Then\\
\begin{eqnarray*}
\frac{d\log \rho}{dt} &=& {\rm Re}\left(\left.\left( \sum_{i=1}^n \frac{z_i}{|z_i|^2} a_{ij}\frac{\partial \rho}{\partial z_j} \right) \right/s\rho\sum_{i=1}^{n} \frac{1}{|z_i|^2} \right)(1 + O(|z|))\\
&\leq& C\left.\left(\sum_{i=1}^{n} \frac{1}{|z_i|} \right)|z|^2 \right/s\rho \sum_{i=1}^{n} \frac{1}{|z_i|^2} \leq C\left/s\left(\sum_{i=1}^{n} \frac{1}{|z_i|^2}\right)^{\frac{1}{2}}\right.\\
&\leq& C\left.\left(\prod_{i=1}^{n} z_i \right)^{\frac{1}{n}} \right/s \leq C\frac{1}{s^{1-\frac{1}{n}}}.
\end{eqnarray*}
Recall that we are considering solution curves along $\{{\rm Im}(s)=0\}$, where we may take $t=s$. Then\\
\[
\lim_{t\rightarrow 0} \log \rho(t) \geq \log \rho(t_0) - C\int_0^{t_0}\frac{1}{t^{1-\frac{1}{m}}}dt \geq C.
\]\\
Namely\\
\[
\lim_{t\rightarrow 0} \rho(t) \not= 0.
\]
\begin{flushright} $\Box$ \end{flushright}
Together with corollary \ref{cg} and the example afterward, we have actually proved the following

\begin{theorem}
With the metric $g_{i\bar{j}} = \delta_{ij} + O(|z|)$, the unstable variety $S^+$ of $V = \frac{\nabla f}{|\nabla f|^2}$ at the origin is homeomorphic and tangent to the cone\\
\[
S_0^+ = \{z| s= \prod_{i=1}^nz_i \in \mathbb{R}_{\geq 0},\ |z_i|=|z_j|=0\},
\]
which is a cone over an $(n-1)$-torus. Any solution curve of $V$ with the origin as limit point is inside $S^+$.
\end{theorem}
\begin{flushright} $\Box$ \end{flushright}
This discussion is not hard to be generalized to the parametrized case. On ${\mathbb{C}^{m+n}}$ with coordinate $(z,w)\in {\mathbb{C}^n}\times \mathbb{C}^{m} \cong {\mathbb{C}^{m+n}}$, consider $f(z) = {\rm Re}(s)$, where\\
\[
s= \prod_{i=1}^nz_i.
\]
Assume that $g^{\bar{i}j} = \delta_{ij} + a_{ij}$, where $a_{ij}=O(|(z,w)|)$. Then

\[
V = \frac{\nabla f}{|\nabla f|^2} = {\rm Re}\left(\left.\left( \sum_{i=1}^n \frac{z_i}{|z_i|^2} g^{\bar{i}j}\frac{\partial \rho}{\partial z_j} \right) \right/s\sum_{i=1}^{n} \frac{1}{|z_i|^2} \right)(1 + O(|z|)).
\]

\begin{theorem}
\label{dd}
With the metric $g_{i\bar{j}} = \delta_{ij} + O(|(z,w)|)$, the unstable variety $S^+$ of $V = \frac{\nabla f}{|\nabla f|^2}$ at the origin is homeomorphic and tangent to the cone\\
\[
S_0^+ = \{(z,0)| s= \prod_{i=1}^nz_i \in \mathbb{R}_{\geq 0},\ |z_i|=|z_j|=0\},
\]
which is a cone over an $(n-1)$-torus. Any solution curve of $V$ with the origin as limit point is inside $S^+$.
\end{theorem}

{\bf Proof:}
Let\\
\[
\rho = \max_{i,j}\left||z_i|^2-|z_j|^2\right| + |w|^2.
\]
Then $\rho$ is a Lipschitz function. $U_\epsilon = \{\rho < \epsilon |(z,w)|^2\}$ is a cone neighborhood of the homogeneous unstable cone $S_0^+ = \{\rho=0\}$. Theorem \ref{cf} proved that $S^+ \cap U_\epsilon$ is tangent and homeomorphic to $S_0^+$. Away from $U_\epsilon$,\\
\[
\epsilon |(z,w)|^2 \leq \rho \leq |(z,w)|^2.
\]
Due to the hyperbolicity of the vector field, $V$ is always pointing toward outside of $U_\epsilon$ along $\partial U_\epsilon$. Therefore a solution curve that starts from outside of $U_\epsilon$ will stay away from $U_\epsilon$.\\\\
Consider a solution curve away from $U_\epsilon$. Assume that $g^{i\bar{j}} = \delta_{ij} + a_{ij}$, where $a_{ij}=O(|(z,w)|)$. Then similarly as in the example we can show\\
\[
\lim_{t\rightarrow 0} \log \rho(t) \geq \log \rho(t_0) - C\int_0^{t_0}\frac{1}{t^{1-\frac{1}{m}}}dt \geq C.
\]\\
Namely\\
\[
\lim_{t\rightarrow 0} \rho(t) \not= 0.
\]

Therefore, no other components of unstable variety $S^+$ exist and any solution curve of $V$ with the origin as limit point is inside $S^+$.
\begin{flushright} $\Box$ \end{flushright}
\begin{theorem}
\label{db}
On $\mathbb{C}^{{n+m}}$ with \k metric $g$ that is toroidal with respect to $X_0$, where

\[
X_0 = \left\{z\left|\prod_{i=1}^nz_i =0 \right.\right\},
\]

for every point in $\{z_i=0,\ {\rm for}\ i\leq n\}$, the unstable variety $S^+$ of $V = \frac{\nabla f}{|\nabla f|^2}$ at this point is topologically a cone over an $(n-1)$-torus, where $f(z) = {\rm Re}(s)$,

\[
s= \prod_{i=1}^nz_i.
\]
\end{theorem}
{\bf Proof:}
Since the \k metric $g$ is toroidal with respect to $X_0$. Near any point $p\in \{z_i=0,\ {\rm for}\ i\leq n\}$, locally we get a situation as in theorem \ref{dd}. By theorem \ref{dd}, the result in this theorem is immediate.
\begin{flushright} $\Box$ \end{flushright}
{\bf Remark:} In this section, for simplicity of notation, we carried out all our arguments only for the unstable variety $S^+$. The corresponding arguments and results for the stable variety $S^-$ are identical up to a sign. We also used the concept of toroidal \k metric that will be discussed in section 7.\\

Summarize our local discussions, we have the following global result. Assume that we have a family of hypersurfaces $\{X_s\}$ in an ambient compact \k manifold $(M, \omega_g)$ of dimension $n+1$. Assume that $X_s$ is smooth for $s\not=0$ and $X_0$ is a divisor in $M$ with only normal crossing singularities. Let the disjoint union

\[
X_0 = \bigcup_{i=0}^n X_0^{(i)}
\]

be the stratification of $X_0$, where $X_0^{(i)}$ denotes the $i$-dimensional strata of $X_0$. Then

\[
{\rm Sing}(X_0) = \bigcup_{i=0}^{n-1} X_0^{(i)}.
\]

\begin{de}
\label{df}
Assume that $(M_1,\omega_1)$ is a smooth symplectic manifold and $(M_2,\omega_2)$ is a symplectic variety. Then a piecewise smooth map $H: M_1 \rightarrow M_2$ is called a symplectic morphism if $H^*\omega_2 = \omega_1$. If $(M_2,\omega_2)$ is also a smooth symplectic manifold and $H$ is a diffeomorphism, then the symplectic morphism $H$ is also called a symplectomorphism.
\end{de}

{\bf Remark:} In this paper, we will only deal with the case of normal crossing symplectic varieties. Therefore, we will not venture into the concept of general symplectic varieties and symplectic forms on them.\\

$s$ defines a meromorphic function on $M$. Consider the flow of $V = \frac{\nabla f}{|\nabla f|^2}$, where $f(z) = {\rm Re}(s)$.

\begin{theorem}
\label{de}
Assume that $X_{\rm inv} = X_s \cap X_0$ is independent of $s$. And assume the \k metric $g$ on $M$ is toroidal along ${\rm Sing}(X_0)$ with respect to $X_0\cup X_s$. Then the flow of $V$ will fix $X_{\rm inv}$ and flow each point in $X_0^{(i)}\backslash X_{\rm inv}$ to a real $(n-i)$-torus in $X_s$ for $s$ real. In another word, the inverse flow of $V$ will induce a symplectic morphism $F_s: X_s \rightarrow X_0$ with respect to the toroidal \k form $\omega_g$ for $s$ real. For $x\in X_0^{(i)}\backslash X_{\rm inv}$, $F_s^{-1}(x)$ is a real $(n-i)$-torus in $X_s$. For $x\in X_{\rm inv}$, $F_s^{-1}(x) = x$.
\end{theorem}

{\bf Proof:} For $x\in X_0\backslash {\rm Sing}(X_0)$, $X_s$ are all smooth near $x$. Therefore, the flow induces diffeomorphisms near $x$.\\

Since $g$ on $M$ is toroidal along ${\rm Sing}(X_0)$ with respect to $X_0\cup X_s$. For $i<n$ and $x\in X_0^{(i)}\backslash X_{\rm inv}$, all smooth components of the normal crossing divisor $X_0\cup X_s \subset M$ are orthogonal to each other at $x$ under the \k metric $g$. It is easy to find a local holomorphic coordinate $z$ such that $z(x)=0$, $\displaystyle s = c\prod_{k=1}^{n+1-i}z_k$, $g_{i\bar{j}} = \delta_{ij} + O(|z|)$. The theorem \ref{db} implies that the flow of $V$ will flow $x\in X_0^{(i)}\backslash X_{\rm inv}$ to a real $(n-i)$-torus in $X_s$.\\

Similarly, for $i<n$ and $x\in X_0^{(i)}\cap X_{\rm inv}$, all smooth components of the normal crossing divisor $X_0\cup X_s \subset M$ are orthogonal to each other at $x$ under the \k metric $g$. It is easy to find a local holomorphic coordinate $z$ such that $z(x)=0$, $\displaystyle s = \frac{c}{z_1}\prod_{k=2}^{n+2-i}z_k$, $g_{i\bar{j}} = \delta_{ij} + O(|z|)$. The theorem \ref{da} implies that the flow of $V$ will fix $x\in X_{\rm inv}$, and under the inverse flow of $V$, no point in $X_s$ other than $x$ will flow to $x$.
\begin{flushright} $\Box$ \end{flushright}

\se{Hamiltonian deformation of submanifolds of symplectic manifold}
For a manifold $X$, consider the submanifolds $S_0,S_1\subset X$ that are isotopic, by which we mean that there exists a map\\
\[
{\cal F}: S\times [0,1] \rightarrow X
\]
such that ${\cal F}_t$ are embeddings, where ${\cal F}_t(x) = {\cal F}(x,t)$, with ${\cal F}_0(S)=S_0$ and ${\cal F}_1(S)=S_1$. Assume that there are two deformation equivalent cohomologous symplectic forms $\omega_0, \omega_1$ on $X$. Namely there is a smooth family $\{\omega_t\}_{t\in [0,1]}$ of symplectic forms of the same cohomology class connecting $\omega_0, \omega_1$ on $X$. We will first define the concept of isotopy of submanifolds in the symplectic context. There are two classes of submanifolds that are particularly interesting for our discussion --- Lagrangian and symplectic submanifolds. Assume further that $S_t={\cal F}_t(S)\subset X$ are Lagrangian (symplectic) submanifolds of symplectic manifold $(X,\omega_t)$ respectively in the two situations for any $t$. In another word, ${\cal F}_t$ is a Lagrangian (symplectic) isotopy between $(S_0\subset X, \omega_0)$ and $(S_1\subset X, \omega_1)$. A natural question is when there will exist a Hamiltonian (symplectic) morphism $h: (X,\omega_0) \rightarrow (X,\omega_1)$ such that $h(S_0)=S_1$. We will mainly concentrate on the symplectic submanifold situation.\\

$\{{\cal F}_t\}_{t\in [0,1]}$ can be viewed as a flow on the family $\{S_t\}_{t\in [0,1]}$. Let $w_t$ be the vector fields generating the flow $\{{\cal F}_t\}_{t\in [0,1]}$. For convenience of arguments, we usually extend $w_t$ to global vector fields on $X$.\\

\begin{de}
Two symplectic submanifold structures $(S_0\subset X, \omega_0)$ and $(S_1\subset X, \omega_1)$ are {\bf symplectic isotopic} if there exists  a smooth hamiltonian equivalent family $\{\omega_t\}_{t\in [0,1]}$ of symplectic forms connecting $\omega_0, \omega_1$ on $X$ and a smooth family of embeddings ${\cal F}_t: S \rightarrow X$ such that ${\cal F}_0(S)=S_0$, ${\cal F}_1(S)=S_1$ and $S_t={\cal F}_t(S)$ are symplectic with respect to $\omega_t$ for any $t$.

If in addition, ${\cal F}_t\circ {\cal F}_0^{-1}: S_0\rightarrow S_t$ are symplectomorphisms for all $t$, then $\{{\cal F}_t\}_{t\in [0,1]}$ is called a {\bf symplectic flow}.

If $\{\omega_t\}_{t\in [0,1]}$ is a constant family, $i(v_t)\omega_t|_{S_t}$ are exact 1-forms for all $t$, then $\{{\cal F}_t\}_{t\in [0,1]}$ is called a {\bf Hamiltonian flow}.

A symplectic (Hamiltonian) flow $\{{\cal F}_t\}_{t\in [0,1]}$ is called $C^{0,1}$ (Lipschitz) if the corresponding vector fields $v_t$ are uniformly $C^{0,1}$ (Lipschitz) for all $t$.
\end{de}

In the proofs, we will need the cutoff function $c_1(y)$ that satisfies\\
\[
c_1(y) = \left\{
\begin{array}{ll}
1,\ &|y|\leq 1\\
0,\ &|y|\geq 2
\end{array}
\right..
\]\\
Naively, the problem of modifying a symplectic isotopy to a symplectic (Hamiltonian) flow can be solved in two steps. The first step is to modify a symplectic isotopy $\{{\cal F}_t\}_{t\in [0,1]}$ to a symplectic (Hamiltonian) flow on the family $\{S_t\}_{t\in [0,1]}$. The second step is to extend $\{{\cal F}_t\}_{t\in [0,1]}$ to a symplectic (Hamiltonian) flow on $X$. In practice, the extension in the second step is not always possible if the construction of symplectic flow in the first step is not done with care. We will first deal with the case when $S$ is a manifold, where the proof is straightforward and clear. Later, we will handle the more general case when $S$ is a union of smooth submanifolds with normal crossing, which is a bit more delicate.\\

Historically, the investigation of Hamiltonian deformation of symplectic manifold started with the famous theorem of Moser. Results in ``the basic case" subsection is a natural generalization to the symplectic submanifold case, which should be well known. We discuss this basic case first to illustrate the basic idea of our method, which is of constructive nature. Later we extend our method and results to several different cases necessary for our applications. One notable difference between our results and Moser type results is that the symplectomorphisms we construct are mostly piecewise smooth (Lipschitz) instead of smooth. To get smooth symplectomorphisms, additional conditions are needed as discussed in ``the smooth case" subsection.\\

{\bf Remark:} There are two places where we will apply results in this section. The first application is to construct a $C^{0,1}$ symplectomorphism from $\mathbb{CP}^4$ with the Fubini-Study metric to $\mathbb{CP}^4$ with a toroidal metric with respect to $X_\infty\cup X_0$ that maps $X_\infty\cup X_0$ to itself (theorem \ref{eg}). The second application is to deform the symplectic curves in $\mathbb{CP}^2$'s to achieve graph image (theorem \ref{ee} and corollary \ref{er}), and extend to $X_\infty \subset \mathbb{CP}^4$ (theorem \ref{ef}).\\

{\bf Remark on notation:} In this section, we will quite often use $v$ and $\alpha$ to denote a vector field and a 1-form that would be modified together and are mutually determined by the equation $i(v)\omega = \alpha$. For simplicity of notation, we will quite oftenly use $b(x)dy$, $b(x)y$ to denote $\displaystyle\sum_{i}b_i(x)dy_i$, $\displaystyle\sum_{i}b_i(x)y_i$ when $y$ is a multivariable coordinate. When $S$ is a normal crossing union of manifolds, $(S_t\subset X, \omega_t)$ is usually refered to as {\bf symplectic submanifold structure} instead of symplectic submanifold, which we reserve to the case that $S$ is a manifold. Unless specified otherwise or obvious from the context, in this section, when we mention manifolds or their submanifolds (with or without boundary), we always assume that they are {\bf compact}.\\

\subsection{The basic case}
In this subsection, we will consider the case when $S_t\subset X$ are smooth submanifolds to illustrate the basic idea. It also serves as the first step of the induction argument of general cases.

\begin{theorem}
\label{ea}
Assume that two symplectic submanifolds $(S_0\subset X, \omega_0)$ and $(S_1\subset X, \omega_1)$ are symplectic isotopic, then there exists a $C^\infty$ symplectomorphism $h: (X,\omega_0) \rightarrow (X,\omega_1)$ such that $h(S_0)=S_1$.
\end{theorem}
{\bf Proof:}
Since $\{\omega_t\}_{t\in [0,1]}$ is a cohomologous family, we may choose a suitable family of 1-forms $\{\alpha_t\}_{t\in [0,1]}$ such that\\
\[
\frac{d\omega_t}{dt} = -d\alpha_t.
\]\\
Let $w_t={\cal F}_*(\frac{d}{dt})$. Since $S_t$ is symplectic, we can find $u_t\in TS_t$ such that\\
\[
i(u_t)\omega_t|_{S_t} =i(w_t)\omega_t|_{S_t} - \alpha_t|_{S_t}.
\]
Let $v_t= w_t - u_t$, then\\
\[
i(v_t)\omega_t|_{S_t} = \alpha_t|_{S_t}.
\]
The flow of $v_t$ will map $(S_0,\omega_0)$ to $(S_t,\omega_t)$ (in particular $(S_1,\omega_1)$) symplectically. To get a globally defined flow, we need to extend $v_t$ to $X$ from $S_t$. For this purpose we need to find a function $f_t$ that satisfies\\
\[
i(v_t)\omega_t|_{T^*_{S_t}X} = (\alpha_t + df_t)|_{T^*_{S_t}X}
\]\\
where $T^*_{S_t}X = T^*X|_{S_t}$. Then $v_t$ defined by\\
\[
i(v_t)\omega_t = \alpha_t + df_t
\]\\
will be our desired extension. Since\\
\[
i(v_t)\omega_t|_{S_t} = \alpha_t|_{S_t}.
\]\\
Locally on $S_t$, assume $x,y$ are tangent and normal coordinates near $S_t$. We may express\\
\[
(i(v_t)\omega_t-\alpha_t)= b(x)dy + O(|y|).
\]\\
We may define\\
\[
f_t = b(x)y c_1(y/\epsilon),
\]\\
where $\epsilon$ is a small positive constant. Globally we may use partition of unity to piece such local $f_t$ together to form a function $f_t$ on $X$ that is supported near $S_t$. Here we will perform this rather standard operation once and for all in this section for readers' convenience.\\

Let $\{U_\gamma\}$ be an open covering of a small tubular neighborhood $U_{S_t}$ of $S_t$, $\{\rho_\gamma\}$ is a partition of unity with respect to $\{U_\gamma\}$ that satisfies $\sum_\gamma \rho_\gamma|_{S_t} = 1$, and $x_\gamma,y_\gamma$ are tangent and normal coordinates near $S_t$ in $U_\gamma$. We may define $f_{\gamma,t} = b_\gamma(x_\gamma)y_\gamma c_1(y_\gamma/\epsilon)$ as above. Then

\[
f_t = \sum_\gamma \rho_\gamma f_{\gamma,t}
\]

is the desired global function supported near $S_t$.\\

The flow of corresponding $v_t$ (that satisfies $i(v_t)\omega_t = \alpha_t + df_t$) will produce the desired symplectomorphism $h$.
\begin{flushright} $\Box$ \end{flushright}
\begin{theorem}
\label{eb}
Assume that $\{{\cal F}_t\}_{t\in [0,1]}$ is a symplectic flow on the family $\{S_t\}_{t\in [0,1]}$ and $H_1(S_t)\rightarrow H_1(X)$ is injective (for example, when $H_1(S)=0$), then there exists a $C^\infty$ symplectomorphism $h: (X,\omega_0) \rightarrow (X,\omega_1)$ such that $h(S_0)=S_1$ and $h|_{S_0} = {\cal F}_1\circ {\cal F}_0^{-1}$.
\end{theorem}
{\bf Proof:}
Let $v_t={\cal F}_*(\frac{d}{dt})$. Since $\{\omega_t\}_{t\in [0,1]}$ is a cohomologous family, under the assumption of the theorem, we may choose a suitable family of 1-forms $\{\alpha_t\}_{t\in [0,1]}$ such that\\
\[
\frac{d\omega_t}{dt} = -d\alpha_t
\]\\
and $\alpha_t|_{S_t}$ is cohomologous to $i(v_t)\omega_t|_{S_t}$. Namely\\
\[
i(v_t)\omega_t|_{S_t} = (\alpha_t + da_t)|_{S_t}.
\]\\
$a_t$ can be easily extended to functions on $X$. Replace $\alpha_t$ by $\alpha_t + da_t$, we may assume\\
\[
\frac{d\omega_t}{dt} = -d\alpha_t
\]\\
and $\alpha_t|_{S_t}=i(v_t)\omega_t|_{S_t}$. To extend our flow, we need to extend $v_t$ to $X$ from $S_t$. For this purpose we need to find a function $f_t$ that satisfies\\
\[
i(v_t)\omega_t|_{T^*_{S_t}X} = (\alpha_t + df_t)|_{T^*_{S_t}X}.
\]\\
Then $v_t$ defined by\\
\[
i(v_t)\omega_t = \alpha_t + df_t
\]\\
will be our desired extension. This can be done similarly as in the proof of theorem \ref{ea}.
\begin{flushright} $\Box$ \end{flushright}
{\bf Remark:} Theorem \ref{eb} is about extending symplectic flow from symplectic submanifold. Naively theorem \ref{ea} is a corollary of theorem \ref{eb} if we can modify ${\cal F}_t$ into a symplectic flow. In practice, the choice of symplectic flow in theorem \ref{ea} is very special and theorem \ref{ea} applies to more general cases than theorem \ref{eb}.\\

\subsection{The general case}
For our application, it is necessary to consider the case when $S_t\subset X$ is a union of several smooth normal crossing submanifolds that are symplectic. We will use induction to construct the symplectic flow in such situation. We will modify the symplectic flow generated by $v_t$ and the symplectic isotopy for $\{S_t\}$ generated by $w_t$ strata by strata. In the end, $v_t$ will coincide with $w_t$ when restricted to $S_t$, therefore finishing our construction. For this purpose, let us consider a structure $(S^1\subset S^2 \subset U_{S^2} \subset X, \omega)$, where $S^1$ is closed in $S^2$, $S^2\subset (X, \omega)$ is a symplectic submanifold, $S^1\subset (X, \omega)$ is a finite union of symplectic submanifolds that are normal crossing to each other and $U_{S^2}$ is a neighborhood of $S^2$ with a natural identification to a neighborhood of zero section of $N_{S^2}$. This identification induces a natural projection $\pi: U_{S^2}\rightarrow S^2$. (For our application, $S_1,S_2$ will be taken as parts of $S_t$.)\\
\begin{lm}
\label{ed}
Assume that we have $(S^1\subset S^2 \subset U_{S^2} \subset X, \omega)$ and vector fields $v,w$ satisfying $v|_{S^1} = w|_{S^1}$. Then $\alpha = i(v)\omega$ can be modified by adding a piecewise smooth $C^{0,1}$ exact 1-form such that for the modified vector field $v$ satisfying $i(v)\omega = \alpha$ (modified $\alpha$), we have $v|_{\pi^{-1}(S^1)}$ is unchanged and $(v-w)|_{S^2}$ is along $S^2$.\\
\end{lm}
{\bf Proof:} Since $S^2$ is a symplectic submanifold, there is a unique vector field $u$ on $S^2$ such that\\
\begin{equation}
\label{ec}
i(w-u)\omega|_{S^2} = \alpha|_{S^2}.
\end{equation}\\
We need to find a function $f$ that satisfies\\
\[
i(w-u)\omega|_{T^*_{S^2}X} = (\alpha + df)|_{T^*_{S^2}X}.
\]\\
Let us first look at the problem locally. Assume that we have local coordinate $(x,y)$ of $U_{S^2}$ such that $\pi(x,y) = x$ and $S^2=\{y=0\}$. Equality (\ref{ec}) implies that\\
\[
i(w-u)\omega - \alpha = b(x)dy + O(|y|).
\]

$v|_{S^1} = w|_{S^1}$ implies that $b|_{S^1}=0$. Let\\
\[
f = b(x)y c_1(y/\rho(x)),
\]\\
where $\epsilon$ is a small positive constant and $\rho(x)\leq \epsilon$ is the distance function to $S^1$ on $S^2$ when $\rho(x)\leq \epsilon/2$. It is easy to see that $f$ is $C^{1,1}$, which is equivalent to the second derivatives being bounded. We will verify the boundedness for $f_{yy}$. The verification for other second derivatives is similar.\\
\[
f_{yy} = 2(b(x)/\rho(x)) c'_1(y/\rho(x)) + (b(x)/\rho(x))(y/\rho(x)) c''_1(y/\rho(x)).
\]

Notice that $b|_{S^1}=0$ and $c'_1(y/\rho(x))$ is nonvanishing only when $1\leq y/\rho(x) \leq 2$. Therefore $f_{yy}$ is bounded.\\

Globally, one can define $f$ by using partition of unity to piece together $f$ defined in each local chart (as in the proof of theorem \ref{ea}). The desired modification is achieved by replacing $\alpha$ by $\alpha + df$. The corresponding modified $v$ will satisfy\\
\[
i(w-v-u)\omega = O(|y|).
\]

Therefore $(v-w)|_{S^2}$ is along $S^2$.\\

Since $c_1(y/\rho(x))=0$ when $y/\rho(x) \geq 2$, we have $f|_{\{y/\rho(x) \geq 2\}} =0$. This together with $\pi^{-1}(S^1) \subset \{y/\rho(x) \geq 2\}$ imply $v|_{\pi^{-1}(S^1)}$ is unchanged.
\begin{flushright} $\Box$ \end{flushright}
{\bf Remark:} For the sake of induction process, it is also necessary to modify $w$ to ensure $w|_{S^2} = v|_{S^2}$. This can be easily done by replacing $w$ with $w - \sum_\gamma \rho_\gamma(x_\gamma,y_\gamma) u(x_\gamma)$, where $\{U_\gamma\}$ is a covering of $U_{S^2}$, $\{\rho_\gamma\}$ is a partition of unity that satisfies $\sum_\gamma \rho_\gamma|_{S^2} = 1$ and $(x_\gamma,y_\gamma)$ is the local coordinate on $U_\gamma$ that satisfies $\pi(x_\gamma,y_\gamma) = x_\gamma$, $S^2 \cap U_\gamma = \{y_\gamma =0\}$. (To ensure $u(x_\gamma)$ to be along $S_t$, it is necessary to choose $y_\gamma$ suitably so that each smooth component of $S_t \cap U_\gamma$ is defined by some components of $y_\gamma$ being zero.)\\

What we are really interested in is a family of such structure $(S^1_t\subset S^2_t \subset U_{S_t^2} \subset X, \omega_t)$. Here $S^1_t\subset S^2_t \subset S_t$ should be understood as the images of a family of embedding ${\cal F}_t: (S^1\subset S^2 \subset S) \rightarrow X$, where ${\cal F}_t(x) = {\cal F}(t,x)$. Then we have\\
\begin{lm}
\label{ej}
Assume that we have a family of structure $(S^1_t\subset S^2_t \subset U_{S_t^2} \subset X, \omega_t)$ defined as above, and there is a family of 1-forms $\{\alpha_t\}_{t\in [0,1]}$ such that\\
\[
\frac{d\omega_t}{dt} = -d\alpha_t,\ \ i(v_t)\omega_t = \alpha_t
\]\\
and the flow generated by $v_t$ restricts to symplectic isotopy on the family $\{S^1_t\}_{t\in [0,1]}$. Then we can modify $\{\alpha_t\}_{t\in [0,1]}$ by adding piecewise smooth $C^{0,1}$ exact forms such that the modified $v_t$ is unchanged on $\pi_t^{-1}(S_t^1)$ and the flow generated by $v_t$ restricts to symplectic isotopy on the family $\{S^2_t\}_{t\in [0,1]}$. In particular, there exists a piecewise smooth $C^{0,1}$ symplectomorphism $h: (X,\omega_0) \rightarrow (X,\omega_1)$ such that $h(S^2_0)=S^2_1$.
\end{lm}
{\bf Proof:} The condition ``the flow generated by $v_t$ restricts to symplectic isotopy on the family $\{S^1_t\}_{t\in [0,1]}$'' implies that\\
\[
v_t|_{S^1_t} = {\cal F}_*(\frac{d}{dt})|_{S^1_t}.
\]\\
Choose a vector field $w_t$ such that\\
\[
w_t|_{S_t} = {\cal F}_*(\frac{d}{dt})|_{S_t}.
\]\\
Then according to lemma \ref{ed}, $\alpha_t$ can be modified by exact 1-form such that $v_t|_{\pi^{-1}(S_t^1)}$ is unchanged and $(v_t-w_t)|_{S_t^2}$ is along $S_t^2$. Namely, the flow generated by $v_t$ is unchanged on $\pi_t^{-1}(S_t^1)$ and restricts to symplectic isotopy on the family $\{S^2_t\}_{t\in [0,1]}$ as we claimed.
\begin{flushright} $\Box$ \end{flushright}
Consider two symplectic isotopic symplectic submanifold structures $(S_0\subset X, \omega_0)$ and $(S_1\subset X, \omega_1)$. Let\\
\[
{\cal F}: S\times [0,1] \rightarrow X
\]\\
be the corresponding isotopic map such that ${\cal F}_t(x) = {\cal F}(x,t)$. We will generalize the concept of submanifold structure and the corresponding symplectic isotopy by allowing $S$ to be a union of smooth manifolds with normal crossings. Then we have\\
\begin{theorem}
\label{eh}
Assume that two symplectic submanifold structures $(S_0\subset X, \omega_0)$ and $(S_1\subset X, \omega_1)$ are symplectic isotopic, where $S$ is a union of smooth manifolds with normal crossings. Then there exists a piecewise smooth $C^{0,1}$ symplectomorphism $h: (X,\omega_0) \rightarrow (X,\omega_1)$ such that $h(S_0)=S_1$.
\end{theorem}
{\bf Proof:} Since $\{\omega_t\}_{t\in [0,1]}$ is a hamiltonian equivalent family, we may choose a suitable family of 1-forms $\{\alpha_t\}_{t\in [0,1]}$ such that\\
\[
\frac{d\omega_t}{dt} = -d\alpha_t,\ \ i(v_t)\omega_t = \alpha_t.
\]\\
The flow of $v_t$ is symplectic, but usually does not restrict to symplectic isotopy to the family $\{S_t\}_{t\in [0,1]}$. On the other hand, ${\cal F}_*(\frac{d}{dt})$ on $S_t$ can be naturally extended to a smooth vector field $w_t$ on $X$ supported in a neighborhood of $S_t$. $S$ induces a natural filtration $\emptyset = S^{(0)}\subset S^{(1)} \subset \cdots \subset S^{(l)}= S$, such that $S^{(k)}\backslash S^{(k-1)}$ are open manifolds. We will modify $\alpha_t$ inductively, strata by strata according to the filtration to ensure that the flow of $v_t$ will restrict to symplectic isotopy to the family $\{S_t\}_{t\in [0,1]}$.\\\\
Assume that the symplectic flow of $v_t$ restricts to symplectic isotopy for the family $\{S^{(k-1)}_t\}_{t\in [0,1]}$. More precisely $w_t|_{S^{(k-1)}_t} = v_t|_{S^{(k-1)}_t}$. $S_t^{(k)}$ is a union of several smooth manifolds with normal crossing. Let $S^2_t$ be one of these component, and let $S^1_t=S^2_t\cap S^{(k-1)}_t$. Choose $U_{S^2_t}$ and corresponding $\pi_t: U_{S^2_t} \rightarrow S^2_t$ suitably so that $(S_t^{(k)}\backslash S^2_t)\cap U_{S^2_t} \subset \pi_t^{-1}(S^1_t)$. Apply lemma \ref{ej} to $(S^1_t\subset S^2_t \subset U_{S_t^2} \subset X, \omega_t)$, we can modify $\alpha_t$ to ensure that the symplectic flow of $v_t$ will restrict to symplectic isotopy for the family $\{S_t^2\}_{t\in [0,1]}$ while the flow on other components of $S_t^{(k)}$ is unaffected. As in the remark after lemma \ref{ej}, we can also correspondingly modify $w_t$ so that $w_t|_{S^2_t} = v_t|_{S_t^2}$ and the flow of $w_t$ will still restrict to symplectic isotopy for the family $\{S_t\}_{t\in [0,1]}$ while the flow of $w_t$ on other components of $S_t^{(k)}$ is unaffected. Repeat this process to other components of $S_t^{(k)}$, we can ensure that the flow of $v_t$ will restrict to symplectic isotopy to the family $\{S_t^{(k)}\}_{t\in [0,1]}$. More precisely $w_t|_{S^{(k)}_t} = v_t|_{S^{(k)}_t}$. By induction on $k$, the theorem can be proved.
\begin{flushright} $\Box$ \end{flushright}
{\bf Remark:} The proof actually proved the result for more general $S$ that possesses a natural filtration $\emptyset = S^{(0)}\subset S^{(1)} \subset \cdots \subset S^{(l)}= S$, such that $S^{(k)}\backslash S^{(k-1)}$ is a disjoint union of open manifolds and the closure of $S^{(k)}\backslash S^{(k-1)}$ is a normal crossing union of closed manifolds for each $k$. Of course the corresponding concepts of symplectic isotopy, symplectic (Hamiltonian) flow, etc., should be adjusted to incorporate compatibility with the filtration.\\\\
We also have extension theorem in the general case as theorem \ref{eb}. Let us start with the infinitesimal version.\\
\begin{lm}
\label{ei}
Assume that we have $(S^1\subset S^2 \subset U_{S^2} \subset X, \omega)$, a vector field $w$ and a 1-form $\alpha$ such that\\
\[
i(v)\omega = \alpha,\ \ \ \ v|_{S^1} = w|_{S^1}
\]\\
and $(\alpha-i(w)\omega)|_{S^2} = 0$ on $S^2$. Then $\alpha$ can be modified by adding a piecewise smooth $C^{0,1}$ exact 1-form such that the modified $v$ is unchanged on $\pi^{-1}(S^1)$ and $(v-w)|_{S^2}=0$.\\
\end{lm}
{\bf Proof:} We need to find a function $f$ that satisfies\\
\[
i(w)\omega|_{T^*_{S^2}X} = (\alpha + df)|_{T^*_{S^2}X}.
\]\\
Let us first look at the problem locally. Assume that we have local coordinate $(x,y)$ of $U_{S^2}$ such that $\pi(x,y) = x$. Since $(\alpha-i(w)\omega)|_{S^2} = 0$ on $S^2$, we may assume that\\
\[
i(w)\omega - \alpha = b(x)dy + O(|y|).
\]\\
Define\\
\[
f = b(x)y c_1(y/\rho(x)),
\]\\
where $\epsilon$ is a small positive constant and $\rho(x)\leq \epsilon$ is the distance function to $S^1$ on $S^2$ when $\rho(x)\leq \epsilon/2$. It is easy to see that $f$ is $C^{1,1}$.\\\\
Globally, one can define $f$ by using partition of unity to piece together $f$ defined in each local chart (as in the proof of theorem \ref{ea}). The desired modification is done by replacing $\alpha$ by $\alpha + d f$.
\begin{flushright} $\Box$ \end{flushright}
Using this lemma combined with the method in theorem \ref{eb} and similar induction argument as in the proof of theorem \ref{eh}, we can show the following extension theorem.\\
\begin{theorem}
Assume that $\{{\cal F}_t\}_{t\in [0,1]}$ is a piecewise smooth $C^{0,1}$ symplectic flow on the family $\{S_t\}_{t\in [0,1]}$, $S$ is a union of smooth manifolds with normal crossing and $H_1(S_t)\rightarrow H_1(X)$ is injective (for example, when $H_1(S)=0$), then there exists a piecewise smooth $C^{0,1}$ symplectomorphism $h: (X,\omega_0) \rightarrow (X,\omega_1)$ such that $h(S_0)=S_1$ and $h|_{S_0} = {\cal F}_1\circ {\cal F}_0^{-1}$.
\end{theorem}
{\bf Proof:}
Let $v_t={\cal F}_*(\frac{d}{dt})$. Since $\{\omega_t\}_{t\in [0,1]}$ is a cohomologous family, under the assumption of the theorem, we may choose a suitable family of 1-forms $\{\alpha_t\}_{t\in [0,1]}$ such that\\
\[
\frac{d\omega_t}{dt} = -d\alpha_t
\]\\
and $\alpha_t|_{S_t}$ is cohomologous to $i(v_t)\omega_t|_{S_t}$. Namely\\
\[
i(v_t)\omega_t|_{S_t} = (\alpha_t + da_t)|_{S_t}.
\]\\
$a_t$ can be easily extended to functions on $X$. Replace $\alpha_t$ by $\alpha_t + da_t$, we may assume\\
\[
\frac{d\omega_t}{dt} = -d\alpha_t
\]\\
and $\alpha_t|_{S_t}=i(v_t)\omega_t|_{S_t}$. To extend our flow, we need to extend $v_t$ to $X$ from $S_t$. For this purpose we need to find a function $f_t$ that satisfies\\
\[
i(v_t)\omega_t|_{T^*_{S_t}X} = (\alpha_t + df_t)|_{T^*_{S_t}X}.
\]\\
Then $v_t$ defined by\\
\[
i(v_t)\omega_t = \alpha_t + df_t
\]\\
will be our desired extension. This can be done similarly as in the proof of theorem \ref{eh} with the help of lemma \ref{ei}.
\begin{flushright} $\Box$ \end{flushright}

\subsection{The smooth case}
In all these results, one might be tempted to try to remove the $C^{0,1}$ conditions and get $C^\infty$ flows and vector fields. This is not really possible in general because of local symplectic obstructions. Simply put, there are more than one orbits of the action of the symplectic group on the space of configurations of normal crossing symplectic subspaces in a symplectic linear space. In particular, two orthogonal symplectic subspaces can not be mapped to two non-orthogonal symplectic subspaces under symplectic transformation. It turns out that these local obstructions are the only obstruction for the smoothness of the flow. We will concentrate on the case of toroidal symplectic form or equivalently, when normal crossing submanifolds intersect orthogonally with respect to the symplectic form (which will be called symplectic normal crossing to distinguish from the concept of normal crossing for submanifolds in differential topology). The general case can be proved by the same method.\\

\begin{de}
Let $V_1,V_2$ be symplectic subspaces of the symplectic vector space $(V,\omega)$. $V_1$ is {\bf orthogonal} to $V_2$ if $V_1^\perp \subset V_2$ or $V_1^\perp = V_2$ or $V_1^\perp \supset V_2$.
\end{de}

\begin{lm}
Assume that we have $(S^1\subset S^2 \subset U_{S^2} \subset X, \omega)$, $(\tilde{S}^1 \subset X, \omega)$, where $S^2$ is a symplectic submanifold in $X$, $\tilde{S}^1$ and $S^1 = \tilde{S}^1\cap S^2$ are symplectic normal crossing unions of symplectic submanifolds in $X$, and locally $\tilde{S}^1 \subset \pi^{-1}(S^1)$ and $\tilde{S}^1$ is orthogonal to $S^2$. Assume further that we have a vector field $w$ and a closed 1-form $\alpha$ such that\\
\[
i(v)\omega = \alpha\ \ {\rm and}\ \ v|_{\tilde{S}^1} = w|_{\tilde{S}^1},
\]\\
and $w$ preserves the orthogonal relation between $\tilde{S}^1$ and $S^2$. Then $\alpha$ can be modified by adding a smooth exact 1-form such that the modified $v$ is unchanged on $\tilde{S}^1$ and $(v-w)|_{S^2}$ is along $S^2$.\\
\end{lm}
{\bf Proof:} Compare to lemma \ref{ed}, the additional toroidal assumption implies that $b(x)y$ vanishes to second order along $\tilde{S}^1$. (In particular, if $\tilde{S}^1 = \pi^{-1}(S^1)$, then $b(x)$ vanishes to the second order along $\tilde{S}^1 = \pi^{-1}(S^1)$.)\\

By assumption, $w$ preserves the orthogonal relation between $\tilde{S}^1$ and $S^2$. Since $v$ is symplectic, it also preserves the orthogonal relation between $\tilde{S}^1$ and $S^2$. $v|_{\tilde{S}^1} = w|_{\tilde{S}^1}$ implies that $v-w$ ($v-w-u$) fixes $\tilde{S}^1$ and keeps $S^2$ orthogonal to $\tilde{S}^1$. For convention on $v,w,u$, please see the proof of lemma \ref{ed}.\\

To see the second order vanishing, refine the coordinate $(x,y) = (x,y_1,y_2)$ such that $S^2 = \{y=0\}$, $V_2=\{x=x_0,y_2=0\}$ is the fiber of $\pi: \tilde{S}^1 \rightarrow S^1$ over $(x_0,0)\in S^1$, $V_1=\{y_1=0\}$ is the orthogonal compliment of $V_2$. Correspondingly we can also decompose $b(x) = (b_1(x),b_2(x))$. $v-w-u$ keeping $S^2$ orthogonal to $\tilde{S}^1$ implies that the image of $S^2$ under the flow of $v-w-u$ will always be tangent to $V_1$ at $S^1$. Therefore, $b_1(x)$ vanishes to the second order and $b_2(x)$ vanishes to the first order along $S^1$. Consequently, $b(x)y = b_1(x)y_1 + b_2(x)y_2$ vanishes to the second order along $\tilde{S}^1 \subset \pi^{-1}(S^1)\cap \{y_2=0\}$.\\

Therefore we may define

\[
f = b(x)y c_1(y/\epsilon),
\]

where $\epsilon$ is a small constant. Since $b(x)y$ vanishes to second order along $\tilde{S}^1$, replacing $\alpha$ by $\alpha + d f$, we still have that $v|_{\tilde{S}^1}$ is unchanged and all the rest. Now clearly $f$ is a smooth function.
\begin{flushright} $\Box$ \end{flushright}
{\bf Remark:} It is easy to see that the corresponding modification of $w$ as in the remark after lemma \ref{ed} will keep $w$ smooth in this case.\\

Using similar induction argument as in the proof of theorem \ref{eh}, this lemma gives us the $C^\infty$ version of our theorem.\\
\begin{theorem}
Assume that two symplectic submanifold structures $(S_0\subset X, \omega_0)$ and $(S_1\subset X, \omega_1)$ are symplectic isotopic, where $S_t = {\cal F}_t(S)$ is a union of smooth symplectic submanifolds of {\bf symplectic} normal crossing. Then there exists a $C^\infty$ symplectomorphism $h: (X,\omega_0) \rightarrow (X,\omega_1)$ such that $h(S_0)=S_1$.
\end{theorem}
\begin{flushright} $\Box$ \end{flushright}

\subsection{Special cases}
There are two special cases of our result that are particularly useful. The first case is when the symplectic form is fixed. As a special case of theorem \ref{eh}, we have\\
\begin{theorem}
\label{ee}
Assume that two symplectic submanifold structures $S_0, S_1 \subset (X,\omega)$ are symplectic isotopic, where $S$ is a union of smooth manifolds of normal crossing. Then there exists a piecewise smooth $C^{0,1}$ Hamiltonian automorphism $h: (X,\omega) \rightarrow (X,\omega)$ such that $h(S_0)=S_1$.
\end{theorem}
\begin{flushright} $\Box$ \end{flushright}
{\bf Remark:} From the proof of theorem \ref{eh}, it is not hard to see that actually $h$ can be made identity in the region away from $\bigcup_t S_t$. The interesting thing is that $h$ can also be made identity in the region of $S$ that is not moved by the original flow $\{{\cal F}_t\}_{t\in [0,1]}$. This can be achieved by modifying $w$ more carefully in the remark after lemma \ref{ed} to ensure that $w$ is unchanged in the region where $w=0$.\\\\
The second case is when the submanifold $S_t$ is fixed. As another special case of theorem \ref{eh}, we have\\
\begin{theorem}
\label{eg}
Assume that $S$ is a union of normal crossing smooth symplectic submanifolds of $(X,\omega_t)$ for all $t$. Then there exist a piecewise smooth $C^{0,1}$ symplectomorphism $h: (X,\omega_0) \rightarrow (X,\omega_1)$ such that $h(S)=S$.
\end{theorem}
\begin{flushright} $\Box$ \end{flushright}
{\bf Remark:} The proof of theorem actually can ensure that $h$ is identity in the region where $\omega_t$ is unchanged. Again, one need to be careful when modifying the vector field $w$.\\
\begin{lm}
Assume that we have $(S \subset X, \omega)$, and $v$ is a section of $T_SX$ such that the 1-form $i(v)\omega|_S$ is exact on $S$. Then $v$ can be extended to $X$ as a piecewise smooth $C^{0,1}$ Hamiltonian vector field supported near $S$.\\
\end{lm}
{\bf Proof:} Since $i(v)\omega|_S$ is exact on $S$, we can find function $h$ supported near $S$ such that $i(v)\omega|_S = dh|_S$.
We need to find a function $f$ that satisfies\\
\[
i(v)\omega|_{T^*_{S}X} = d(h + f)|_{T^*_{S}X}.
\]\\
Let us first look at the problem locally. Assume that we have local coordinate $(x,y)$ of $U_{S}$ such that $\pi(x,y) = x$. Assume that\\
\[
i(v)\omega - dh = a(x)dx + b(x)dy + O(|y|).
\]\\
$i(v)\omega|_S = dh|_S$ implies that $a(x)=0$. We may define\\
\[
f = b(x)y c_1(y/\rho(x)),
\]\\
where $\epsilon$ is a small positive constant and $\rho(x)\leq \epsilon$ is the distance function to $S^1$ on $S^2$ when $\rho(x)\leq \epsilon/2$. It is easy to see that $f$ is of $C^{1,1}$.\\\\
Globally, one can define $f$ by using partition of unity to piece together $f$ defined in each local chart (as in the proof of theorem \ref{ea}). Then $v$ can be extended according to $i(v)\omega = d(h + f)$ and $v$ is surported near $S$.
\begin{flushright} $\Box$ \end{flushright}
With this lemma, the following extension theorem is immediate.\\
\begin{theorem}
\label{ef}
Assume that $\{{\cal F}_t\}_{t\in [0,1]}$ is a piecewise smooth $C^{0,1}$ Hamiltonian flow on the family $\{S_t\}_{t\in [0,1]}$, then there exists a piecewise smooth $C^{0,1}$ Hamiltonian morphism $h: (X,\omega) \rightarrow (X,\omega)$ such that $h(S_0)=S_1$ and $h|_{S_0} = {\cal F}_1\circ {\cal F}_0^{-1}$.
\end{theorem}
\begin{flushright} $\Box$ \end{flushright}
Although our purpose in this paper is to deform our Lagrangian fibration to a desired one, where we need the above results on deformation of symplectic submanifolds, classification of symplectic submanifolds via Hamiltonian deformation is very interesting in its own right within symplectic geometry. For example, the above results imply the following well known result.\\
\begin{co}
Assume that two complex submanifolds $S_1$, $S_2$ of a K\"{a}hler manifold $(X, \omega_g)$ are complex deformation equivalent. Then there exists a Hamiltonian diffeomorphism $h$ of $X$ such that $h(S_1) = S_2$.
\end{co}
\begin{flushright} $\Box$ \end{flushright}

\subsection{The piecewise smooth case}
For our application, we need to deal with the situation that the ambient symplectic manifold $(X,\omega)$ is smooth while the symplectic submanifold $S$ is piecewise smooth. More precisely, $S = \displaystyle\bigcup_{k=1}^l S^{[k]}$ is a union of smooth symplectic submanifolds with boundaries and corners, by which we assume that each $S^{[k]}$ is inside a open symplectic manifold of the same dimension as a symplectic submanifold in $(X,\omega)$. (By corners, we mean real normal crossing singularities for boundary.) It is straightforward to see that all of the results in previous subsections have analogous statements in the piecewise smooth case. For simplicity, we will only prove the special result we need.

\begin{lm}
\label{em}
Assume $Y$ is a manifold with smooth boundary. Let $\alpha$ be a smooth 1-form on $Y$ such that $\alpha|_{\partial Y}$ is exact on $\partial Y$. Then there exists a smooth function $f$ on $Y$ such that $\alpha = df$ along $\partial Y$.
\end{lm}
{\bf Proof:} Since $\alpha|_{\partial Y}$ is exact on $\partial Y$, there exists a function $f_1$ on $Y$ such that $\alpha|_{\partial Y} = df_1|_{\partial Y}$. Take a partition of unity $\{\rho_\beta\}$ along $\partial Y$ with respect to the open covering $\{U_\beta\}$ of the tubular neighborhood of $\partial Y$, by which we mean the support of $\rho_\beta$ is in $U_\beta$ for each $\beta$ and $(\sum_\beta \rho_\beta)|_{\partial Y}=1$. Let $x^\beta$ be the defining function of $\partial Y$ in $U_\beta$. Then $\alpha|_{\partial Y} = df_1|_{\partial Y}$ implies that

\[
(\alpha - df_1)|_{U_\beta} = f_2^\beta dx^\beta + O(|x^\beta|).
\]

Then it is easy to see

\[
f = f_1 + \sum_\beta \rho_\beta f_2^\beta x^\beta
\]

will satisfy our need.
\begin{flushright} $\Box$ \end{flushright}
\begin{lm}
\label{en}
In lemma \ref{em}, if we allow $\partial Y$ to have corners, then there exists a $C^{1,1}$ function $f$ on $Y$ such that $\alpha = df$ along $\partial Y$.
\end{lm}
{\bf Proof:} We will follow the proof of lemma \ref{em} closely and only mention the part that need modification. Assume that locally in $U_\beta$, $\partial Y \cap U_\beta = \{x^\beta_1x^\beta_2 =0,x^\beta_1\geq 0,x^\beta_2\geq 0\}$. (The generalization to corner of higher codimension is very straightforward.) Then $\alpha|_{\partial Y} = df_1|_{\partial Y}$ implies that $(\alpha - df_1)|_{U_\beta}$ will contain the term like $f_3^\beta x^\beta_1dx^\beta_2$. We will make the contribution of such term to $f$ to be $\rho_\beta f_3^\beta x^\beta_1 x^\beta_2 c_1(x^\beta_2/x^\beta_1)$, which is $C^{1,1}$. Then it is easy to see that the resulting $f$ will satisfy our need.
\begin{flushright} $\Box$ \end{flushright}
\begin{lm}
\label{eo}
Assume $Y$ is a symplectic submanifold with boundaries and corners in $(X,\omega)$. Let $\alpha$ be a smooth section of $T^*_X|_Y$ such that $\alpha|_{\partial Y}$ is exact on $\partial Y$. Then there exists a $C^{1,1}$ function $f$ on $X$ supported near $Y$ such that $\alpha = df$ when restricted to $T^*_X|_{\partial Y}$. If $\alpha$ as a section of $T^*_X|_Y$ vanishes on $(\partial Y)_0$, which is a connected subset in $\partial Y$, then $f$ can be made to vanish in the set

\[
\{p\in X|{\rm Dist}(p,Y) \geq \min(\epsilon_1, \epsilon_2 {\rm Dist}(p,(\partial Y)_0))\},
\]
where $\epsilon_1,\epsilon_2 >0$ are small.
\end{lm}
{\bf Proof:} As we mentioned earlier, $Y$ can be extended somewhat as symplectic submanifold beyond the boundary. By lemma \ref{en}, there exists a $C^{1,1}$ function $\hat{f}$ on $Y$ such that $\alpha|_Y = d\hat{f}|_Y$ along $\partial Y$. Extend $\hat{f}$ as a $C^{1,1}$ function on $X$ supported near $Y$. Take a partition of unity $\{\rho_\beta\}$ along $\partial Y$ with respect to the open covering $\{U_\beta\}$ of the tubular neighborhood of $\partial Y$. Assume $(x^\beta,y^\beta)$ is a local coordinate of $U_\beta$ such that $x^\beta$ is coordinate on $Y$ and locally $Y \cap U_\beta =\{y^\beta=0, x^\beta_1 \geq 0\}$. Then locally

\[
\alpha - d\hat{f} = b^\beta(x^\beta)dy^\beta + O(|y^\beta|) + O(|x^\beta_1|).
\]

We get the desired function

\[
f = \hat{f} + \sum_\beta \rho_\beta b^\beta(x^\beta)y^\beta.
\]

If $\alpha$ as a section of $T^*_X|_Y$ vanishes on $(\partial Y)_0$, then $d\hat{f}$ as a section of $T^*_X|_Y$ vanishes on $(\partial Y)_0$ and $b^\beta(x^\beta)$ vanishes on $(\partial Y)_0 \cap U_\beta$. Since $(\partial Y)_0$ is connected, we can adjust $\hat{f}$ by a constant so that $\hat{f}$ vanishes to the second order along $(\partial Y)_0$. $b^\beta(x^\beta)y^\beta$ also vanishes to the second order along $(\partial Y)_0 \cap U_\beta$. Let $\rho(p)$ be a smoothing of

\[
c_1\left(\frac{{\rm Dist}(p,Y)}{2\min(\epsilon_1, \epsilon_2 {\rm Dist}(p,(\partial Y)_0))}\right)
\]

away from $(\partial Y)_0$. It is easy to see that $\rho f$ is $C^{1,1}$. Replace $f$ by $\rho f$, we get the $C^{1,1}$ function with the desired vanishing condition.
\begin{flushright} $\Box$ \end{flushright}
Recall that a piecewise smooth symplectic submanifold $S = \displaystyle\bigcup_{k=1}^l S^{[k]}$ in $(X,\omega)$ is a union of smooth symplectic submanifolds with boundaries and corners in $(X,\omega)$. A section $\alpha$ of $T^*_X|_S$ is piecewise smooth if $\alpha$ is smooth when restricted to $T^*_X|_{S^{[k]}}$. $\alpha|_{\partial S}$ is said to be exact on $\partial S = \displaystyle\bigcup_{k=1}^l \partial S^{[k]}$ if there exists a piecewise smooth continuous function $f_{\partial}$ on $\partial S$ such that $\alpha|_{\partial S} = df_{\partial}$.

\begin{lm}
For a piecewise smooth section $\alpha$ of $T^*_X|_S$, assume that $\partial S^{[k]}$, $\displaystyle\left(\bigcup_{i=1}^{k-1} \partial S^{[i]}\right) \cap \partial S^{[k]}$ are connected and $\alpha|_{\partial S^{[k]}}$ is exact on $\partial S^{[k]}$ for each $k$. Then $\alpha|_{\partial S}$ is exact on $\partial S = \displaystyle\bigcup_{k=1}^l \partial S^{[k]}$.
\end{lm}
\begin{flushright} $\Box$ \end{flushright}
\begin{lm}
\label{ep}
Assume $S = \displaystyle\bigcup_{k=1}^l S^{[k]}$ is a piecewise smooth symplectic submanifold in $(X,\omega)$. Let $\alpha$ be a piecewise smooth section of $T^*_X|_S$ such that $\alpha|_{\partial S}$ is exact on $\partial S$, where $\partial S = \displaystyle\bigcup_{k=1}^l \partial S^{[k]}$. Then there exists a $C^{1,1}$ function $f$ on $X$ supported near $S$ such that $\alpha = df$ when restricted to $T^*_X|_{\partial S}$. If $\alpha$ as a section of $T^*_X|_S$ vanishes on $(\partial S)_0$, which is a connected subset in $\partial S$, then $f$ can be made to vanish in the set

\[
\{p\in X|{\rm Dist}(p,S) \geq \min(\epsilon_1, \epsilon_2 {\rm Dist}(p,(\partial S)_0))\},
\]
where $\epsilon_1,\epsilon_2 >0$ are small.
\end{lm}
{\bf Proof:} According to lemma \ref{eo}, there exists a $C^{1,1}$ function $f_{[k]}$ on $X$ supported near $S^{[k]}$ such that $\alpha = df_{[k]}$ when restricted to $T^*_X|_{\partial S^{[k]}}$ for each $k$. We may further adjust, so that $f_{[k]}|_{\partial S^{[k]}} = f_{\partial}|_{\partial S^{[k]}}$ for each $k$.\\

It is easy to constract a tubular neighborhood $U_\epsilon$ of $S$ with projection $\pi: U_\epsilon \rightarrow S$, so that the fibers of $\pi$ vary smoothly and $f_{[k]}$ is supported in $U_\epsilon$. Piece together $f_{[k]}|_{\pi^{-1}(S^{[k]})}$, we get a discontinuous function

\[
\hat{f} = \sum_{k=1}^lf_{[k]}\chi_{\pi^{-1}(S^{[k]})}
\]

supported in $U_\epsilon$, satisfying $\alpha = d\hat{f}$ when restricted to $T^*_X|_{\partial S}$. The variation functions of the discontinuous function $\hat{f}$ and its derivatives are supported in $\pi^{-1}(\partial S)$ and satisfy

\[
{\rm Var}[\hat{f}](p) \leq C{\rm Dist}^2(p,\partial S),\ \ C{\rm Var}[D\hat{f}](p) \leq {\rm Dist}(p,\partial S),\ \ {\rm Var}[D^2\hat{f}](p) \leq C.
\]

For each $k$, it is straightforward to construct a function $\hat{\rho}_{[k]}\geq 0$ on $U_\epsilon$ that is smooth away from $\partial S^{[k]}$ satisfying

\[
\hat{\rho}_{[k]}|_{\pi^{-1}(S^{[k]})} = 1,\ \ \hat{\rho}_{[k]}(p) = 0\ {\rm when}\ {\rm Dist}(p,S) < {\rm Dist}(p,\pi^{-1}(S^{[k]})),
\]

\[
|D\hat{\rho}_{[k]}(p)|\leq {\rm Dist}^{-1}(p,\partial S^{[k]}),\ \ |D^2\hat{\rho}_{[k]}(p)|\leq {\rm Dist}^{-2}(p,\partial S^{[k]}).
\]

Normalize $\hat{\rho}_{[k]}$, we get

\[
\rho_{[k]} = \hat{\rho}_{[k]}\left(\sum_{i=1}^l\hat{\rho}_{[i]}\right)^{-1},\ \ \sum_{k=1}^l\rho_{[k]} = 1\ \ {\rm on}\ U_\epsilon.
\]

Define

\[
f = \sum_{k=1}^lf_{[k]}\rho_{[k]}.
\]

$f$ is a smoothing of $\hat{f}$. To verify that $f$ is $C^{1,1}$, it is sufficient to show that $|Df|$ and $|D^2f|$ are bounded.

\[
|Df| \leq \sum_{k=1}^l|Df_{[k]}|\rho_{[k]} + C{\rm Var}[\hat{f}]\max_k |D\rho_{[k]}| \leq C.
\]

\[
|D^2f| \leq \sum_{k=1}^l|D^2f_{[k]}|\rho_{[k]} + C{\rm Var}[D\hat{f}]\max_k |D\rho_{[k]}| + C{\rm Var}[\hat{f}]\max_k |D^2\rho_{[k]}| \leq C.
\]

It is easy to verify that $df = d\hat{f}$ when restricted to $T^*_X|_S$. Therefore $\alpha = d\hat{f} = df$ when restricted to $T^*_X|_{\partial S}$. The vanishing condition will be satisfied as direct consequence of the vanishing condition in lemma \ref{eo}.
\begin{flushright} $\Box$ \end{flushright}
\begin{theorem}
\label{eq}
Assume that two piecewise smooth symplectic submanifold structures $(S_0\subset X, \omega)$ and $(S_1\subset X, \omega)$ are symplectic isotopic with $\alpha_t|_{\partial S_t}$ being exact, where $\alpha_t = i({\cal F}_{*t}(\frac{d}{dt}))\omega$, then there exists a piecewise smooth $C^{0,1}$ symplectomorphism $h: (X, \omega) \rightarrow (X, \omega)$ so that $h(S_0)=S_1$. If for all $t$, ${\cal F}_t \circ {\cal F}_0^{-1}$ retricts to identity on $(\partial S_0)_0$, which is a connected subset in $\partial S_0$, then $h$ can be made to retrict to identity on

\[
\{p\in X|\min_t({\rm Dist}(p,S_t)) \geq \min(\epsilon_1, \epsilon_2 {\rm Dist}(p,(\partial S_0)_0))\},
\]
where $\epsilon_1,\epsilon_2 >0$ are small.
\end{theorem}
{\bf Proof:}
Since $\alpha|_{\partial S_t}$ is exact, by lemma \ref{ep}, there exists a smooth familly of $C^{1,1}$ functions $\{\hat{f}_t\}$ on $X$ supported near $\{S_t\}$ such that $\alpha_t - d\hat{f}_t$ vanishes when restrict to $T^*_X|_{\partial S_t}$. If ${\cal F}_t \circ {\cal F}_0^{-1}$ retricts to identity on $(\partial S_0)_0$ for all $t$, then $\alpha_t$ as a section of $T^*_X|_{S_t}$ vanishes on $(\partial S_0)_0$. By lemma \ref{ep}, $\hat{f}_t$ can be made to vanish on the set

\[
\{p\in X|{\rm Dist}(p,S_0) \geq \min(\epsilon_1, \epsilon_2 {\rm Dist}(p,(\partial S_0)_0))\}.
\]

Let $w_t={\cal F}_*(\frac{d}{dt})$. Since $S_t$ is symplectic, we can find $u_t\in TS_t$ such that\\
\[
i(u_t)\omega|_{S_t} = i(w_t)\omega|_{S_t} - d\hat{f}_t|_{S_t}.
\]
Since $\alpha_t - d\hat{f}_t$ vanishes along $\partial S_t$, we have $u_t|_{\partial S_t} = 0$. Let $v_t=w_t - u_t$ then\\
\[
i(v_t)\omega|_{S_t} = d\hat{f}_t,\ \ v_t|_{\partial S_t}=w_t|_{\partial S_t}.
\]
The flow of $v_t$ will map $(S_0,\omega)$ to $(S_t,\omega)$ (in particular $(S_1,\omega)$) symplectically. To get a globally defined flow, we need to extend $v_t$ to $X$ from $S_t$. For this purpose we need a function $f_t$ on $X$ (extending $\hat{f}_t$ on $S_t$) that satisfies\\
\[
i(v_t)\omega|_{T^*_{S_t}X} = (df_t)|_{T^*_{S_t}X}.
\]\\
Then $v_t$ defined by\\
\[
i(v_t)\omega = df_t
\]\\
will be our desired extension. Since\\
\[
i(v_t)\omega|_{S_t} = d\hat{f}_t|_{S_t}.
\]\\
Locally on $S_t$, assume $x,y$ are tangent and normal coordinates near $S_t$. We may express\\
\[
i(v_t)\omega - d\hat{f}_t= b(x)dy + O(|y|),
\]\\
where $b(x)=0$ for $x\in \partial S_t$. We may define\\
\[
f_t = \hat{f}_t + b(x)y c_1(y/\epsilon),
\]\\
where $\epsilon$ is a small positive constant. Clearly $f_t \in C^{1,1}$ and $(df_t - d\hat{f}_t)|_{T^*_{\partial S_t}X} =0$. Globally we may use partition of unity to piece such local $f_t$ together to form a function $f_t$ on $X$ that is supported near $S_t$. The flow of corresponding $v_t$ (that satisfies $i(v_t)\omega = df_t$ and $v_t|_{\partial S_t} = w_t|_{\partial S_t}$) will produce the desired piecewise smooth $C^{0,1}$ symplectomorphism $h$. For the last statement of the theorem, it is necessary to use the more refined

\[
f_t = \hat{f}_t + b(x)y c_1\left(\frac{{\rm Dist}(p,S_t)}{2\min(\epsilon_1, \epsilon_2 {\rm Dist}(p,(\partial S_0)_0))}\right).
\]

It is easy to check that $v_t$ vanishes on

\[
\{p\in X|{\rm Dist}(p,S_t) \geq \min(\epsilon_1, \epsilon_2 {\rm Dist}(p,(\partial S_0)_0))\}.
\]

Consequently, $h$ retricts to identity on

\[
\{p\in X|\min_t({\rm Dist}(p,S_t)) \geq \min(\epsilon_1, \epsilon_2 {\rm Dist}(p,(\partial S_0)_0))\}.
\]

\begin{flushright} $\Box$ \end{flushright}
There are two corollaries that are particularly useful.

\begin{co}
Assume that two piecewise smooth symplectic submanifold structures $(S_0\subset X, \omega)$ and $(S_1\subset X, \omega)$ are symplectic isotopic with $\partial S_t$ fixed, then there exists a piecewise smooth $C^{0,1}$ symplectomorphism $h: (X, \omega) \rightarrow (X, \omega)$ so that $h(S_0)=S_1$ and $h|_{\partial S_0} = id$.
\end{co}
\begin{flushright} $\Box$ \end{flushright}
\begin{co}
\label{er}
Assume that two piecewise smooth symplectic submanifold structures $(S_0\subset X, \omega)$ and $(S_1\subset X, \omega)$ are symplectic isotopic with $\dim_\mathbb{R}S = 2$ and the area of each component of $S_t = {\cal F}_t(S)$ is constant when $t$ varies, then there exists a piecewise smooth $C^{0,1}$ symplectomorphism $h: (X, \omega) \rightarrow (X, \omega)$ so that $h(S_0)=S_1$. If ${\cal F}_t \circ {\cal F}_0^{-1}$ retricts to identity on $(\partial S_0)_0$ for all $t$, which is a connected subset in $\partial S_0$, then $h$ can be made to retrict to identity on

\[
\{p\in X|\min_t({\rm Dist}(p,S_t)) \geq \min(\epsilon_1, \epsilon_2 {\rm Dist}(p,(\partial S_0)_0))\},
\]
where $\epsilon_1,\epsilon_2 >0$ are small.
\end{co}
{\bf Proof:}
Without loss of generality, assume $S$ is smooth. For $\alpha = i({\cal F}_*(\frac{d}{dt}))\omega$, $d\alpha = {\cal L}_{{\cal F}_*(\frac{d}{dt})}\omega$. Since $\dim_\mathbb{R}\partial S = 1$, $\alpha|_{\partial S_t}$ is automatically closed.

\[
\int_{\partial S_t} \alpha = \int_{\partial S} {\cal F}_t^* \alpha = \int_S {\cal F}_t^*d\alpha = \int_S {\cal L}_{{\cal F}_*(\frac{d}{dt})}\omega  = \int_S \frac{d}{dt}\left({\cal F}_t^*\omega\right) = \frac{d}{dt}\int_{S_t} \omega = 0
\]

implies that $\alpha|_{\partial S_t}$ is exact. By theorem \ref{eq}, we get our conclusion.
\begin{flushright} $\Box$ \end{flushright}

\se{Toroidal manifold and toroidal \k metric}
To ensure that our gradient flow is well behaved, it is necessary to use the so-called toroidal \k metrics for the gradient flow. In this section, we will discuss the construction of toroidal \k metric as small perturbation of arbitrary \k metric.\\

Let $(X, D, \omega_g)$ be a \k manifold $X$ with a divisor $D$ and a \k metric $g$ with corresponding \k form $\omega_g$. $(X,D)$ is called {\bf toroidal} if at every point $p\in D$ there exists a local chart $i_p: V_p \rightarrow P_\Sigma$, where $P_\Sigma$ is a toric variety with the big torus $T \subset P_\Sigma$, such that $i_p^{-1}(T) = V_p \cap (X\backslash D)$. A local coordinate on $V_p$ is called {\bf toroidal coordinate} if each coordinate function is a pullback by $i_p$ of a toric monoidal function on $P_\Sigma$.\\

Our gradient flow approach more generally can be applied to toroidal situations. To ensure good behavior of the gradient flow in such situations, it is necessary to adopt suitable \k metrics that are compatible with toroidal structure. (There are counter-examples for general metrics.) In this work, we will mainly concern the case of normal crossing $D$.\\
\begin{lm}
If $D$ is of normal crossing, then $(X,D)$ is toroidal.
\end{lm}
\begin{flushright} $\Box$ \end{flushright}
From now on, in this section, we always assume that $D$ is of normal crossing.\\

{\bf Remark:} Assume $\displaystyle D = \bigcup_{i=1}^l D_i$ is a normal crossing divisor. A coordinate $z = (z_1,\cdots, z_n)$ of a neighborhood $V_p$ of $p\in D$ that satisfies $D_i \cap V_p = \{z_i =0\}$ will usually be refered to as a toroidal coordinate in the sense that $z$ is a toroidal coordinate with respect to the toroidal structure defined by $i_p (= z): V_p \rightarrow {\mathbb{C}^n}$.\\

\begin{de}
A \k form $\omega$ is called {\bf toroidal} (along $S\subset D$) with respect to $(X,D)$ if different components of $D$ intersect orthogonally (at $S\subset D$) with respect to $\omega$.\\
\end{de}

{\bf Remark:} The concept of toroidal defined above probably should be called infinitesimally toroidal in more detailed discussion. However, this concept of toroidal \k form is enough for our gradient flow to work.\\

Assume $\displaystyle D = \bigcup_{i=1}^l D_i$ and each $D_i$ is a smooth divisor. For each index set $I\in \{1,\cdots,l\}$ define $\displaystyle D_I=\bigcap_{i\in I} D_i$ when the intersection is non-empty. $D_I$ so defined is an $(n-|I|)$-submanifold in $X$. Let

\[
D^{(k)} = \bigcup_{|I|=n-k} D_I,\ \ D^{(k)}_0 = D^{(k)}\backslash D^{(k-1)}.
\]

Then we have the filtration $D = D^{(n-1)} \supset \cdots \supset D^{(1)} \supset D^{(0)}$ and

\[
D = \bigcup_{i=0}^{n-1} D^{(i)}_0.
\]

Let $U_i$ denote a tubular neighborhood of $D_i$ then

\[
U_I=\bigcap_{i\in I} U_i,\ \ U^{(k)} = \bigcup_{|I|=n-k} U_I
\]

form tubular neighborhoods of $D_I$, $D^{(k)}$. In particular, $U = U^{(n-1)}$ is a tubular neighborhood of $D$.\\

It is not hard to construct a suitable open covering $\{U^\alpha\}_{\alpha\in A}$ of $U$ satisfying the following properties. For any $\alpha$, there exists a unique index set $I_\alpha$ such that $U^\alpha\cap D_{I_\alpha} \not= \emptyset$ and $U^\alpha\cap D_{J} = \emptyset$ when $|J|>|I_\alpha|$, and $U^\alpha$ is a tubular neighborhood of $U^\alpha\cap D_{I_\alpha}$ with coordinate $(w,z)$. More precisely, $w$ is the coordinate on  $U^\alpha\cap D_{I_\alpha}$ and $(w,z)$ defines an identification of $U^\alpha$ with the product of $U^\alpha\cap D_{I_\alpha}$ and an open neighborhood in $\mathbb{C}^{|I|}$. $w$ defines a holomorphic fibration $\pi_\alpha: U^\alpha \rightarrow U^\alpha\cap D_{I_\alpha}$ and $z$ is the coordinate on fibers.\\

For any index set $I$, let $A_I = \{\alpha\in A|I_\alpha \supset I\}$. Then $\{U^\alpha\}_{\alpha\in A_I}$ forms a covering of $U_I$, and $\{U^\alpha\cap D_I\}_{\alpha\in A_I}$ forms a covering of $D_I$. Let $(w^\alpha,z^\alpha)$ be the coordinate of $U^\alpha$. (Here when $I_\alpha \not= I$, $(w^\alpha,z^\alpha)$ has to be modified so that $z^\alpha = (z^\alpha_i)_{i\in I}$ and $(z^\alpha_i)_{i\in I_\alpha\backslash I}$ is switched to be part of $w^\alpha$, $\pi_\alpha$ will also need to be adjusted accordingly.) $\pi_\alpha$ defines the local holomorphic fibration of $U^\alpha$ into $D_I$. We would like to combine these local holomorphic fibrations $\{\pi_\alpha\}_{\alpha\in A_I}$ into a global smooth fibration $\pi_I: U_I \rightarrow D_I$ with holomorphic fibres using a partition of unity $\{\rho_\alpha\}_{\alpha\in A_I}$ on $D_I$ with respect to the open covering $\{U^\alpha\cap D_I\}_{\alpha\in A_I}$. We would also like a holomorphic toroidal coordinate $z^I$ on each fibre $\pi_I^{-1}(x)$ of $\pi_I$ that vary smoothly on $x\in D_I$. We start with the following lemma that constructs smooth varying toroidal holomorphic local charts.\\

\begin{lm}
\label{fd}
For any $x\in D_I$, there exists a holomorphic coordinate $(w^I_x,z^I_x)$ in a neighborhood of $x$, where $z^I_x = (z^I_{x,i})_{i\in I}$, such that locally around $x$, $D_i = \{z^I_{x,i}=0\}$ for any $i\in I$, $w^I_x$ is holomorphic coordinate on $D_I = \{z^I_x=0\}$ and the family $\{(w^I_x,z^I_x)\}_{x\in D_I}$ varies smoothly when $x$ varies (modulo linear transformation on $w^I_x$ and non-zero multiple for $z^I_{x,i}$ ($i\in I$) that only depends on $x$).\\
\end{lm}
{\bf Proof:} Fix $x\in D_I$. For each $U^\alpha$ containing $x$, adjust $w^\alpha$ by an affine transformation that only depends on $x$ such that $w^\alpha(x)=0$ and $dw^\alpha$ are all equal at $x$ for all $\alpha$. (Namely $w^\alpha$ all agree to 1st order at $x$.) For $z^\alpha = (z^\alpha_i)_{i\in I}$, one may adjust each component $z^\alpha_i$ of $z^\alpha$ by a non-zero multiple that only depends on $x$ such that $dz^\alpha$ are all equal at $x$ for all $\alpha$. Then define

\[
w^I_x = \sum_{\alpha} \rho_\alpha(x)w^\alpha,\ \ z^I_x = \sum_{\alpha} \rho_\alpha(x)z^\alpha.
\]

$(w^I_x,z^I_x)$ so defined clearly satisfies the requirements of the lemma.
\begin{flushright} $\Box$ \end{flushright}
{\bf Remark:} It is easy to observe from our construction (with possible shrinking of tubular neiborhoods $\{U_I\}$ when necessary) that for $I\subset J$, $x\in D_J$ and $x'\in D_I$ near $x$ satisfying $x = w^J_x(x')$, $(w^I_{x'},z^I_{x'})$ and $(w^J_x,z^J_x)$ are compatible in the sense that $w^I_{x'}=(w^J_x,(z^J_{x,i})_{i\in J\backslash I})$ and $z^I_{x',i}/z^J_{x,i}$ ($i\in I$) are non-zero and only depend on $(z^J_{x,i})_{i\in J\backslash I}$. In particular, fibers of $w^J_x$ are unions of fibers of $w^I_{x'}$.\\

{\bf Remark on terminology:} The family $\{(w^I_x,z^I_x)\}_{x\in D_I}$ is said to vary smoothly when $x$ varies if there is a smooth varying family of small neighborhoods $\{U_x\}_{x\in D_I}$ such that $(w^I_x,z^I_x)$ define a smooth map from the open set $\displaystyle \bigcup_{x\in D_I} U_x \subset X \times D_I$ to $\mathbb{C}^n$. The family $\{(w^I_x,z^I_x)\}_{x\in D_I}$ is said to vary smoothly when $x$ varies (modulo linear transformation on $w^I_x$ and non-zero multiple for $z^I_{x,i}$ ($i\in I$) that only depend on $x$), if there is an open covering $\{U_\beta \}$ of $D_I$, on each $U_\beta$, there is a smooth varying family $\{(w^\beta_x,z^\beta_x)\}_{x\in U_\beta}$ such that for $x\in U_\beta \cap U_{\beta'}$, $w^\beta_x$, $w^{\beta'}_x$ and $w^I_x$ are related by linear transformations that only depend on $x$, and $z^\beta_{x,i}$, $z^{\beta'}_{x,i}$ and $z^I_{x,i}$ are related by non-zero multiples that only depend on $x$ for $i\in I$.\\

\begin{prop}
\label{fe}
For each index set $I$, one may construct a smooth fibration $\pi_I: U_I \rightarrow D_I$, whose fibres are holomorphic. On each fibre $\pi_I^{-1}(x)$, there is a holomorphic toroidal coordinate $z^I = (z^I_i)_{i\in I}$ that varies smoothly when $x\in D_I$ varies (modulo non-zero multiple for $z^I_i$ ($i\in I$) that only depend on $x$). One can also make such $\{(\pi_I,z^I)\}$ compatible in the sense that for $I\subset J$, fibers of $\pi_J$ are unions of fibers of $\pi_I$, moreover, $\pi_I$ restricted to each fiber of $\pi_J$ is a holomorphic fibration and for any $i\in I$, $z^J_i/z^I_i$ only depends on $z^J_j$ for $j\in J\backslash I$.\\
\end{prop}
{\bf Proof:} We will start by constructing the holomorphic fiber of $\pi_I$ over a point $x\in D_I$. Recall from lemma \ref{fd}, we have the local holomorphic coordinate $(w^I_x, z^I_x)$. $w^I_x$ defines a holomorphic fibration near $x$. Define $\pi_I^{-1}(x) = (w^I_x)^{-1}(x)$. Then clearly $\pi_I^{-1}(x)$ is holomophic and varies smoothly when $x$ changes. The toric coordinate on the fiber $\pi_I^{-1}(x)$ can be defined as

\[
z^I|_{\pi_I^{-1}(x)} = z^I_x|_{\pi_I^{-1}(x)}.
\]

The compatibility of $\{(\pi_I,z^I)\}_I$ is a direct consequence of the compatibility of $\{(w^I_{x},z^I_{x})\}_I$ discussed in the remark after lemma \ref{fd}.
\begin{flushright} $\Box$ \end{flushright}
The toroidal coordinate $z^I$ naturally determines a rank $|I|$ real torus $T^{|I|}_\mathbb{R}$-action on the fiber $\pi_I^{-1}(x)$ that varies smoothly when varying $x$. Together we get a smooth $T^{|I|}_\mathbb{R}$-action on $U_I$. These actions are compatible in the sense that for $I\subset J$, $T^{|I|}_\mathbb{R}$-action on $U_I$ restricted to $U_J$ is a subaction of $T^{|J|}_\mathbb{R}$-action on $U_J$.\\

{\bf Remark:} Our discussion so far does not involve \k form. When there is a \k form, it is desirable that the fibres of $\pi_I$ are normal to $D_I$. Recall the local holomorphic coordinate $(w^I_x,z^I_x)$ in lemma \ref{fd}. $w^I_x$ can be modified uniquely by linear functions on $z^I_x$ to ensure that for correspondingly modified $\pi_I$, $\pi_I^{-1}(x)$ is normal to $D_I$. Such modified $\{\pi_I\}$ are compatible (in the sense that for $I\subset J$ and any $x\in D_J\subset D_I$, $\pi_I^{-1}(x) = \{(w^J,z^J)\in \pi_J^{-1}(x)|(z^J_i)_{i\in J\backslash I}=0\}$) if the \k form is toroidal along $D_J$. The compatibility condition for the induced real torus actions is also weaker in this case.\\

To construct toroidal \k metrics, let us first look at some local constructions. Let $\omega_1$ and $\omega_2$ be two flat \k forms on $\mathbb{C}^n$. We are interested in constructing \k form $\omega$ such that $\omega=\omega_2$ near origin and $\omega=\omega_1$ away from a compact set around origin. Without loss of generality, we may assume that\\
\[
\omega_1 = \sum_{i=1}^n dz^i\wedge d\bar{z}^i = \sum_{i=1}^n \partial\bar{\partial} |z^i|^2,
\]
and\\
\[
\omega_2 = \sum_{i=1}^n \lambda_i dz^i\wedge d\bar{z}^i = \sum_{i=1}^n \partial\bar{\partial} \lambda_i |z^i|^2.
\]\\
Let\\
\[
\omega = \sum_{i=1}^n \partial\bar{\partial} h_i(z).
\]\\
Define two cut off functions $\rho$ and $\rho_c$ such that\\
\[
\rho(r)=1\ {\rm for}\ |r|\leq 1,\ \rho(r)=0\ {\rm for}\ |r|\geq 2,
\]
\[
\rho_c(r)=1\ {\rm for}\ |r|\leq 1,\ \rho_c(r)=0\ {\rm for}\ |r|\ {\rm large\ and\ } |\rho_c'|(r), |\rho_c''|(r)\leq c\ {\rm for\ all}\ r.
\]\\
Let $\hat{z}^i= z-z^ie_i$. Then we may take\\
\[
h_i(z) = \rho_\epsilon(|\hat{z}^i|/a)f_{\lambda_i}(|z^i|^2)+ (1- \rho_\epsilon(|\hat{z}^i|/a))|z^i|^2,
\]
where
\[
f_\lambda(r) = c_1 + \int_0^r\left(\lambda - c_2\int_0^s \rho\left(\frac{\log t-2\log a}{\lambda-1}\right)\frac{1}{2t}dt\right)ds,
\]
where\\
\[
c_2 = (\lambda-1)\left/\int_0^\infty\rho\left(\frac{\log t-2\log a}{\lambda-1}\right)\frac{1}{2t} dt\right..
\]\\
The constant $c_1$ should be chosen to ensure $f_\lambda(r) =r$ for $r$ large. Clearly $|c_2|\leq 1$ and $|c_1| \leq Ca^2$. $\epsilon$ should be taken to be sufficiently small.\\
\begin{lm}
\label{fc}
$\omega$ defined above is a \k form and satisfies:\\
\[
\omega=\omega_2 \ \ {\rm for}\ |z|\leq R_1 = \min_i (a_ie^{-|\lambda_i-1|}),
\]
\[
\omega=\omega_1 \ \ {\rm for}\ |z|\geq R_2 = \max_i (a_ie^{|\lambda_i-1|},a_i/\epsilon_i).
\]
\end{lm}
{\bf Proof:} The only non-trivial part is to verify that $\omega$ is a \k form. Notice that $f_\lambda(r)$ satisfies:\\
\[
f'_\lambda(r) = \left\{\begin{array}{ll}\lambda,\ \ &r\leq a^2e^{-2|\lambda-1|}\\{\rm monotone},\ \ &a^2e^{-2|\lambda-1|} \leq r\leq a^2e^{2|\lambda-1|}\\1,\ \ &r\geq a^2e^{2|\lambda-1|}\end{array}\right.
\]\\
\[
f_\lambda(r) = \left\{\begin{array}{ll}c_1 + \lambda r,\ \ &r\leq a^2e^{-2|\lambda-1|}\\r,\ \ &r\geq a^2e^{2|\lambda-1|}\end{array}\right.
\]\\
$f_\lambda(r) - r$ is supported in $r\leq a^2e^{2|\lambda-1|}$. We have\\
\[
\partial\bar{\partial} h_i(z) = \rho_\epsilon(|\hat{z}^i|/a)\partial\bar{\partial} f_{\lambda_i}(|z^i|^2)+ (1- \rho_\epsilon(|\hat{z}^i|/a))\partial\bar{\partial} |z^i|^2. + \rho'_\epsilon(|\hat{z}^i|/a)\beta_1 + \rho''_\epsilon(|\hat{z}^i|/a)\beta_2.
\]\\
It is not hard to see that $\beta_1,\beta_2$ are bounded with compact support.\\
\[
\partial\bar{\partial} f_{\lambda_i}(|z^i|^2) = (f'_{\lambda_i}(|z^i|^2) + |z^i|^2f''_{\lambda_i}(|z^i|^2))\partial\bar{\partial} |z^i|^2.
\]\\
When $\lambda_i \leq 1$, $f''_{\lambda_i}(r)\geq 0$ for all $r$. When $\lambda_i \geq 1$, $|rf''_{\lambda_i}(r)|\leq \frac{1}{2}$ for all $r$. Therefore\\
\[
f'_{\lambda_i}(|z^i|^2) + |z^i|^2f''_{\lambda_i}(|z^i|^2) \geq \min(\lambda_i,\frac{1}{2}).
\]\\
\[
\partial\bar{\partial} f_{\lambda_i}(|z^i|^2) \geq \min(\lambda_i,\frac{1}{2})\partial\bar{\partial} |z^i|^2 + O(\epsilon_i).
\]\\
When $|\epsilon|$ is sufficiently small, we have\\
\[
\omega = \sum_{i=1}^n \partial\bar{\partial} h_i(z) \geq (\min_i(\lambda_i,\frac{1}{2})- C|\epsilon|)\partial\bar{\partial} |z|^2
\]\\
is a \k form.
\begin{flushright} $\Box$ \end{flushright}
{\bf Remark:} Observe that $\omega$ can be written as\\
\[
\omega = \omega_1 + \partial\bar{\partial} R,
\]
where
\[
R = \sum_{i=1}^n \rho_\epsilon(|\hat{z}^i|/a)(f_{\lambda_i}(|z^i|^2) - |z^i|^2)
\]
has compact support.
\begin{flushright} $\Box$ \end{flushright}
The following corollary is a generalization of lemma \ref{fd} to non-flat \k metrics that we do not really need in this paper.

\begin{co}
\label{ff}
For a general (not necessarily flat) \k form $\tilde{\omega}_1 = \omega_1 + \omega_3$ with flat part $\omega_1$ and higher order term $\omega_3$, one can construct a \k form $\tilde{\omega}$ such that

\[
\tilde{\omega} = \omega_2\ {\rm for}\ |z|\leq R_1,\ \ {\rm and}\ \ \tilde{\omega} = \tilde{\omega}_1\ {\rm for}\ |z|\geq R_2.
\]
\end{co}

{\bf Proof:} Recall the $\omega$ constructed in lemma \ref{fc} satisfies

\[
\omega = \omega_2\ {\rm for}\ |z|\leq R_1,\ \ {\rm and}\ \ \omega = \omega_1\ {\rm for}\ |z|\geq R_2,\ \ {\rm and}\ \ \omega \geq c\omega_1\ (c>0).
\]

We can write $\omega_3 = \partial\bar{\partial} h_3$, where $h_3 = O(|z|^3)$. Define

\[
\tilde{\omega}_3 = \partial\bar{\partial} \left(1-\rho\left(\frac{2|z|}{R_2}\right)\right)h_3.
\]

It is easy to see that

\[
\tilde{\omega}_3 = O(|z|),\ \ \tilde{\omega}_3 = 0\ {\rm for}\ |z|\leq \frac{R_2}{2},\ \ {\rm and}\ \ \tilde{\omega}_3 = \omega_3\ {\rm for}\ |z|\geq R_2.
\]

Therefore $\tilde{\omega} = \omega + \tilde{\omega}_3$ will be \k when $R_2$ is small enough (which can be easily achieved in lemma \ref{fc}) and satisfies all the requirements.
\begin{flushright} $\Box$ \end{flushright}
{\bf Remark:} Observe that $\tilde{\omega}$ can be written as\\
\[
\tilde{\omega} = \tilde{\omega}_1 + \partial\bar{\partial} \tilde{R},
\]
where
\[
\tilde{R} = \sum_{i=1}^n \rho_\epsilon(|\hat{z}^i|/a)(f_{\lambda_i}(|z^i|^2) - |z^i|^2) - \rho\left(\frac{2|z|}{R_2}\right)h_3
\]
has compact support.
\begin{flushright} $\Box$ \end{flushright}
The construction of toroidal metrics will be done through induction on strata $\{D_I\}$ of $D$ starting from the lowest strata. During the induction process, a typical situation is $S^0 \subset S \subset D$, $S$ as a closed submanifold of X is a strata of $D$. $S^0$ is an open set in $S$ such that $D$ near $S^0$ is a product structure. $S\backslash S^0$ is a finite union of lower dimensional stratas of $D$. We need the following lemma\\
\begin{lm}
\label{fb}
Assume that $\omega_g$ is toroidal in a neighborhood of $S\backslash S^0$, then $\omega_g$ can be perturbed near $S^0$ so that $\omega_g$ is toroidal in a neighborhood of $S$.\\
\end{lm}
{\bf Proof:} Consider a normal neighborhood construction $\pi: U_S \rightarrow S$, where $U_S$ is a neighborhood of $S$ such that fibres of $\pi$ are holomorphic and intersect $S$ orthogonally with respect to the \k metric. (Here we are using the construction in the remark after proposition \ref{fe}. Since $\omega_g$ is toroidal in a neighborhood of $S\backslash S^0$, $\pi$ so constructed is compatible with previous fibration onto lower strata.) Let $w$ be holomorphic coordinate on $S$, $z$ be holomorphic toroidal coordinate on $\pi^{-1}(w)$ that depends smoothly on $w$. $\omega_g$ naturally determines a function $h$ such that $h$ restricted to each fibre $\pi^{-1}(w)$ is $(1,1)$-quadratic on $z$ and $\omega_g|_{\pi^{-1}(w)}=\partial \bar{\partial} h|_{\pi^{-1}(w)}$ at $\pi^{-1}(w)\cap S$. By assumption, $\omega_g|_{\pi^{-1}(w)}$ is toroidal along $S\backslash S^0$. Therefore, $\partial \bar{\partial} h|_{\pi^{-1}(w)}$ is toroidal along $S\backslash S^0$. We may choose $h_0$ that is $(1,1)$-quadratic on $z$ as toroidal extension of $h|_{\pi^{-1}(S\backslash S^0)}$ to $\pi^{-1}(S)$. (A canonical way to get $h_0$ is to consider the natural real torus action on $\pi^{-1}(w)$ determined by the toroidal coordinate $z$ and take $h_0$ to be the average function of $h$ with repect to the real torus action.)\\\\
For any $w\in S$, change coordinate $z$ on $\pi^{-1}(w)$ such that\\
\[
h|_{\pi^{-1}(w)} = \sum_{i=1}^n |z^i|^2,\ \ \ h_0|_{\pi^{-1}(w)} = \sum_{i=1}^n \lambda |z^i|^2.
\]
Define
\[
R|_{\pi^{-1}(w)} = \sum_{i=1}^n \rho_\epsilon\left(|\hat{z}^i|/a\right)\left(f_{\lambda_i}(|z^i|^2)- |z^i|^2\right).
\]
as in the previous lemma. Notice that since $z|_S=0$, $\frac{\partial z}{\partial w}|_S=0$. $w$-components of $\partial \bar{\partial} R$ will vanish along $S$. Then according to lemma \ref{fc} it is easy to see that when $\epsilon_i$ and $a_i/\epsilon_i$ are taken to be sufficiently small, $\omega_g + \partial \bar{\partial} R$ corresponds to a \k metric that is toroidal in a neighborhood of $S$.\\

For example, to see that terms in $\partial \bar{\partial} R$ involving $dw$ are small, one can look into the proof of lemma \ref{fc} and imagine that $z$ also depends on $w$. The terms involving $dw$ come from the $\frac{\partial z}{\partial w}dw$ component of $dz$. It is easy to see that $\frac{\partial z}{\partial w},\frac{\partial^2 z}{\partial w^2} = O(|z|)$. Then by the same reason that $\partial \bar{\partial} R$ is bounded in the proof of lemma \ref{fc}, terms in $\partial \bar{\partial} R$ involving $dw$ are of order $O(|z|)$. Recall that $\partial \bar{\partial} R$ is supported in $|z|\leq \max_i (a_ie^{|\lambda_i-1|},a_i/\epsilon_i)$. Therefore the terms in $\partial \bar{\partial} R$ involving $dw$ can be made arbitrarily small, if $\epsilon_i$ and $a_i/\epsilon_i$ are taken to be sufficiently small.
\begin{flushright} $\Box$ \end{flushright}
Using corollary \ref{ff} instead of lemma \ref{fc}, we have the following improvement of lemma \ref{fb} that we do not really need in this paper.

\begin{lm}
\label{fg}
The toroidal metric $\omega_g$ constructed in lemma \ref{fb} can be made flat in a small neighborhood of the origin in each fiber of $\pi$.\\
\end{lm}
{\bf Proof:} The only change necessary from the proof of lemma \ref{fb} is to replace $R$ by $\tilde{R}$ satisfying

\[
\tilde{R}|_{\pi^{-1}(w)} = \sum_{i=1}^n \rho_\epsilon\left(|\hat{z}^i|/a\right)\left(f_{\lambda_i}(|z^i|^2)- |z^i|^2\right)- \rho\left(\frac{2|z|}{R_2}\right)\hat{h},
\]

where $\omega_g|_{\pi^{-1}(w)}=\partial \bar{\partial} (h+\hat{h})|_{\pi^{-1}(w)}$ on $\pi^{-1}(w)$, $h$ is the quadratic term and $\hat{h}$ is higher order term. The key point here is how to make $\hat{h}$ vary smoothly according to $w$. One way to do this is to choose $(h+\hat{h})|_{\pi^{-1}(w)}$ to be the canonical \k potential of $\omega_g|_{\pi^{-1}(w)}$ with respect to the origin $\pi^{-1}(w)\cap S$ as discussed in \cite{berg}.
\begin{flushright} $\Box$ \end{flushright}
By lemma \ref{fb}, we can easily see that
\begin{theorem}
\label{fa}
When $D$ is of normal crossing, any \k metric $\omega_g$ can be perturbed locally near ${\rm Sing}(D)$ to become a global toroidal metric for $(X,D)$.
\end{theorem}
{\bf Proof:} Start with the strata $S$ with the lowest dimension in $D$. Then $S_0 = S$. By lemma \ref{fb}, we can make $\omega_g$ toroidal near $S$. By induction on the dimension of the strata $S\subset D$, and use lemma \ref{fb} repeatedly in each step, we can extend the construction to whole $D$.
\begin{flushright} $\Box$ \end{flushright}
By applying lemma \ref{fg} instead of lemma \ref{fb}, we have the following improvement of theorem \ref{fa} that we do not really need in this paper, but we will need in the case of Calabi-Yau complete intersections.

\begin{theorem}
\label{fh}
The toroidal metric $\omega_g$ constructed in theorem \ref{fa} can be made flat in a small neighborhood of the origin in each fiber of $\pi_I$ for all index set $I$. Moreover, by possibly shinking tubular neighborhoods $\{U_I\}$, $\{\pi_I\}$ can be made compatible in the sense of proposition \ref{fe}.
\end{theorem}
\begin{flushright} $\Box$ \end{flushright}

\se{The construction of Lagrangian torus fibration}
In this section, we will formulate a general theorem on construction of Lagrangian torus fibration via gradient flow. Then we will apply it to our special case of Fermat type quintic family. Assume that we have a family of hypersurfaces $\{Y_t\}$ in an ambient compact \k manifold $(M, \omega_g)$. Assume that $Y_t$ is smooth for $t\not=0$ and $D=Y_0$ is a divisor in $M$ with only normal crossing singularities. We will also assume that

\[
Y_{\rm inv} = Y_t \cap Y_0,\ \ {\rm in\ particular}\ \ C = Y_t \cap {\rm Sing}(Y_0)
\]

is independent of $t$. We will follow the notation in section 7 on normal crossing divisor $D$. In particular $Y_0 = D$ has the stratification

\[
Y_0 = D = \bigcup_{i=0}^{n-1} D^{(i)}_0.
\]

Let $B$ be a smooth real manifold with the stratification

\[
B = \bigcup_{i=0}^{n-1} B^{(i)}_0.
\]

\begin{de}
A map $\pi_0: Y_0 \rightarrow B$ is called a (topologically) smooth Lagrangian torus fibration of $Y_0$ if for all $i$, $\pi_0(D^{(i)}_0) = B^{(i)}_0$, and $\pi_0|_{D^{(i)}_0}: D^{(i)}_0 \rightarrow B^{(i)}_0$ is a Lagrangian torus fibration with each fibre being real $i$-torus.
\end{de}

\begin{theorem}
\label{ha}
Start with a (topologically) smooth Lagrangian torus fibration $\pi_0: Y_0 \rightarrow B$, we can construct a symplectic morphism $F_t: Y_t \rightarrow Y_0$ such that $\pi_t = \pi_0\circ F_t: Y_t \rightarrow B$ is a Lagrangian torus fibration. $\Gamma = \pi_0(C)$ is the singular locus of $\pi_t$. For $b\not\in \Gamma$, $\pi_t^{-1}(b)$ is a real $(n-1)$-torus. For $b \in \Gamma$, $\pi_t^{-1}(b)$ is singular. For $b \in \Gamma \cap B^{(i)}_0$, $\pi_t^{-1}(b)\cap C = \pi_0^{-1}(b)\cap C$ and $F_t: \pi_t^{-1}(b)\backslash C \rightarrow \pi_0^{-1}(b)\backslash C$ is a topologically smooth $(n-i-1)$-torus fibration.
\end{theorem}

{\bf Proof:} According to theorem \ref{fa}, the \k metric $g$ can be perturbed locally near ${\rm Sing}(D)$ to become a global toroidal \k metric $\hat{g}$ on $M$ that is toroidal along ${\rm Sing}(D)$ with respect to $D\cup Y_t$. (Since $Y_t \cap Y_0$ is independent of $t$, it is easy to observe that being toroidal along ${\rm Sing}(D)$ with respect to $D\cup Y_t$ is equivalent to being toroidal along ${\rm Sing}(D)$ with respect to $D\cup Y_{t'}$ for $t,t'\not= 0$.)\\

According to theorem \ref{eg}, there exists a $C^{0,1}$ symplectomorphism $H: (M,\omega_g) \rightarrow (M,\omega_{\hat{g}})$ such that $H(D\cup Y_t) = D\cup Y_t$ for a fixed $t$. Define $\hat{\pi}_0 = H\circ \pi_0 \circ H^{-1}: Y_0 \rightarrow B$. Then $\hat{\pi}_0$ defines a topologically smooth Lagrangian torus fibration for $Y_0$ with respect to the toroidal \k form $\omega_{\hat{g}}$.\\

According to theorem \ref{de}, the inverse gradient flow will induce a symplectic morphism $\hat{F}_t: Y_t \rightarrow Y_0$ with respect to the toroidal \k form $\omega_{\hat{g}}$. $\hat{F}_t$ fixes $Y_{\rm inv}$, and for $x\in D_0^{(i)}\backslash Y_{\rm inv}$, $\hat{F}_t^{-1}(x)$ is a real $(n-i-1)$-torus. Let $F_t = H^{-1}\circ \hat{F}_t \circ H$. Then $F_t: Y_t \rightarrow Y_0$ is a symplectic morphism with respect to the \k form $\omega_g$. Clearly, $\pi_t = \pi_0\circ F_t: Y_t \rightarrow B$ is a Lagrangian torus fibration, and $F_t$, $\pi_t$ so constructed satisfy all conditions in the theorem.
\begin{flushright} $\Box$ \end{flushright}
We can apply this theorem to the situation of the fermat type quintic family in $\mathbb{CP}^4$ with the Fubini-Study metric. We may let $s=1/\psi$. Then $Y_0 = X_\infty$. The natural moment map of the Fubini-Study metric naturally define the topologically smooth Lagrangian torus fibration $\pi_0: X_\infty \rightarrow \partial \Delta$. Let $C = X_\psi \cap {\rm Sing}(X_\infty)$, $\tilde{\Gamma} = \pi_0(C)$. $\tilde{\Gamma}$ is a union of 10 curved triangles. $\tilde{\Gamma}^2$ denotes the union of the interior of all these triangles, $\tilde{\Gamma}^1$ denotes the union of the interior of the edges of all these triangles and $\tilde{\Gamma}^0$ denotes the union of the vertices of all these triangles. Then apply theorem \ref{ha}, we have the following, which is the theorem 3.1 in \cite{lag1}.

\begin{theorem}
\label{hb}
The flow of $V$ will produce a \l fibration $F: X_{\psi} \rightarrow \partial\Delta$. There are 4 types of fibers.\\
(i). For $p\in \partial\Delta\backslash \tilde{\Gamma}$, $F^{-1}(p)$ is a smooth \l 3-torus.\\
(ii). For $p\in \tilde{\Gamma}^2$, $F^{-1}(p)$ is a \l 3-torus with $50$ circles collapsed to $50$ singular points.\\
(iii). For $p\in \tilde{\Gamma}^1$, $F^{-1}(p)$ is a \l 3-torus with $25$ circles collapsed to $25$ singular points.\\
(iv). For $p\in \tilde{\Gamma}^0$, $F^{-1}(p)$ is a \l 3-torus with $5$ 2-torus collapsed to $5$ singular points.
\end{theorem}
\begin{flushright} $\Box$ \end{flushright}
{\bf Remark:} Compare with the constructions with codimension 2 singular locus in the next section, the construction in theorem \ref{hb} is more natural and technically much easier. (The Lagrangian torus fibration for $X_\infty$ basically comes for free.) More importantly, according to recent work of Joyce \cite{Joyce}, which we suspected all along, the actual special Lagrangian fibration for Calabi-Yau 3-fold probably should have singular locus of codimension 1. We believe the Lagrangian torus fibration structure constructed for Fermat type quintic Calabi-Yau in theorem \ref{hb} should be the correct symplectic topological model for the actual special Lagrangian torus fibration.

\se{Deforming to codimension 2 singular locus}
In the previous paper \cite{lag1} we constructed a Lagrangian torus fibration with codimension 1 singular locus as fattening of a graph and discussed possible structure of the related Lagrangian torus fibration with graph singular locus (``expected special Lagrangian fibration structure") based on monodromy information. Following the discussion in Section 1 (Background), we will refer to these two kinds of Lagrangian fibrations as physical model and mathematical model respectively for the sake of distinguishing them. The two models are different. For the mathematical model, the singular locus in $\partial\Delta$ is supposed to be a one-dimensional graph $\Gamma$ and singular fibres have singularity of dimension one. For the physical model, the singular locus $\tilde{\Gamma}$ is a two-dimensional object with boundary that can be viewed as some fattened version of $\Gamma$, and singular fibres have only isolated point singularities. The two models are closely related. The total singular set of the two fibrations in $X_{\psi}$ are both 10 genus six Riemann surfaces and the monodromy of the two torus fibrations are the same. We intend to use tools developed in this section to modify our Lagrangian torus fibration with codimension 1 singular locus (the physical model) to get a Lagrangian torus fibration with graph singular locus (the mathematical model) that coincides with our proposed ``special Lagrangian torus fibration" differential topologically.\\

\subsection{The piecewise smooth argument}
Consider $ {\mathbb{CP}^2}$ with the Fubini-Study metric and the curve $S_0: z_0^5 + z_1^5 + z_2^5 =0$ in ${\mathbb{CP}^2}$. We have the torus fibration $F: {\mathbb{CP}^2}\rightarrow \mathbb{R}^+\mathbb{P}^2$ defined as $F( [z_1,z_2,z_3]) = [|z_1|,|z_2|,|z_3|]$. Choose inhomogenuous coordinate $x_i= z_i/z_0$. Then locally we have $F: {\mathbb{C}^2}\rightarrow {\bf (R^+)^2}$, $F(x_1,x_2)=(r_1,r_2)$, where $x_k=r_ke^{i\theta_k}$. The image of $S_0: x_1^5 + x_2^5 +1=0$ under $F$ is\\
\[
\tilde{\Gamma} =\{(r_1,r_2)|r^5_1 + r^5_2 \geq 1, r^5_1 \leq r^5_2 + 1, r^5_2 \leq r^5_1 + 1 \}.
\]
$S_0$ is a symplectic submanifold. We want to deform $S_0$ symplectically to $S_1$ whose image under $F$ is expected to be\\
\[
\Gamma = \{(r_1,r_2)| 0\leq r_2\leq r_1=1 \ {\rm or} \ 0\leq r_1\leq r_2 =1 \ {\rm or} \ r_1=r_2\geq 1\}.
\]
This is not hard, since $\tilde{\Gamma}$ and $\Gamma$ are very close. A moment of thought suggests the following: When $|x_1|\geq |x_2|\geq 1,$\\
\[
S_t =\left\{{\cal F}_t(x) = \left.\left(\left(\frac{r_2}{r_1}\right)^tx_1,x_2\right)\right|x_1^5 + x_2^5 +1=0\right\};
\]
when $|x_2|\geq |x_1|\geq 1,$\\
\[
S_t =\left\{{\cal F}_t(x) = \left.\left(x_1,\left(\frac{r_1}{r_2}\right)^tx_2\right)\right|x_1^5 + x_2^5 +1=0\right\};
\]
when $|x_2|\geq 1\geq |x_1|,$\\
\[
S_t =\left\{{\cal F}_t(x) = \left.\left(x_1,\left(\frac{1}{r_2}\right)^tx_2\right)\right|x_1^5 + x_2^5 +1=0\right\};
\]
when $1\geq |x_2| \geq |x_1|,$\\
\[
S_t =\left\{{\cal F}_t(x) = \left.\left(\left(\frac{1}{r_2}\right)^tx_1, \left(\frac{1}{r_2}\right)^tx_2\right)\right|x_1^5 + x_2^5 +1=0\right\};
\]
when $|x_1|\geq 1\geq |x_2|,$\\
\[
S_t =\left\{{\cal F}_t(x) = \left.\left(\left(\frac{1}{r_1}\right)^tx_1,x_2\right)\right|x_1^5 + x_2^5 +1=0\right\};
\]
when $1\geq |x_1| \geq |x_2|,$\\
\[
S_t =\left\{{\cal F}_t(x) = \left.\left(\left(\frac{1}{r_1}\right)^tx_1, \left(\frac{1}{r_1}\right)^tx_2\right)\right|x_1^5 + x_2^5 +1=0\right\}.
\]
It is easy to verify that these definitions coincide on the common boundaries. In particular, $F(S_0) = \Gamma$ is a 1-dimensional graph. $S_t$ can also be defined uniformly as

\[
S_t = \left\{{\cal F}_t(x) = \left.\left(\left(\frac{\max(1,r_2)}{\max(r_1,r_2)}\right)^tx_1, \left(\frac{\max(1,r_1)}{\max(r_1,r_2)}\right)^tx_2\right)\right|x_1^5 + x_2^5 +1=0\right\}.
\]

The K\"{a}hler form of the Fubini-Study metric can be written as\\
\[
\omega = \frac{ (1+|x_2|^2)dx_1\wedge d\bar{x}_1 + (1+|x_1|^2) dx_2\wedge d\bar{x}_2 - \bar{x}_1x_2 dx_1\wedge d\bar{x}_2 - \bar{x}_2x_1 dx_2\wedge d\bar{x}_1}{(1+|x|^2)^2}.
\]\\
The K\"{a}hler form of the Fubini-Study metric can also be written as\\
\[
\omega = \frac{ dx_1\wedge d\bar{x}_1 + dx_2\wedge d\bar{x}_2 + (x_2 dx_1 - x_1 dx_2)\wedge (\bar{x}_2d\bar{x}_1-\bar{x}_1d\bar{x}_2)}{(1+|x|^2)^2}.
\]\\
Due to its symmetric nature, to verify that $S_t$ is symplectic, we only need to varify for one region out of six. Consider $1\geq |x_2|\geq |x_1|,$\\
\[
S_t =\left\{{\cal F}_t(x) = \left.\left( \left(\frac{1}{r_2}\right)^tx_1,\left(\frac{1}{r_2}\right)^tx_2\right)\right|x_1^5 + x_2^5 +1=0\right\}.
\]
$x_1^5 + x_2^5 +1=0$ implies that\\
\[
dx_1 = -\left(\frac{x_2}{x_1}\right)^4 dx_2.
\]
Recall that\\
\[
\frac{dr_k}{r_k}= {\rm Re}\left(\frac{dx_k}{x_k}\right),\ d\theta_k = {\rm Im}\left(\frac{dx_k}{x_k}\right).
\]
We have\\
\[
d\left(\left(\frac{1}{r_2}\right)^tx_1\right) = \left(\frac{1}{r_2}\right)^t\left(dx_1 - tx_1\frac{dr_2}{r_2}\right).
\]
\[
\frac{dr_2}{r_2}= {\rm Re}\left(\frac{dx_2}{x_2}\right)= -{\rm Re}\left(\left(\frac{x_1}{x_2}\right)^5\frac{dx_1}{x_1}\right).
\]
\[
d\left(\left(\frac{1}{r_2}\right)^tx_1\right)\wedge d\left(\left(\frac{1}{r_2}\right)^t\bar{x}_1\right) = \left(\frac{1}{r_2}\right)^{2t}\left(dx_1\wedge d\bar{x}_1 +t(x_1d\bar{x}_1 - \bar{x}_1dx_1)\wedge \frac{dr_2}{r_2}\right)
\]
\[
=\left(\frac{1}{r_2}\right)^{2t} \left(1 + t{\rm Re}\left(\left(\frac{x_1}{x_2}\right)^5\right)\right)dx_1\wedge d\bar{x}_1.
\]
\[
d\left(\left(\frac{1}{r_2}\right)^tx_2\right) = \left(\frac{1}{r_2}\right)^t\left(dx_2 - tx_2\frac{dr_2}{r_2}\right).
\]
\[
\frac{dr_2}{r_2}= {\rm Re}\left(\frac{dx_2}{x_2}\right).
\]
\[
d\left(\left(\frac{1}{r_2}\right)^tx_2\right)\wedge d\left(\left(\frac{1}{r_2}\right)^t\bar{x}_2\right) = \left(\frac{1}{r_2}\right)^{2t}\left(dx_2\wedge d\bar{x}_2 +t(x_2d\bar{x}_2 - \bar{x}_2dx_2)\wedge \frac{dr_2}{r_2}\right)
\]
\[
=\left(\frac{1}{r_1}\right)^{2t} (1 - t)dx_2\wedge d\bar{x}_2,
\]
\[
\left( \left(\frac{1}{r_2}\right)^tx_2\right)d\left( \left(\frac{1}{r_2}\right)^tx_1\right) - \left( \left(\frac{1}{r_2}\right)^tx_1\right)d\left( \left(\frac{1}{r_2}\right)^tx_2\right)
= \left(\frac{1}{r_2}\right)^{2t}(x_2dx_1 - x_1dx_2).
\]
\[
= -\left(\frac{1}{r_2}\right)^{2t}x_2\left(1+ \left(\frac{x_1}{x_2}\right)^5\right)dx_1 = \left(\frac{1}{r_2}\right)^{2t} \left(\frac{1}{x_2^4}\right)dx_1
\]\\
By restriction to $S_t$ we get

\[
\frac{\omega|_{S_t}}{dx_1\wedge d\bar{x}_1} = \frac{(1-t)\left(\frac{1}{r_2}\right)^{2t}\left(\frac{r_1}{r_2}\right)^8 + \left(\frac{1}{r_2}\right)^{2t} \left(1 + t{\rm Re}\left(\left(\frac{x_1}{x_2}\right)^5\right)\right) + \left(\frac{1}{r_2}\right)^{4t+8}}{\left(1+r_2^{2-2t} + \left(\frac{r_1}{r_2}\right)^{2t} r_1^{2-2t}\right)^2}
\geq \frac{1}{9}.
\]

These computations show that $S_t$ is symplectic in the region $r_1 < r_2 < 1$. By symmetry, we can see that $S_t$ is symplectic in the other five regions. Namely

\begin{lm}
\label{gd}
$S_t$ is symplectic when $(r_1-1)(r_2-1)(r_1-r_2) \not= 0$.
\end{lm}
\begin{flushright} $\Box$ \end{flushright}
\begin{theorem}
\label{ge}
There exists a family of piecewise smooth Lipschitz continuous Hamiltonian diffeomorphism $h_t: {\mathbb{CP}^2}\rightarrow{\mathbb{CP}^2}$ such that $h_t(S_0)=S_t$ and $h_t$ is identity away from an arbitrarily small neighborhood of $\bigcup_tS_t$. In particular $h_t$ leaves the three coordinate ${\mathbb{CP}^1}$'s invariant.\\
\end{theorem}
{\bf Proof:} Lemma \ref{gd} implies that $S_t$'s are piecewise smooth symplectic submanifolds in ${\mathbb{CP}^2}$. Each $S_t$ is a union of 6 pieces of smooth symplectic submanifolds with boundaries and corners. The 6 pieces have equal area (equal to one-sixth of the total area of $S_t$), which is independent of $t$. $S_0$ is symplectic isotopic to $S_1$ via the family $\{S_t\}$. By corrollary \ref{er}, we may construct a piecewise smooth Lipschitz Hamiltonian diffeomorphism $h_t: {\mathbb{CP}^2}\rightarrow{\mathbb{CP}^2}$ such that $h_t(S_0)=S_t$. According to the notation in corrollary \ref{er}, the part of $\partial S_0$ that is fixed by the original symplectic isotopy flow is $(\partial S_0)_0 = F^{-1}(\Gamma)\cap S_0$. Corrollary \ref{er} asserts that $h_t$ can be made identity on

\[
U = \{p\in {\mathbb{CP}^2}|\min_t({\rm Dist}(p,S_t)) \geq \min(\epsilon_1, \epsilon_2 {\rm Dist}(p,(\partial S_0)_0))\}.
\]

By suitably adjust $\epsilon_1, \epsilon_2$, $U$ will contain all the three coordinate ${\mathbb{CP}^1}$'s. Therefore $h_t$ can be made to leave the three coordinate ${\mathbb{CP}^1}$'s invariant as desired.
\begin{flushright} $\Box$ \end{flushright}

\subsection{The smooth argument}
Notice that $S_t$ is not smooth in the common boundaries $(r_1-1)(r_2-1)(r_1-r_2) = 0$ of the six regions. To modify the definition of $S_t$ to make it smooth, consider real function $b(a)\geq 0$ such that $b(a)+b(-a)=1$ for all $a$ and $b(a)=0$ for $a\leq -\epsilon$. Then consequently, $b(a)=1$ for $a\geq \epsilon$ and $b(a)\leq 1$.\\\\
In general, we may modify the definition of $S_t$ to consider\\
\[
S_t = \left\{{\cal F}_t(x) = \left.\left( \left(\frac{\rho_2}{\rho_0}\right)^tx_1,\left(\frac{\rho_1}{\rho_0}\right)^tx_2\right)\right|x_1^5 + x_2^5 +1=0\right\},
\]
where\\
\[
\rho_1 = r_1^{b(\log r_1)},\ \rho_2 = r_2^{b(\log r_2)},\ \rho_0 = r_1\left(\frac{r_2}{r_1}\right)^{b(\log (r_2/r_1))}.
\]\\
$S_t$ is now smooth and is only modified in an $\epsilon$-neighborhood of $(r_1-1)(r_2-1)(r_1-r_2) = 0$. Let $\tilde{\Gamma}_\epsilon = F(S_1)$, then $\tilde{\Gamma}_\epsilon$ coincides with the graph $\Gamma$ outside of the $\epsilon$-neighborhood of the vertex $r_1=r_2=1$ of $\Gamma$. In the $\epsilon$-neighborhood of the vertex $r_1=r_2=1$ of $\Gamma$, $\tilde{\Gamma}_\epsilon$ is an $\epsilon$ fattening of $\Gamma$. This is the price to pay if we want $S_1$ to be smooth.\\\\
$x_1^5 + x_2^5 +1=0$ implies that\\
\[
dx_1 = -\left(\frac{x_2}{x_1}\right)^4 dx_2.
\]
Recall that
\[
\frac{dr_k}{r_k}= {\rm Re}\left(\frac{dx_k}{x_k}\right),\ d\theta_k = {\rm Im}\left(\frac{dx_k}{x_k}\right).
\]
Consequently\\
\[
\frac{dr_1}{r_1}= {\rm Re}\left(\frac{dx_1}{x_1}\right)= -{\rm Re}\left(\left(\frac{x_2}{x_1}\right)^5\frac{dx_2}{x_2}\right).
\]
\[
\frac{dr_2}{r_2}= {\rm Re}\left(\frac{dx_2}{x_2}\right)= -{\rm Re}\left(\left(\frac{x_1}{x_2}\right)^5\frac{dx_1}{x_1}\right).
\]\\
Assume $\lambda(a) = b(\log a) + b'(\log a)\log a$, $\lambda_0 = \lambda(r_2/r_1)$, $\lambda_1 = \lambda(r_1)$, $\lambda_2 = \lambda(r_2)$. Then\\
\[
\frac{d\rho_1}{\rho_1}= (b + b'\log (r_1))\frac{dr_1}{r_1} = \lambda_1\frac{dr_1}{r_1},\ \frac{d\rho_2}{\rho_2}= (b + b'\log (r_2))\frac{dr_2}{r_2} = \lambda_2\frac{dr_2}{r_2}.
\]\\
\[
\frac{d\rho_0}{\rho_0}= \frac{dr_1}{r_1} + (b + b'\log (r_2/r_1))\left(\frac{dr_2}{r_2} - \frac{dr_1}{r_1}\right)= \frac{dr_1}{r_1} + \lambda_0\left(\frac{dr_2}{r_2} - \frac{dr_1}{r_1}\right).
\]\\
\[
d\left(\left(\frac{\rho_2}{\rho_0}\right)^tx_1\right) = \left(\frac{\rho_2}{\rho_0}\right)^t\left(dx_1 + tx_1\left(\frac{d\rho_2}{\rho_2} - \frac{d\rho_0}{\rho_0}\right)\right).
\]
\[
\frac{d\rho_2}{\rho_2} - \frac{d\rho_0}{\rho_0} = -(1 - \lambda_0) \frac{dr_1}{r_1} - (\lambda_0 - \lambda_2)\frac{dr_2}{r_2}.
\]
\[
d\left(\left(\frac{\rho_2}{\rho_0}\right)^tx_1\right)\wedge d\left(\left(\frac{\rho_2}{\rho_0}\right)^t\bar{x}_1\right) = \left(\frac{\rho_2}{\rho_0}\right)^{2t}\left(dx_1\wedge d\bar{x}_1 -t(x_1d\bar{x}_1 - \bar{x}_1dx_1)\wedge \left(\frac{d\rho_2}{\rho_2} - \frac{d\rho_0}{\rho_0}\right)\right)
\]
\[
=\left(\frac{\rho_2}{\rho_0}\right)^{2t} \left(1 - (1-\lambda_0)t + (\lambda_0 - \lambda_2)t{\rm Re}\left(\left(\frac{x_1}{x_2}\right)^5\right)\right)dx_1\wedge d\bar{x}_1.
\]\\
\[
d\left(\left(\frac{\rho_1}{\rho_0}\right)^tx_2\right) = \left(\frac{\rho_1}{\rho_0}\right)^t\left(dx_2 + tx_2\left(\frac{d\rho_1}{\rho_1} - \frac{d\rho_0}{\rho_0}\right)\right).
\]
\[
\frac{d\rho_1}{\rho_1} - \frac{d\rho_0}{\rho_0} = -(1 - \lambda_0 - \lambda_1) \frac{dr_1}{r_1} - \lambda_0 \frac{dr_2}{r_2}.
\]
\[
d\left(\left(\frac{\rho_1}{\rho_0}\right)^tx_2\right)\wedge d\left(\left(\frac{\rho_1}{\rho_0}\right)^t\bar{x}_2\right) = \left(\frac{\rho_1}{\rho_0}\right)^{2t}\left(dx_2\wedge d\bar{x}_2 +t(x_2d\bar{x}_2 - \bar{x}_2dx_2)\wedge \left(\frac{d\rho_1}{\rho_1} - \frac{d\rho_0}{\rho_0}\right)\right)
\]
\[
=\left(\frac{\rho_1}{\rho_0}\right)^{2t} \left(1+ (1-\lambda_0 - \lambda_1)t{\rm Re}\left(\left(\frac{x_2}{x_1}\right)^5\right) - \lambda_0 t\right)dx_2\wedge d\bar{x}_2.
\]\\
\[
\alpha = \left( \left(\frac{\rho_1}{\rho_0}\right)^tx_2\right)d\left( \left(\frac{\rho_2}{\rho_0}\right)^tx_1\right) - \left( \left(\frac{\rho_2}{\rho_0}\right)^tx_1\right)d\left( \left(\frac{\rho_1}{\rho_0}\right)^tx_2\right)
\]
\[
= \left(\frac{\rho_1\rho_2}{\rho_0^2}\right)^{t}\left(x_2dx_1 - x_1dx_2 + tx_1x_2\left(\frac{d\rho_2}{\rho_2} - \frac{d\rho_1}{\rho_1}\right)\right)
\]
\[
= \left(\frac{\rho_1\rho_2}{\rho_0^2}\right)^{t} \left(-\left(\frac{1}{x_2^4}\right)dx_1 + tx_1x_2\left(\lambda_2\frac{dr_2}{r_2} - \lambda_1\frac{dr_1}{r_1}\right)\right).
\]\\
\[
\alpha \wedge \bar{\alpha} = \left(\frac{\rho_1\rho_2}{\rho_0^2}\right)^{2t}\left(\left(\frac{1}{r_2}\right)^8 - t\left(\frac{\bar{x}_1\bar{x}_2}{x_2^4}dx_1 - \frac{x_1x_2}{\bar{x}_2^4}d\bar{x}_1\right)\left(\lambda_2\frac{dr_2}{r_2} - \lambda_1\frac{dr_1}{r_1}\right)\right)
\]
\[
= \left(\frac{\rho_1\rho_2}{\rho_0^2}\right)^{2t}\left(\left(\frac{1}{r_2}\right)^8 + t\frac{1}{r_2^8}\left(\lambda_2{\rm Re}(x_1^5) + \lambda_1{\rm Re}(x_2^5)\right)\right)dx_1d\bar{x}_1.
\]
By restricting to $S_t$ we get\\
\[
\frac{\omega|_{S_t}}{dx_1\wedge d\bar{x}_1}= \left[\left(\frac{\rho_1}{\rho_0}\right)^{2t}\left(1+ (1-\lambda_0 - \lambda_1)t{\rm Re}\left(\left(\frac{x_2}{x_1}\right)^5\right) - \lambda_0 t\right)\left(\frac{r_1}{r_2}\right)^8\right.
\]
\[
+ \left(\frac{\rho_2}{\rho_0}\right)^{2t}\left(1 - (1-\lambda_0)t + (\lambda_0 - \lambda_2)t{\rm Re}\left(\left(\frac{x_1}{x_2}\right)^5\right)\right)
\]
\[
+ \left.\left.\left(\frac{\rho_1\rho_2}{\rho_0^2}\right)^{2t}\left(\left(\frac{1}{r_2}\right)^8 + t\frac{1}{r_2^8}\left(\lambda_2{\rm Re}(x_1^5) + \lambda_1{\rm Re}(x_2^5)\right)\right)\right]\right/\left(1+\left(\frac{\rho_1}{\rho_0}\right)^{2t}r_2^2 + \left(\frac{\rho_2}{\rho_0}\right)^{2t}r_1^2\right)^2.
\]\\
\[
= \left[\left(\frac{\rho_1}{\rho_0}\right)^{2t}\left(1+ t{\rm Re}\left(\left(\frac{x_2}{x_1}\right)^5\right)\right)\left(\frac{r_1}{r_2}\right)^8 + \left(\frac{\rho_2}{\rho_0}\right)^{2t}(1 - t) + \left(\frac{\rho_1\rho_2}{\rho_0^2}\right)^{2t}\left(\frac{1}{r_2}\right)^8 \right.
\]
\[
+\lambda_0t\left(\left(\frac{\rho_2}{\rho_0}\right)^{2t}\left(1 + {\rm Re}\left(\left(\frac{x_1}{x_2}\right)^5\right)\right) - \left(\frac{\rho_1}{\rho_0}\right)^{2t}\left(1+ {\rm Re}\left(\left(\frac{x_2}{x_1}\right)^5\right)\right)\left(\frac{r_1}{r_2}\right)^8\right)
\]
\[
+\lambda_1t\left(\frac{\rho_1}{\rho_0}\right)^{2t}\left(\rho_2^{2t} {\rm Re}(x_2^5)\left(\frac{1}{r_2}\right)^8 - {\rm Re}\left(\left(\frac{x_2}{x_1}\right)^5\right)\left(\frac{r_1}{r_2}\right)^8\right)
\]
\[
+ \left.\left.\lambda_2t\left(\frac{\rho_2}{\rho_0}\right)^{2t}\left(\rho_1^{2t} {\rm Re}(x_1^5)\left(\frac{1}{r_2}\right)^8 - {\rm Re}\left(\left(\frac{x_1}{x_2}\right)^5\right)\right)\right]\right/\left(1+\left(\frac{\rho_1}{\rho_0}\right)^{2t}r_2^2 + \left(\frac{\rho_2}{\rho_0}\right)^{2t}r_1^2\right)^2.
\]\\
{\bf Remark:}
With the help of symmetry, the cases remain to be varified are $\epsilon$-neighborhood of $\{r_1=r_2\leq 1-\epsilon\}$, $\{r_2=1,0\leq r_1\leq 1-\epsilon\}$ and $\{r_1=r_2=1\}$. Observe that in these regions the coefficients of $\lambda_i$ are of order $\epsilon$. Since $|\lambda_i|$ are bounded, it is not hard to see that the $\epsilon$-variation of $\frac{\omega|_{S_t}}{dx_1\wedge d\bar{x}_1}$ when $r_1,r_2$ are bounded is of order $\epsilon$. Therefore we can reduce our task to verifying that $S_t$ is symplectic on $\{r_1=r_2\leq 1-\epsilon\}$, $\{r_2=1,0\leq r_1\leq 1-\epsilon\}$ and $\{r_1=r_2=1\}$.\\\\
When $r_1=r_2\leq 1-\epsilon$, $\rho_1 = \rho_2 =1$, $\lambda_1=\lambda_2=0$, $\lambda_0=\frac{1}{2}$.\\
\[
\frac{\omega|_{S_t}}{dx_1\wedge d\bar{x}_1}= \frac{\left(\frac{1}{\rho_0}\right)^{2t}\left(2+ t{\rm Re}\left(\left(\frac{x_2}{x_1}\right)^5\right) - t\right) + \left(\frac{1}{\rho_0}\right)^{4t}\left(\left(\frac{1}{r_2}\right)^8 \right)}{\left(1+\left(\frac{1}{\rho_0}\right)^{2t}(r_2^2 + r_1^2)\right)^2}.
\]\\
When $r_1=r_2\leq 1-\epsilon$,\\
\[
-\frac{1}{2} \leq {\rm Re}\left(\left(\frac{x_2}{x_1}\right)^5\right) \leq 1,\ \ \rho_0 = \sqrt{r_1r_2}<1.
\]
Therefore\\
\[
\frac{\omega|_{S_t}}{dx_1\wedge d\bar{x}_1} \geq \frac{3-\frac{3}{2}t}{9} \geq \frac{1}{6}.
\]

When $r_1=r_2=1$, $\rho_1 = \rho_2 = \rho_0 =1$, $\lambda_1=\lambda_2=\lambda_0=\frac{1}{2}$.\\\\
$x_1^5 + x_2^5 +1 =0$ implies,\\
\[
{\rm Re}(x_1^5) = {\rm Re}(x_2^5) = {\rm Re}\left(\left(\frac{x_2}{x_1}\right)^5\right) = -\frac{1}{2}.
\]
Therefore\\
\[
\frac{\omega|_{S_t}}{dx_1\wedge d\bar{x}_1}= \frac{1 + t{\rm Re}\left(\left(\frac{x_2}{x_1}\right)^5\right) + (1-t)+1}{9} = \frac{3-\frac{3t}{2}}{9} = \frac{2-t}{6} \geq \frac{1}{6}.
\]

Another boundary we need to check is when $0\leq r_1 \leq 1-\epsilon$ and $r_2=1$. Then $\rho_1 = \rho_2 = \rho_0 =1$, $\lambda_1=0$, $\lambda_0=1$, $\lambda_2=\frac{1}{2}$. We have

\[
\frac{\omega|_{S_t}}{dx_1\wedge d\bar{x}_1}= \frac{1}{\left(2 + r_1^2\right)^2}\left(\left(1+ (1-\lambda_0 - \lambda_1)t{\rm Re}\left(\left(\frac{x_2}{x_1}\right)^5\right) - \lambda_0 t\right)r_1^8\right.
\]
\[
+ \left.\left(1 - (1-\lambda_0)t + (\lambda_0 - \lambda_2)t{\rm Re}\left(\left(\frac{x_1}{x_2}\right)^5\right)\right) + \left(1 + t\left(\lambda_2{\rm Re}(x_1^5) + \lambda_1{\rm Re}(x_2^5)\right)\right)\right).
\]
\[
= \frac{1}{\left(2 + r_1^2\right)^2}\left((1 - t)r_1^8+\left(1 + \frac{1}{2}t{\rm Re}\left(\left(\frac{x_1}{x_2}\right)^5\right)\right) + \left(1 + t\frac{1}{2}{\rm Re}(x_1^5) \right)\right).
\]
\[
= \frac{1}{\left(2 + r_1^2\right)^2}\left((1 - t)r_1^8+ \left(2 + t{\rm Re}(x_1^5) \right)\right) \geq \frac{1}{9}.
\]

Now we have shown

\begin{lm}
\label{gb}
$S_t$ is symplectic for $t\in [0,1]$. Namely, $S_0$ is symplectic isotopic to $S_1$ via the family $S_t$.
\end{lm}
\begin{flushright} $\Box$ \end{flushright}
\begin{co}
\label{gc}
There exists a family of Hamiltonian diffeomorphism $h_t: {\mathbb{CP}^2}\rightarrow{\mathbb{CP}^2}$ such that $h_t(S_0)=S_t$ and $h_t$ is identity away from an arbitrary small neighborhood of $\bigcup_tS_t$. In particular $h$ leaves the three coordinate ${\mathbb{CP}^1}$'s invariant.\\
\end{co}
{\bf Proof:} Lemma \ref{gb} implies that $S_0$ is symplectic isotopic to $S_1$ via the family $S_t$. By theorem \ref{ee} above conclusion is immediate.
\begin{flushright} $\Box$ \end{flushright}

\subsection{Deforming to codimension 2 singular locus}
Recall that\\
\[
X_{\infty} =\bigcup_{\tiny{\begin{array}{c}I \subset \{1,2,3,4,5\}\\ 0<|I|<5\end{array}}} D_I
\]
where\\
\[
D_I = \{z:\ z_i=0,z_j\not= 0,\ {\rm for}\ i\in I,j\in\{1,2,3,4,5\}\backslash I\}
\]\\
is a ($4-|I|$)-dimensional complex torus. There is a natural torus fibration\\
\[
F: X_{\infty} \rightarrow \partial\Delta = \bigcup_{\tiny{\begin{array}{c} I \subset \{1,2,3,4,5\}\\ 0<|I|<5\end{array}}} \Delta_I.
\]\\
The fibres over $\Delta_I$ are ($4-|I|$)-dimensional torus.\\\\
$X_{\infty}$ is a union of 5 ${\mathbb{CP}^3}$'s ($\overline{D_k}$). They intersect in 10 ${\mathbb{CP}^2}$'s ($\overline{D_{ij}}$). In each of these ${\mathbb{CP}^2}$'s ($\overline{D_I}$, $|I|=2$) there is a quintic curve ($\Sigma_{I^c}$), where $I^c$ is the set of compliment of $I \subset \{1,2,3,4,5\}$.\\
\[
\Sigma_{ijk} =\{[z]\in {\mathbb{CP}^4}| z_i^5 + z_j^5 + z_k^5 =0, z_l=0\ \ {\rm for}\ l\in\{1,2,3,4,5\}\backslash\{i,j,k\}\}.
\]\\
$\Sigma_{ijk}$ is a genus 6 curve. Let\\
\[
\Sigma = \bigcup_{\{i,j,k\}\subset \{1,2,3,4,5\}} \Sigma_{ijk}
\]\\
$\Sigma = X_{\infty}\cap X_{\psi}$ for any $\psi$.\\\\
$\Sigma$ is the singular set of our Lagrangian construction. Its image under $F$ is $\tilde{\Gamma}$. We want to modify the Lagrangian fibration so that the image of $\Sigma$ is $\Gamma$.\\

\begin{lm}
\label{gf}
One can construct a topologically smooth Lagrangian fibration $F_\infty: X_\infty \rightarrow \partial \Delta$ with respect to $\omega_{\rm FS}$ such that $F_\infty(\Sigma) = \Gamma$.\\
\end{lm}
{\bf Proof:} Let\\
\[
X^{(k)}_{\infty} =\bigcup_{\tiny{\begin{array}{c}I \subset \{1,2,3,4,5\}\\ 0<|I|\leq k\end{array}}} D_I
\]
Then $X^{(2)}_{\infty}$ is the union of the 10 ${\mathbb{CP}^2}$'s. By corollary \ref{ge}, we can easily construct a $C^{0,1}$ symplectic flow $h_t: X^{(2)}_{\infty} \rightarrow X^{(2)}_{\infty}$, such that $F(h_1(\Sigma)) = \Gamma$. Apply theorem \ref{ef}, we can extend $h_t$ as a $C^{0,1}$ symplectic flow $h_t: X_{\infty} \rightarrow X_{\infty}$. $F_\infty =F\circ h_1$ induces a Lagrangian fibration of $X_{\infty}$ with respect to the Fubini-Study metric such that $F_\infty(\Sigma) = \Gamma$.
\begin{flushright} $\Box$ \end{flushright}
Let $\Gamma = \Gamma^1\cup \Gamma^2\cup \Gamma^3$, where $\Gamma^1$ is the smooth part of $\Gamma$,\\
\[
\Gamma^2 = \bigcup_{i,j}P_{ij},\ \Gamma^3 = \bigcup_{i,j,k}P_{ijk},
\]
then we have\\
\begin{theorem}
\label{ga}
Start with Lagrangian fibration $F_\infty$ the gradient method will produce a Lagrangian fibration $F_{\psi}: X_{\psi} \rightarrow \partial\Delta$. There are 4 types of fibres.\\
(i). For $p\in \partial\Delta\backslash \Gamma$, $F_{\psi}^{-1}(p)$ is a Lagrangian 3-torus.\\
(ii). For $p\in \Gamma^1$, $F_{\psi}^{-1}(p)$ is a type $I_5$ singular fibre.\\
(iii). For $p\in \Gamma^3$, $F_{\psi}^{-1}(p)$ is a type $II_{5\times 5}$ singular fibre.\\
(iv). For $p\in \Gamma^2$, $F_{\psi}^{-1}(p)$ is a type $III_5$ singular fibre.\\
\end{theorem}
{\bf Proof:} According to lemma \ref{gf}, we have a topologically smooth Lagrangian fibration $F_\infty: X_\infty \rightarrow \partial \Delta$ with respect to $\omega_{\rm FS}$ such that $F_\infty(\Sigma) = \Gamma$.\\

According to theorem \ref{ha}, one can construct a symplectic morphism $H_\psi: X_\psi \rightarrow X_\infty$ such that $X_{\rm inv} = X_\infty \cap X_0$ is fixed by $H_\psi$. For $2\leq |I|\leq 4$, the inverse image of each point in $D_I\backslash X_{\rm inv}$ under $H_\psi$ is a torus of dimension $|I|-1$ in $X_\psi\backslash X_{\rm inv}$. When $|I|=1$, $H_\psi$ is 1-1 on $D_I$. Define $F_\psi = F_\infty \circ H_\psi$.\\

Let $\Delta^{(k)}$ denote the $k$-skeleton of $\Delta$. Then $\partial \Delta = \Delta^{(3)}$. For $p\in \partial \Delta \backslash \Delta^{(2)}$, $F_{\infty}^{-1}(p)\subset \displaystyle \bigcup_{|I|=1} D_I$ is a 3-torus that under $H_\psi^{-1}$ is mapped 1-1 to a 3-torus $F_{\psi}^{-1}(p) \subset X_\psi$. For $k\leq 2$ and $p\in \Delta^{(k)} \backslash (\Delta^{(k-1)} \cup \Gamma)$, $F_{\infty}^{-1}(p)\subset \displaystyle \bigcup_{|I|=4-k} (D_I\backslash \Sigma)$ is a k-torus. Each point in $F_{\infty}^{-1}(p)$ under $H_\psi^{-1}$ is mapped to a $(3-k)$-torus. The whole fibre $F_{\infty}^{-1}(p)$ under $H_\psi^{-1}$ is mapped to a 3-torus $F_{\psi}^{-1}(p) \subset X_\psi$. Now we have proved the statement of the theorem in the case (i).\\

For $p\in \Gamma^1$, $F_{\infty}^{-1}(p)\subset \displaystyle \bigcup_{|I|=2} D_I$ is a 2-torus. Each point in $F_{\infty}^{-1}(p)\backslash \Sigma$ under $H_\psi^{-1}$ is mapped to a circle. $F_{\infty}^{-1}(p)\cap \Sigma$, which is a union of 5 circles, will be fixed by $H_\psi$. Then it is very easy to see that the whole fibre $F_{\infty}^{-1}(p)$ under $H_\psi^{-1}$ is mapped to a type $I_5$ singular fibre $F_{\psi}^{-1}(p) \subset X_\psi$, that is a product of a circle and a $I_5$ Kodaira singular fibre.\\

For $p\in \Gamma^2$, $F_{\infty}^{-1}(p)\subset \displaystyle \bigcup_{|I|=1} D_I$ is a circle. Each point in $F_{\infty}^{-1}(p)\backslash \Sigma$ under $H_\psi^{-1}$ is mapped to a 2-torus. $F_{\infty}^{-1}(p)\cap \Sigma$, which is a union of 5 points, will be fixed by $H_\psi$. Then it is very easy to see that the whole fibre $F_{\infty}^{-1}(p)$ under $H_\psi^{-1}$ is mapped to a type $III_5$ singular fibre $F_{\psi}^{-1}(p) \subset X_\psi$.\\

For $p\in \Gamma^3$, $F_{\infty}^{-1}(p)\subset \displaystyle \bigcup_{|I|=2} D_I$ is a 2-torus that is parametrized by $(\theta_1,\theta_2)$. Each point in $F_{\infty}^{-1}(p)\backslash \Sigma$ under $H_\psi^{-1}$ is mapped to a circle. It is not hard to figure out that

\[
F_{\infty}^{-1}(p)\cap \Sigma = \left\{ (\theta_1,\theta_2)\left|
\begin{array}{ll}&e^{5i(\theta_1+\theta_2)}=1,\cos 5\theta_1 \leq -\frac{1}{2}\\{\rm or}&e^{5i(2\theta_1-\theta_2)}=1,\cos 5\theta_1 \leq -\frac{1}{2}\\{\rm or}&e^{5i(2\theta_2-\theta_1)}=1,\cos 5\theta_2 \leq -\frac{1}{2}\end{array}
\right.\right\}
\]

is the graph in the 2-torus $F_{\infty}^{-1}(p)$ corresponding to type $II_{5\times 5}$ singular fibre that will be fixed by $H_\psi$. Then it is very easy to see that the whole fibre $F_{\infty}^{-1}(p)$ under $H_\psi^{-1}$ is mapped to a type $II_{5\times 5}$ singular fibre $F_{\psi}^{-1}(p) \subset X_\psi$.
\begin{flushright} $\Box$ \end{flushright}

{\bf Acknowledgement:} I would like to thank Qin Jing for many very stimulating discussions during the course of my work, and helpful suggestions while carefully reading my early draft. I would also like to thank Prof. S.-T. Yau for his constant encouragement. This work was done while I was in Columbia University. I am very grateful to Columbia University for excellent research environment.\\\\

\end{document}